\magnification=1000
\baselineskip12.7pt
\hsize=13.7truecm
\vsize=22.6truecm
\leftskip=0truecm
\rightskip=0truecm
\hoffset=1,22truecm
\voffset=0truecm
\hfuzz=10pt
\pretolerance=500
\tolerance=1000
\brokenpenalty=5000

\topskip=2truecm

\def\D{\displaystyle}
\def\Ev{\vskip0,5cm}
\def\ev{\vskip0,25cm}
\def\adh{\hbox{adh}}
\def\fr{\hbox{fr}}
\def\inte{\hbox{int}}

\def\min{\hbox{min}}
\def\max{\hbox{max}}

\def\N{\hbox{\t N}}
\def\T{\hbox{T}}

\def\u{\epsilon}
\def\L{\Lambda}
\def\Ll{\Lambda^{\ell oc}}
\def\lg{\longrightarrow}
\def\n{\noindent}
\def\N{\hbox{\bbbb N}}
\def\R{\hbox{\bbbb R}}
\def\T{\hbox{T}}
\def\é{\hbox{\'e}}
\def\ê{\hbox{\^e}}
\def\è{\hbox{\`e}}
\def\à{\hbox{\`a}}
\def\ù{\hbox{\`u}}
\def\K{\hbox{\rm K}}
\def\ti{\widetilde}

\def\picture #1 by #2 (#3){\vbox to #2{\hrule width #1 height
0pt
\vfill\special{picture #3}}}
\def\scaledpicture #1 by #2 (#3 scaled
#4){{\dimen0=#1\dimen1=#2\divide\dimen0 by 1000
\multiply\dimen0 by #4\divide\dimen1 by 1000 \multiply\dimen1
by #4
\picture \dimen0 by \dimen1 (#3 scaled #4)}}

\font\bb=msbm10

\font\bbbb=msbm10 at 7pt
\def\d #1{\hbox{\bb #1}}

\def\r #1{\hbox{$\cal #1 $}}

\font\tm=cmr10
\font\ce=eufm10 at 10pt

\font\Tm=cmr10 at 14pt
\font\ttt=cmr10 at 9pt
\font\td=cmr10 at 9pt
\font\tt=cmr9 at 8pt
\font\t=cmr10 at 7pt
\font\tcq=cmr10 at 5pt

\font\tttsl=cmsl10 at 7pt
\font\Tmind=cmr10 at 8pt

\font\titt=cmb10 at 8pt
\font\LL=lasy10 at 14pt
\def\f{{\hbox{\hskip 5mm \LL 2}}}

\font\acc=cmsl10 at 7pt

\def\lkl{\L_i^{\ell oc}}

\vglue0cm
\rm

\def\date{\the\day\  \ifcase\month\or janvier\or f\'evrier\or
mars\or avril\or mai \or juin \or juillet\or ao\^ut\or
septembre\or octobre\or novembre\or d\'ecembre\fi
\  {\the\year} }

\hbox to 13.5cm{\hrulefill}
\hskip-7mm\vtop{\hsize=13,5truecm  {\tttsl \baselineskip8pt
\noindent
\vfill} }\ \kern0cm\hfill\

\vbox{ \vskip1,5cm}
\centerline{\hskip0mm \hbox{\Tm \'EQUISINGULARIT\'E
R\'EELLE II :  }}
\vskip3,5mm
\centerline{\hskip0mm \hbox{\Tm INVARIANTS LOCAUX ET CONDITIONS
DE R\'EGULARIT\'E
       } }
\vskip10mm
\centerline{\hskip0mm Georges COMTE \&
 Michel MERLE}
 \Ev
 \centerline{\hskip0mm (\ttt Appendice de G. Comte, P. Graftieaux
 \& M. Merle)}
\vskip5mm
 \centerline{ \hbox to 1.5cm{\hrulefill}}
\vskip3mm
\hskip-7mm \vtop{\hsize=13,5truecm \tolerance10000
{ \tt \baselineskip=9pt
R{\tcq \'ESUM\'E}. -- On d\'efinit, pour un germe d'ensemble sous-analytique,
 deux nouvelles suites finies d'invariants num\'eriques.
La premi\`ere a pour termes les localisations des
courbures de Lipschitz-Killing classiques, la seconde est l'\'equivalent r\'eel des
caract\'eristiques \'evanescentes complexes introduites par M.
Kashiwara.
On montre que chaque terme d'une de ces suites est combinaison lin\'eaire des termes de l'autre, puis
on relie ces invariants \`a la g\'eom\'etrie des discriminants des
projections du germe sur des plans de toutes les dimensions.
Il apparait alors que ces invariants sont continus le long de
strates de Verdier d'une stratification sous-analytique d'un ferm\'e.

\vfill} }\ \kern0cm\hfill\
\vskip5mm
 \hskip-7mm \vtop{\hsize=13,5truecm
{\tt \baselineskip=9pt
A{\tcq BSTRACT}. --  For germs of subanalytic sets, we define two
finite sequences of
new numerical invariants.
The first one is  obtained by localizing the classical
Lipschitz-Killing
curvatures, the second one is the real analogue of the
evanescent characteristics introduced by M. Kashiwara.
We show that each invariant of one sequence is a linear combination
of the invariants of the other sequence.
We then connect our invariants to the geometry of the
discriminants of all dimension.
Finally we prove that these invariants are continuous along
Verdier strata of a  closed subanalytic set.
\vfill} }\ \kern0cm\hfill\
\Ev

\vskip1,5cm

\centerline{\bf 0. Introduction}
\Ev

D\'esignons par $X$ un ensemble sous-analytique de  $\d R^n$
et pour $y\in X$, notons $X_y$ le germe de
$X$ en $y$ puis $d$ sa dimension.

Dans [Co2] il est prouv\'e que la fonction densit\'e locale
$y\mapsto \Theta_d(X_y)$
(la localisation du volume) est une fonction continue le
long de strates de Verdier ou m\^eme $(b^*)$-r\'eguli\`eres de $X$ 
(voir [Va1], [Va2]
pour une preuve du caract\`ere lipschitzien de la
densit\'e le long de strates de Verdier et pour la
continuit\'e le long de strates de Whitney). Il s'agissait
dans [Co2] de
ramener l'\'etude de la densit\é \`a celle d'un invariant
$\sigma_d(X_y)$   associ\'e
aux discriminants des
projections du germe $X_y$  sur des plans de dimension
$d$, via une formule du type Cauchy-Crofton localis\'ee
\'egalant $\Theta_d(X_y)$ et $\sigma_d(X_y)$
([Co1], [Co2] Th\'eor\`emes 1.10 et 1.16).
Ce type de r\'esultats reliant le comportement d'un invariant et
la g\'eom\'etrie des discriminants s'inscrit dans le programme
de l'\'equisingularit\'e inaugur\'e par Zariski pour les ensembles
alg\'ebriques complexes (Voir [Za1,2,3,4]).
La condition g\'eom\'etrique portant sur les discrimants
$(d-1)$-dimensionnels qui
assure la continuit\'e de la densit\'e le long d'une strate
$Y$ est que ceux-ci aient le long de $Y$ un c\^one normal dont la
dimension des
fibres soit major\'ee par la dimension g\'en\'erique, ie
$d-\dim(Y)-1$. Il est prouv\'e dans [Co2], Proposition 3.7, que
cette condition est satisfaite le long de strates
de Verdier.
Dans le cadre analytique complexe densit\'e locale et multiplicit\'e
 co\"\i ncident ([Dra]), de m\^eme que conditions de Whitney
et de Verdier ([He-Me1], [Te2]).
Ce faisant, le r\'esultat de [Co2] \'etend \`a la
g\'eom\'etrie r\'eelle celui de [Hi1] selon lequel la condition de
Whitney assure l'\'equimultiplicit\'e.
\ev
Dans le pr\'esent article, nous \'etendons  l'\'etude faite dans
[Co2] \`a  toutes les dimensions de projection :
nous pr\'esentons ici des invariants polaires et de courbures
dont le comportement traduit la g\'eom\'etrie des discriminants
de $X_y$ associ\'es \`a des projections sur des plans de dimension
$i$, pour $0\le i\le d$.
 Pour cela, d'une part
nous  consid\'erons  le volume comme le dernier invariant
de la suite :
$$ \L_*(X)= \big( \L_0(X)=\chi(X), \cdots, \L_d(X)\big)$$
des courbures de
Lipschitz-Killing de $X$. L'\'egalit\'e
 $  \L_d(X)=Vol_d(X)$ \'etant la formule de Cauchy-Crofton classique.
 La localisation en $y$ de $\L_*(X)$ permet
alors de d\'efinir une
suite de nouveaux invariants  attach\'es au germe $X_y$ :
$$\Ll_*(X_y) =\big(\Ll_0(X_y)=1, \dots, \Ll_d(X_y)=\Theta_d(X_y)\big).
$$
Nous appelons ces invariants les {\sl courbures de Lipschitz-Killing
locales}.
 D'autre part la d\'efinition de
l'invariant $\sigma_d(X_y)$ se g\'en\'eralise \`a
toutes les dimensions, ce qui donne lieu \`a la suite des
{\sl invariants polaires} :
 $$ \sigma_*(X_y)=\big(\sigma_0(X_y)=1, \cdots, \sigma_d(X_y)\big).$$
L'\'egalit\'e
 $  \sigma_d(X_y)=\Ll_d(X_y)(=\Theta_d(X_y))$
est la formule de Cauchy-Crofton
locale de [Co2]. Celle-ci
pose la question de la d\'ependance des termes d'une
des deux suites $\Ll_*$ et $\sigma_*$ en
fonction des termes de l'autre. Cette possibilit\'e est
notamment appuy\'ee par le probl\`eme pos\'e dans [Gr-Sc]
(Probl\`eme $14$ ou  [Sc-McMu] $14.3$)  consistant \`a se demander
si le th\'eor\`eme de Hadwiger  poss\`ede un analogue en
g\'eom\'etrie convexe sph\'erique. En effet, si pour un
convexe $K$ de la sph\`ere unit\'e $S^{n-1}$ on note
$\widehat K$ le c\^one de sommet $0$ support\'e par $K$,
les suites $\sigma_*$ et $\Ll_*$
d\'efinissent les valuations $K\mapsto \sigma_*(\widehat K_0)$
et $K\mapsto \Ll_*(\widehat K_0)$ sur les convexes sph\'eriques.
Celles-ci sont continues relativement \`a la m\'etrique de Hausdorff
et invariantes sous les rotations de la sph\`ere.
Dans le cas euclidien il est connu depuis
[Had] (voir aussi [Kl]) que les valuations sur les convexes de $\d R^n$
 continues
 et  invariantes sous les
isom\'etries forment un espace vectoriel de dimension $n+1$
dont une base est la famille $\L_*$.
Le probl\`eme toujours en suspens (r\'esolu seulement pour $n\le 3$)
pos\'e dans [Gr-Sc]
est celui de la validit\'e (d'une formulation \'equivalente) de ce
th\'eor\`eme de finitude dans le cas sph\'erique.
La solution positive de ce probl\`eme aurait pour cons\'equence 
imm\'ediate
que les invariants de la suite $\Ll_*$
sont des combinaisons lin\'eaires de ceux de $\sigma_*$.
\ev
Dans cette direction, pour $X_y$ un germe d'ensemble
sous-analytique quelconque,
nous montrons  que chaque \'el\'ement d'une
des suites de $\Ll_*(X_y)$ ou $\sigma_*(X_y)$ est en effet
combinaison
lin\'eaire (\`a coefficients universels) des \'el\'ements de l'autre.
Plus pr\'ecis\'ement (Th\'eor\`eme
$3.1$) il existe une matrice triangulaire sup\'erieure
$\r M$ n'ayant que des $1$ sur sa diagonale, telle que :
$$\Ll_*=\r M\cdot \sigma_* .$$
Puisque la derni\`ere ligne de cette \'egalit\'e matricielle redonne
 la formule de  Cauchy-Crofton locale $\Theta_d=\sigma_d$,
on peut voir  $\Ll_*=\r M\cdot \sigma_* $ comme une {\sl formule de
Cauchy-Crofton locale multi-dimensionnelle}.

Pour une pr\'esentation plus large des questions li\'ees aux
valuations d\'efinies sur les convexes ainsi que pour le calcul
des coefficients de la matrice $\r M$, nous renvoyons le lecteur
\`a l'appendice.
\ev

En tant que moyenne, sur les projections g\'en\'eriques,
d'int\'egrales (relativement \`a la densit\'e locale)
de fonctions constructibles sur les projet\'es de $X_y$,
les termes de la suite $\sigma_*(X_y)$ sont
 les analogues r\'eels des
caract\'eristiques \'evanescentes introduites dans
le cadre analytique complexe par M. Kashiwara  dans
[Ka1] et \'etudi\'ees dans [Du1], [Du2], [Br-Du-Ka],
[L\^e-Te1,2,3] (voir aussi [Me]). Dans le cas des hypersurfaces
\`a singularit\'e isol\'ee
il s'agit (\`a coefficients pr\`es) de la suite $\mu_*$ des nombres de
Milnor des sections planes de $X_y$.
Il est montr\'e dans [Br-Sp] que pour une famille analytique
d'hypersurfaces analytiques complexes $(X_y)_{y\in \hbox{\bbbb C}}$
\`a singularit\'e isol\'ee et v\'erifiant la condition de Whitney
le long de l'axe des param\`etres $y$, la suite $\mu_*$ est constante.
Dans le cas g\'en\'eral, et non plus seulement
dans celui des hypersurfaces, d'apr\`es [Du2], [L\^e-Te1,3], chaque terme
de la suite $\ti \sigma_*(X_y)$ des invariants polaires complexes
est combinaison lin\'eaire des \'el\'ements de
la suite $\ti  m_*(X_y)$ des multiplicit\'es des vari\'et\'es polaires
de $X_y$. Ces derni\`eres \'etant constantes le long de strates de
Whitney  ([He-Me1], [Na2], [Te2]),
la suite $\ti \sigma_*(X_y)$ est elle-m\^eme constante
le long de strates de Whitney d'un ensemble analytique complexe.

\ev
Nous donnons ici  la version r\'eelle de ce r\'esultat
(Th\'eor\`emes $4.9$ et $4.10$), en montrant
que le long d'une strate de Verdier $Y$
d'un ensemble sous-analytique ferm\'e
$X$, l'application  $Y\ni y\mapsto \sigma_*(X_y)$ est continue,
g\'en\'eralisant \`a toutes les
dimensions la continuit\'e, prouv\'ee dans [Co2], de
la seule fonction
$Y\ni y\mapsto \sigma_d(X_y)$.
En particulier, chaque courbure localis\'ee $\Ll_i$ \'etant
combinaison lin\'eaire des invariants $\sigma_j$,
 nous en d\'eduisons la continuit\'e de $y\mapsto \Ll_*(X_y)$ le long
des strates $Y$ d'une stratification de Verdier de $X$.
\ev
Notons que bien qu'en g\'eom\'etrie complexe
il y ait \'equivalence entre la constance des
carac\-t\'eristiques \'evanescentes le long d'une strate et
le fait que celle-ci soit une strate d'une stratification de
Whitney, on ne peut
esp\'erer une telle r\'eciproque en g\'eom\'etrie r\'eelle.
Par exemple dans $\d R^3$ si $Y$ est l'axe $Oy$ et si
$X$ est le semi-alg\'ebrique suivant
(une fronce le long de l'axe $Oy$, pinc\'ee sur l'axe $Oz$) :
$$X=\{ (x,y,z)\in \d R^3; yz^3=x^3-3z^3x, \ z\ge 0 \},$$
avec les notations de la D\'efinition 2.8
on v\'erifie que :

- $\forall y\in Y$, $\sigma_1(X_y)=1$.
Ceci est clair pour $y\not=0$. Mais pour $y=0$ on observe,
suivant les notations du Th\'eor\`eme $2.8$, que pour des
projections g\'en\'eriques sur des droites $P$ de $\d R^3$
(celles telles que $Oz\not\subset P^\perp$), $n_P=2$,
$\chi^P_1=\chi^P_2=1$ et $\Theta_1(\r K_1^P)=\Theta_1(\r K_2^P)=1/2$.
C'est-\`a-dire que $\sigma_1(X_0)=1$.

- $\forall y\in Y$, $\sigma_2(X_y)=\Theta_2(X_y)=1/2$. Encore une fois
ceci est clair pour $y\not=0$, puisque dans ce cas $X_y$
est une hypersurface \`a bord de bord $Oy$. Pour $y=0$,
on calcule $\Theta_2(X_0)$ en calculant la $2$-densit\'e
des composantes du  c\^one tangent pur de $X_0$ affect\'es
de leur multiplicit\'e. On obtient $\Theta_2(X_0)=1/2.1/2+
1/2.1/2+3.0=1/2$ (cf [Ku-Ra]).

Les suites $\Ll_*$ et $\sigma_*$ sont par cons\'equent continues
(et m\^eme constantes) le long de $Y$, sans pour autant que
$X^{\hbox{\acc r\'eg}}$ soit $(w)$ ou $(b)$-r\'egulier le long de $Y$,
puisque pas m\^eme $(a)$-r\'egulier. Notons de plus que sur cet
exemple les discriminants g\'en\'eraux (de dimension~$1$) ne
satisfont pas la condition
g\'eom\'etrique du Th\'eor\`eme $4.9$, dont on montre qu'elle
suffit
\`a la continuit\'e de nos invariants locaux.

\ev
{\bf Contenu. } Cet article s'organise de la fa\c con suivante :
\ev
Dans la Section 1,
 nous montrons comment les courbures de
Lipschitz-Killing $\L_i(X)$ d'un ensemble $X$ d\'efinissable
dans une structure
o-minimale sur les r\'eels se localisent en des invariants
$\Ll_i(X_x)$ attach\'es au germe de $X$ en $x$. Nous traitons
explicitement le cas des ensembles sous-analytiques \`a l'aide
du th\'eor\`eme d'isotopie de Thom-Mather, ce qui n'est pas restrictif
puisque les ensembles d\'efinissables dans une structure
o-minimale sur $(\d R,+,.)$ admettent des
stratifications de Whitney ([Lo], [Sh] par exemple), mais on pourrait
aussi bien invoquer le th\'eor\`eme de d\'ecomposition cellulaire
([Dri], [Pi-St], [Kn-Pi-St]), dont
les cons\'equences en termes de finitude uniforme suffisent
pour nos preuves.

Dans la Section 2,
 nous d\'efinissons les {\sl invariants polaires}
$\sigma_i(X_x)$ (Section 2-a). Il s'agit essentiellement de remarquer que
les directions de projection g\'en\'erales ne rencontrent pas
les vari\'et\'es polaires  du germe $X_x$ qu'elles d\'efinissent
 (Proposition 2.2), pas plus que les
lieux critiques du link de $X_x$ associ\'es \`a ces projections
(Proposition 2.6). Nous interpr\'etons ensuite (Section 2-b) les
invariants $\sigma_i$ dans le cas analytique complexe et nous rappelons que leur
constance le long des strates d'une stratification \'equivaut \`a la
$(b)$-r\'egularit\'e de celle-ci (Th\'eor\`eme 2.12).

Dans la Section 3,
nous montrons le Th\'eor\`eme 3.1 \'evoqu\'e en introduction,
c'est-\`a-dire l'existence d'une
matrice triangulaire $\cal M$ telle que : $ \Ll_*=\cal M \cdot\sigma_*$.
Nous commen\c cons par montrer cette \'egalit\'e pour les ensembles
d\'efinissables (sous-analytiques dans le texte) coniques (Section 3-a)
en utilisant les techniques de calcul de [Br-Ku].
Nous montrons que dans le cas conique $\Ll_i(X_x)$ et
$\sigma_i(X_x)$ sont des combinaisons lin\'eaires des courbures
de Lipschitz-Killing $\L_j(L)$ du link $L=X\cap S_{(x,1)}$ de $X$.
  Le cas g\'en\'eral
(Section 3-b) s'en d\'eduit par   d\'eformation
tranverse sur le c\^one tangent.

Dans la Section 4,
nous \'etablissons un lien entre la variation de $\Ll_*$ et $\sigma_*$
le long d'une strate d'une stratification de $X$ et la r\'egularit\'e
de cette stratification : la $(w)$-r\'egularit\'e assure la continuit\'e
 des
invariants polaires. Pour cela nous montrons dans un premier temps
comment la condition $(b^*)$ permet de calculer $\sigma_i(X_x)$, 
pour $x$ au voisinage
de $0$ et le long d'une strate,  \`a l'aide, pour chaque direction
de projection
g\'en\'erique, d'un bon repr\'esentant du seul
germe $X_0$ en $0$ (Proposition 4.4). Au passage, l'existence de
tels bons voisinages assure imm\'ediatement dans le cas complexe la
constance
des caract\'eristiques \'evanescentes le long de strates de Whitney
(Corollaire 4.5) et ce sans mentionner la constance des multiplicit\'es
des vari\'et\'es polaires et la formule qui les relie aux
 caract\'eristiques \'evanescentes (voir [Du2], [L\^e-Te3]).
 Dans un second temps,
et en s'appuyant sur l'existence de bons voisinages,
nous montrons qu'une condition suffisant \`a la continuit\'e des
invariants polaires r\'eels $\sigma_i$ porte sur le contr\^ole
de la dimension des fibres des c\^ones normaux aux discriminants
(Th\'eor\`eme~4.9)~: il s'agit de la version g\'eom\'etrique de
l'\'equimuliplicit\'e des vari\'et\'es polaires et des discriminants,
 encore
pertinente en g\'eom\'etrie r\'eelle, et le Th\'eor\`eme 4.10
montre que cette condition est en
g\'eom\'etrie r\'eelle une cons\'equence de la
$(w)$-r\'egularit\'e, comme c'est d\'ej\`a le cas en complexe.

\ev
{\bf Notations. }
 Nous noterons $B_{(x,r)}$ ou au besoin
$B^n_{(x,r)}$ la boule ferm\'ee de centre $x$ et de rayon $r$
dans $\d R^n$, $S^{n-1}_{(x,r)}$ sa fronti\`ere, $\r H^i$ la
mesure $i$-dimensionnelle de Hausdorff, $\alpha_i=\r H^i(B^i_{(0,1)})$,
avec la convention $\alpha_0=1$,
$C_i^j=\D{j! \over (j-i)!i!}$ pour deux entiers
$0\le i\le j$, 
$\chi $ la caract\'eristique d'Euler Poincar\'e,  $1_E$ la fonction
caract\'eristique d'un sous-ensemble $E$ de $\d R^n$,
$e(X,x)$ la multiplicit\'e en $x$ d'un ensemble analytique complexe $X$,
$\r O_n(\d R)$ le groupe orthogonal de $\d R^n$,
$U_n(\d C)$ le groupe unitaire de $\d C^n$, $G(i,n)$ la grassmannienne
des $i$-plans vectoriels de $\d R^n$, $\bar G(i,n)$ la
grassmannienne
des $i$-plans affines de $\d R^n$, $\widetilde G(i,n)$ la grassmannienne
des $i$-plans vectoriels (complexes) de $\d C^n$, $\gamma_{i,n}$ la
mesure unitaire naturelle sur $G(i,n)$, $\bar\gamma_{i,n}$ la
mesure naturelle sur $\bar G(i,n)$,
et pour $P$ un $i$-plan de $\d R^n$ ou $\d C^n$,
$\pi_P$ la projection orthogonale  sur~$P$.

\Ev\vskip1cm
\centerline{\bf Sommaire}
\Ev

\centerline{ \hbox to 12cm{
{\titt 0. Introduction } \dotfill \hskip0,2cm {\titt 1}}}

\centerline{ \hbox to 12cm{
{\titt1. Invariants de Lipschitz-Killing locaux } \dotfill \hskip0,2cm
{\titt 5}}}

\centerline{ \hbox to 12cm{
{\titt 2. Invariants polaires  } \dotfill \hskip0,2cm
{\titt 8}}}

\centerline{ \hbox to 12cm{\hskip1cm
{\titt 2-a. Invariants polaires r\éels  } \dotfill \hskip0,2cm
{\titt 8}}}

\centerline{ \hbox to 12cm{\hskip1cm
{\titt 2-b. Invariants polaires complexes  } \dotfill \hskip0,2cm
{\titt 16}}}

\centerline{ \hbox to 12cm{
{\titt 3. Courbures de Lipschitz-Killing locales et
invariants polaires  } \dotfill \hskip0,2cm
{\titt 18}}}

 \centerline{ \hbox to 12cm{\hskip1cm
{\titt 3-a. Le cas conique  } \dotfill \hskip0,2cm
{\titt 18}}}

\centerline{ \hbox to 12cm{\hskip1cm
{\titt 3-b. Le cas g\'en\'eral  } \dotfill \hskip0,2cm
{\titt 27}}}

\centerline{ \hbox to 12cm{
{\titt 4. Conditions de r\'egularit\'e et invariants locaux }
\dotfill \hskip0,2cm
{\titt 28}}}

\centerline{ \hbox to 12cm{
{\titt Appendice.   Poly\`edres et valuations sph\'eriques
- Calcul des coefficients $m_i^j$ }
\dotfill \hskip0,2cm
{\titt 38}}}

\centerline{ \hbox to 12cm{
{\titt R\'ef\'erences }
\dotfill \hskip0,2cm
{\titt 45}}}

\vskip1cm
\Ev
\centerline{\bf 1. Invariants de Lipschitz-Killing locaux}
\Ev

Dans cette partie nous montrons comment les propri\'et\'es de
finitude uniforme locale des ensembles sous-analytiques
(ou plus largement des ensembles d\'efinissables dans une structure
o-minimale [Dr-Mi], [Sh]) permettent de localiser les courbures de
Lipschitz-Killing. Nous mentionnons [Be-Br2], Section 5, pour une notion
comparable de localisation des courbures de Lipschitz-Killing :
les localisations de [Be-Br2] sont des combinaisons lin\'eaires
des notres, comme cela est indiqu\'e dans [Be-Br2].
Nous   rappelons auparavant les d\'efinitions
et propri\'et\'es essentielles de ces courbures.
Sur ce sujet nous renvoyons
sans plus le r\'ep\'eter dans la suite,
 entre autres r\'ef\'erences \`a [Be-Br1,2], [Bl],  [Br-Ku],
[Ch-M\"u-Sc],  [Fe3], [Fu1,$\cdots$, 5], [La1,3],  [Laf],
[McMu-Sc], [Sc3,4],  [St1,2],  [We].

Consid\'erons $X$ un ensemble sous-analytique compact de
$\d R^n$ et
notons, pour $r$ un r\'eel~$>0$, $\r T_r(X) $ le voisinage
tubulaire de
rayon $r$ de  $X$, c'est-\`a-dire :
$$\r T_r(X)=\bigcup_{x\in X}B_{(x,r)}. $$

\n
Nous d\'efinissons alors la quantit\é $\r V_X(r) $
de la fa\c con suivante :
$$\r V_X(r)=
             \int_{x \in {\cal T}_r(X)  } \chi(X \cap B_{(x,r)})
\ d \r H^n(x).$$
Bien s\^ur, lorsque $X$ est lisse et pour $r$ suffisament
petit :
 $$\D \r V_X(r)=  \int_{x\in {\cal T}_r(X)  }  \
d \r H^n(x)=\r H^n(\r T_r(X)),$$
On appelle  $\r V_X(r) $  le {\sl volume modifi\'e}
de $\r T_r(X) $.
\n
D'apr\`es
[Fu4], [Be-Br2], [Br-Ku], quel que soit $r\ge 0$,
$\r V_X(r) $  est un polyn\^ome de degr\'e $n$
en la variable $r$, que nous notons :
$$\r V_X(r)=\L_n(X)+ \L_{n-1}(X).\alpha_1.r+ \ldots +\L_1(X).\alpha_{n-1}.r^{n-1}+
\L_0(X).\alpha_n.r^n .$$

\n
On dispose de plus
d'une formule de repr\'esentation int\'egrale pour chaque $\L_i(X)$
 :
\ev
{\bf Th\'eor\`eme 1.1. ---}
{\sl Soit $X$ un ensemble sous-analytique born\'e
de $\d R^n$. Pour tout $i \in \{0, \ldots , n \} $, on a l'\'egalit\'e :
$$\L_i(X)= \int_{\bar P \in \bar G(n-i,n)} \chi(X \cap \bar P) \ \
{d \bar\gamma_{n-i,n}(\bar P) \over \beta(n,i)},$$
o\`u $ \D \beta(n,i)$ est la constante universelle $ \D
\Gamma({n-i+1\over 2})\Gamma({i+1\over 2}) /
\Gamma({n+1\over 2}) \Gamma({1\over 2})  $ et   $ \Gamma$ est  la
fonction d'Euler.}
\ev

{\bf D\'efinition 1.2.}
 Les quantit\'es $\L_i(X)$ sont appel\'es les
{\sl courbures de Lipschitz-Killing de} $X$.

\ev
Supposons que l'origine de $\d R^n $ soit dans $X$, et notons $X_0$ le
germe \`a l'origine que d\'efinit~$X$.
Nous allons montrer dans cette section (Th\'eor\`eme $1.3$) que l'on peut
localiser les courbures $\L_i(X) $,
en \'etablissant que les limites
$\D \lim_{\u \to 0}{1\over \alpha_i.\u^i} \L_i(X \cap B^n_{(0,\u)}) $
existent. Nous  noterons  ces limites $\L_i^{\ell oc}(X_0) $, et nous les
appellerons les
{\sl invariants de Lipschitz-Killing locaux, ou les
courbures  de Lipschitz-Killing locales  du germe $X_0$}.
 Pour certaines valeurs de $i$ et de $d$, l'existence des limites $\Ll_i $
est claire. En effet, lorsque :
\ev
$ \bullet $  $i=0$ :
$\L_i(X\cap B^n_{(0,\u)})=\chi (X\cap B^n_{(0,\u)})=1, $
pour $\u  $ suffisamment petit, du fait de la structure conique de $X_0 $.
 De sorte que l'on peut poser  $\L_0^{\ell oc}(X_0)=1 $.
\ev

$\bullet$ $i>\dim(X_0) $ : l'intersection $X \cap \bar P $ est
g\'en\'eriquement vide, pour $\bar P \in \bar G(n-i,n) $, ce qui donne :
$\L_i(X\cap B^n_{(0,\u)})=0, \ \hbox{ et donc } \ \L_i^{\ell oc}(X_0)=0. $

\ev

$\bullet$  $i=d=\dim(X_0)$ :  l'intersection $X \cap \bar P $ est,
pour $\bar P \in \bar G(n-d,n) $ g\'en\'eral, un nombre fini de points
$N(\bar P, X)=N(P,X,y) $, lorsque $\bar P=\pi_P^{-1}(y) $.
On a alors :
$$\L_d(X\cap B^n_{(0,\u)})=
\int_{\bar P \in \bar G(n-d,n)} N(\bar P, X  \cap B^n_{(0,\u)}) \ \
{d \bar\gamma_{n-d,n}(\bar P) \over \beta(n,d)}$$
$$= \int_{ P \in  G(d,n)} \int_{y \in P} N( P, X  \cap B^n_{(0,\u)} ,y) \
d\r H^d(y) \ \
{d \gamma_{d,n}( P) \over \beta(n,d)} , $$
ce qui donne, d'apr\`es la formule classique de Cauchy-Crofton pour le
volume ([Fe1] 5.11, [Fe2] 2.10.15, [Sa] 14.69) :
$$\L_d(X\cap B^n_{(0,\u)})= \r H^d(X\cap B^n_{(0,\u)})  \ \hbox{ et donc  }$$
$$  \lim_{\u \to 0}{1\over \alpha_d.\u^d} \L_d(X \cap B^n_{(0,\u)})=
 \lim_{\u \to 0}{1\over \alpha_d.\u^d} \r H^d(X\cap B^n_{(0,\u)}).$$
Or d'apr\`es le th\'eor\`eme de Kurdyka et Raby ([Ku-Ra],
[Ku-Ra-Po]), qui sera par cons\équent obtenu  comme corollaire
  du
Th\'eor\`eme $1.3$ (Corollaire $1.4$), cette limite existe bien, il
s'agit de la $d$-densit\'e du germe $X_0 $
(cf [Fe2],  [Le] pour le cas analytique complexe). On la note
traditionnellement :
$$ \Theta_d(X_0)= \lim_{\u \to 0}{1\over \alpha_d.\u^d}
 \r H^d(X\cap B^n_{(0,\u)})
.$$
\ev
En r\'esum\'e, dans le Th\'eor\`eme $1.3$,
en localisant les courbures de Lipschitz-Killing de
$X$, nous d\'efinissons une suite finie d'invariants du germe $X_0$ :
$$(\L_0^{\ell oc}(X_0)=1, \L_1^{\ell oc}(X_0), \ldots , \L_d^{\ell oc}(X_0)=
\Theta_d(X_0),
0, \ldots, 0 ),$$  dont le terme de rang $d$ est
la $d$-densit\'e de $X_0$, c'est-\`a-dire
la localisation habituelle  du volume $d$-dimensionnel $\L_d=\r H^d$.

\Ev
{\bf Th\'eor\`eme 1.3. --- }
{\sl Soit $X$ un ensemble sous-analytique de
$\d R^n $, repr\'esentant quelconque
du germe $X_0$. Avec les notations pr\'ec\'edentes, quel que soit
$i \in \{0, \ldots, n \} $, la limite suivante existe :
$$ \lim_{\u \to 0}{1\over \alpha_i.\u^i} \L_i(X \cap B^n_{(0,\u)}). $$
Nous notons ces limites $ (\L_i^{\ell oc}(X_0))_{i \in \{0, \ldots
, n \} }$ et nous
les appelons les invariants de Lipschitz-Killing locaux ou les
courbures de
 Lipschitz-Killing locales du germe
$X_0$. De plus :
$$\L_0^{\ell oc}(X_0)=1, \
\L_d^{\ell oc}(X_0)=\Theta_d(X_0), \hbox{ pour } d=\dim(X_0)
\hbox{ et } \L_i^{\ell oc}(X_0)=0, \ \hbox{ pour } i>d.  $$}

{\bf Preuve.}
 Soit $i \in \{0, \ldots, n \} $.
Commen\c cons par d\'eformer $X$ sur son c\^one tangent  dans un produit.
Notons :
 $$X_\u=\D {1 \over \u}. (X \cap B^n_{(0,\u)})\subset B^n_{(0,1)} \hbox{ et }
[\bar G(n-i,n)]_1=\{ \bar P \in G(n-i,n); \ \bar P \cap B^n_{(0,1)} \not=
\emptyset\}.$$
 Le sous-analytique $[\bar G(n-i,n)]_1$ de $\bar G(n-i,n)$ est compact.
Consid\'erons maintenant les sous-analytiques suivants de
$G=[\bar G(n-i,n)]_1\times
[0,1]\times  B^{n-i}_{(0,1)} $ :
$$ E=\{ (\bar P,\u, x); \ \bar P  \in [\bar G(n-i,n)]_1,
\ \u \in ]0,1], \ x \in  X_\u \cap \bar P \},$$
$$F=\{ (\bar P,\u,x);
\ \bar P  \in [\bar G(n-i,n)]_1, \ \u \in ]0,1], \ x \in
 B^{n-i}_{(0,1)} \setminus (X_\u \cap \bar P)\}.$$
On a bien s\^ur : $G= \adh(E \cup F ).$

 La projection naturelle $p : G \lg [\bar G(n-i,n)]_1\times [0,1] $ est un
morphisme
sous-analytique propre. La fibre de $p$ au-dessus de $(\bar P,0)$ est
l'intersection
$\r C_0X \cap \bar P$ du c\^one tangent \`a l'origine de $X $ et de $\bar P $,
autrement dit au-dessus de $\{\bar P\} \times [0,1]$, $p$ est la d\'eformation de
$X \cap \bar P$ sur $\r C_0X \cap \bar P$.

Le morphisme $p$ est stratifiable d'apr\`es [Th], [Ha1,2], [Hi2] ou
[Go-McPh]~: il existe une stratification  de Whitney sous-analytique  $\Sigma$ de $G$
compatible avec
$E$ et $F$,  une stratification sous-analytique $\Sigma'$ de
$ [\bar G(n-i,n)]_1\times [0,1] $, telles que la pr\'eimage par $p$ de toute
strate
$\sigma' $ de $\Sigma'$ soit r\'eunion de strates de $\Sigma$, et la restriction
de $p$ aux strates de $ p^{-1}(\sigma')$ soit une submersion au-dessus de
$\sigma'$.

  Dans ces conditions le premier lemme d'isotopie de Thom-Mather
([Th], [Ma],
[Ha2]) assure que $p$ est topologiquement
triviale au-dessus de chaque strate $\sigma'$ de $\Sigma'$, de fa\c con
compatible
avec $\Sigma $. Comme $ G$ est compact, le nombre de strates de $\Sigma'$ est
fini, et au-dessus de chacune d'elles, deux fibres quelconques
$p^{-1}(\{ (\bar P, \u) \}) $, $p^{-1}(\{ (\bar Q, \eta) \}) $  sont
hom\'eomorphes
par un hom\'eomorphisme respectant les fibres de $p$ dans $E$ et $F$.
Or ces fibres dans $E$ sont pr\'ecis\'ement $X_\u \cap \bar P $ et
$X_\eta \cap \bar Q$; leur caract\'eristique d'Euler-Poincar\'e est ainsi
la m\^eme.

\n
 {\bf -} On en conclut que la famille d'entiers
$(\chi(X_\u \cap \bar P))_{\bar P \in
 [\bar G(n-i,n)]_1, \u \in ]0,1] } $ est une famille finie.

 \n
   De plus,  pour $\bar P$ fix\'e dans $[\bar G(n-i,n)]_1$,
le segment $\{\bar P\}\times [0,1] $
rencontre chaque strate de $\Sigma'$ un nombre fini de fois.

\n
 {\bf -}  On en conclut que la fonction
$]0,1] \ni \u \lg \chi (X_\u \cap \bar P)\in \d Z $ converge.

Pour terminer la preuve remarquons que par homog\'en\'eit\'e de $\L_i$, on a :
$$ {1 \over \u^i}\L_i(X \cap B^n_{(0,\u)})=
 {1 \over \u^i}\int_{\bar P \in \bar G(n-i,n)}\chi(X \cap B^n_{(0,\u)} \cap \bar
P)
\ { d\bar \gamma_{n-i,n}(\bar P) \over \beta(n,i) }$$
$$ =  \int_{\bar P \in [\bar G(n-i,n)]_1} \chi(X_\u  \cap \bar P)
\ { d\bar \gamma_{n-i,n}(\bar P) \over \beta(n,i) },$$
or on vient de voir que la famille de fonctions $\bigg([\bar G(n-i,n)]_1\ni
\bar P \lg
\chi(X_\u \cap \bar P)\bigg)_{\u \in ]0,1]}$ est domin\'ee par une fonction
born\'ee
 et  converge simplement. Le th\'eor\`eme de convergence domin\'ee permet de
conclure.
\f
\ev
  Le r\'esultat de [Ku-Ra] devient ainsi
un corollaire du th\'eor\`eme pr\'ec\'edent (voir [Li] pour
la premi\`ere mise en \oe uvre de la formule de Crofton
dans la preuve de l'existence de la densit\'e en tout point d'un ensemble
semi-pfaffien).
\ev
{\bf Corollaire 1.4. ---  }
{\sl Soit $X$ un ensemble sous-analytique de dimension $d$
dans $\d R^n$. La densit\'e $d$-dimensionnelle de $X$ en  tout point
de  $\d R^n$ existe. }
\ev
{\bf Preuve.}
On applique le Th\'eor\'eme $1.3$, avec $i=d$, en remarquant
d'apr\`es la
formule
de Cauchy-Crofton pour le volume,
que $\L_d(X \cap B^n_{(0,\u)})=\r H^d(X \cap B^n_{(0,\u)}) $. \f

\Ev
\centerline{\bf 2. Les invariants  polaires $\sigma_i(X_0) $}
\Ev

\centerline{\bf 2-a. Les invariants polaires r\éels}
\Ev

Comme pr\'ec\'edemment $X$ est un ensemble sous-analytique de $\d R^n$ qui
contient l'origine et qui repr\'esente le germe $X_0$ de dimension $d$.
Si $P$ est un  $i$-plan vectoriel de  $\d R^n$, nous identifierons
au besoin $\pi_P(0)$ et $0$.
Notons qu'il n'est  pas s\^ur que la projection  de $X_0$ sur
un $i$-plan vectoriel $P$ de $\d R^n$ soit bien d\'efinie; c'est-\`a-dire
que le germe en $0$ du projet\'e de $X\cap B^n_{(0,r)}$ sur $P$ n'est
peut-\^etre
pas ind\'ependant de $r$. Il suffit de penser  par
exemple \`a l'\'eclatement $X$ de $\d R^2 $ de centre l'origine, projet\'e
suivant la direction
$P^\perp= \d P^1(\d R)$ : si $B^3_{(x,r)} $ est une boule de $\d R^3$
centr\'ee en un point
$x$ de l'axe $P^\perp$ de l'h\'elice $X$, la projection sur $P$ de
$X \cap B^3_{(x,r)}$ n'est pas ind\'ependante de $r$.
Cependant cette situation n'est pas g\'en\'erique; on peut d\'efinir,
pour des projections
g\'en\'erales $\pi_P $  sur des  $i$-plans vectoriels $P$ de $\d R^n $,
le projet\'e du germe $X_0$,
c'est l'objet de la Proposition $2.4$.

\ev
{\bf D\'efinition 2.1.}
Soit $X$ un ensemble sous-analytique de $\d R^n$ de dimension
$d$, $i\in \{ 0, \ldots n \}$ et $P \in G(i,n)$. On note $ \r P_X(P)$ la
{\sl vari\'et\'e polaire de $X$ associ\'ee \`a $P$}. Il s'agit de
l'ensemble
des points critiques de $ \pi_{P| X^{\hbox{\acc r\'eg}} } :
X^{\hbox{\acc r\'eg}} \to P $.
C'est-\`a-dire que lorsque $i\le d $ :
$$ \r P_X(P)=\adh \{x \in  X^{\hbox{\acc r\'eg}}; \
\dim( \T_x X^{\hbox{\acc r\'eg}} \cap P^\perp) \ge d-i+1 \}, $$
 et lorsque $d<i$,
$ \r P_X(P)=  \adh(X)$.
Nous noterons $\r D_X(P)$ l'image de $ \r P_X(P)$ par la projection $\pi_P$ qui
d\'etermine $ \r P_X(P)$ et nous dirons que $\r D_X(P)$ est l'{\sl image
polaire de
$X$ associ\'ee \`a~$P$}.
\ev
Nous rappelons maintenant le  r\'esultat
qui stipule que les vari\'et\'es polaires sont transverses aux
directions
auxquelles elles sont associ\'ees, ce qui nous permettra de d\'efinir les
germes
$\big(\pi_P(X\cap B^n_{(0,r)})\big)_0$ et
$\big(\r D_{X \cap B^n_{(0,r)}}(P)\big)_0$ ind\'ependamment de $r>0$
(Propositions $2.4$ et $2.5$).
\ev
{\bf Proposition 2.2.} {\bf (Transversalit\'e pour les vari\'et\'es
polaires absolues)}
{ \bf --- } {\sl  Avec les notations de la d\'efinition
pr\'ec\'edente, il existe  un ouvert sous-analytique
dense $\r F_X^i\subset G(i,n)$ tel que quel que soit $P\in \r F_X^i$,
il existe un voisinage ouvert
$\r U$ de $0$ dans $\d R^n$ tel que :
$$ (\r U\cap \r P_X(P) \cap P^\perp) \setminus \{ 0\} =\emptyset.$$}
\ev
{\bf Preuve.}
On proc\`ede par r\'ecurrence sur $i$.
Lorsque $ i\ge d$, on peut poser : $\r F_X^i =\{P\in G(i,n);
P^\perp\cap C_0X=\{0\}\},$
qui est bien dense (cf par exemple [Co2], Lemme $1.4$).
Soit $i_0\le d$ et
supposons montr\ée l'existence de  $ \r F_X^j $ pour $ i_0\le j\le d$.
 Soit $P\in \r F_X^{i_0}$. Par hypoth\`ese, il existe un voisinage
ouvert $\r U$  de $0$ dans $\d R^n$ tel que :
$ (\r U\cap \r P_X(P) \cap P^\perp) \setminus \{ 0\} =\emptyset.$
Soit $\ell$ une droite (g\'en\'erale) de $P$ ne coupant pas $ \r D_X(P)$, dans
un voisinage $\r U'$ de $0$ dans $P$. On en d\'eduit que
$\r U \cap \pi_P^{-1}(\r U')\cap (\ell\oplus P^\perp)$ ne coupe pas
$\r P_X(P)$. En notant $P'=(\ell \oplus P^\perp)^\perp$ et en
remarquant que $\r P_X(P')\subset \r P_X(P)$, on obtient:
$$ (\r U  \cap \pi_P^{-1}(\r U') \cap \r P_X(P') \cap P'^\perp)
\setminus \{ 0\} =\emptyset ,$$
ce qui termine la preuve, puisque $\dim(P')=i_0-1$. \f

\ev
On peut maintenant prouver que l'image d'un germe est encore un germe pour des projections
g\'en\'eriques.
\ev
Soit $X$ un
 ensemble sous-analytique born\'e de $\d R^n$ dont l'adh\'erence contient l'origine
et soient  $i \in \{ 0, \ldots , n\}$, $(X^j)_{j\in\{0, \cdots, k\}}$ une stratification
sous-analytique finie  de $\adh(X)$, sans condition de r\'egularit\'e
particuli\`ere.
Rappelons que d'apr\`es la Proposition $2.2$ et avec les notations
de celle-ci,
quel que soit $P$ dans  l'ouvert  sous-analytique dense
$\displaystyle \bigcap_{j=0}^k \r F^i_{X^j}$ de $G(i,n)$, quel que soit $j\in \{0, \cdots, k\}$,  il existe
$r_P>0$ tel que :  $(B_{(0,r_P)}\cap \r P_{X^j}(P)\cap P^\perp)\setminus \{0\}=\emptyset$.
\ev
{\bf Notation.}
 Dans la suite on notera $\r E^i_X$ pour
$\displaystyle \bigcap_{j=0}^k \r F^i_{X^j}$. Cet ensemble
est associ\'e \`a une stratification de $\adh(X)$, que
l'on pourra choisir selon la situation. On pr\'ecisera,
quand cela sera n\'ecessaire, la stratification \`a laquelle  $\r E^i_X$ est attach\'e.

\ev
{\bf Proposition 2.4. (Projections sans \'eclatement)  --- }
{\sl Soit $X$ un
 ensemble sous-analytique born\'e de $\d R^n$ dont l'adh\'erence contient l'origine
et  $i \in \{ 0, \ldots , n\}$.  Avec les notations qui pr\'ec\`edent,
quel que soit $P$ dans  l'ouvert  sous-analytique dense
$\r E^i_X$ de $G(i,n) $,  pour tout   $r \in  ]0,r_P]$, il existe $s>0$ v\'erifiant :
 $$ B_{(0,s)}\cap \pi_P(X \cap B_{(0,r)})=
B_{(0,s)}\cap \pi_P(X \cap B_{(0,r_P )}).  $$
Autrement dit, une projection g\'en\'erique de $X_0$ sur un $i$-plan de $ G(i,n)$
d\'efinit bien un germe dans $P_0 $. }
\ev
{\bf Remarque.} La preuve qui suit montre en r\'ealit\'e que d\`es que $\r U$ est un voisinage de $0$ dans
$\d R^n$ tel que pour tout $j\in \{0,\cdots, k\}$, $\big(\r U\cap P^\perp \cap \r P_{X_j}(P)\big)\setminus\{0\}
=\emptyset$, le germe $\pi_P(\r U\cap X) $ est bien d\'efini.
\ev
{\bf Preuve.}
Soit $ i \in\{0, \ldots, n \}$ et $P \in G(i,n)$.
S'il existe $r\in ]0,r_P[$ tel que quel que soit
$s>0$, $ B_{(0,s)}\cap \pi_P(X \cap B_{(0,r_P)})\not=
B_{(0,s)}\cap \pi_P(X \cap B_{(0,r)}) $, on obtient une suite de points
$(x^r_\ell)_{\ell\in \N} $ dans $X $ telle que, si $y^r_\ell=\pi_p(x^r_\ell) $,
$r<||x^r_\ell||<r_P$, $\D \lim_{\ell\to \infty}
y^r_\ell= 0 $, et $\pi_P^{-1}(\{ y^r_\ell \})\cap B^n_{(0,r)}=\emptyset $.
Quitte \`a extraire une sous-suite de $(x^r_\ell)_{\ell\in \N} $,
on peut supposer que $(x^r_\ell)_{\ell \in \N}$
converge vers $x_r \in \Big(P^\perp \cap  \adh(X)\cap B^n_{(0,r_P)} \Big)
\setminus \hbox{int}(B_{(0,r)})$.
Il existe alors un indice $j\in \{ 0, \ldots, k\}$ pour lequel $X^j$ contient
$x_r$.
Maintenant, si $P$ est dans $ \r E_{X^j}^i$, pour tout
 $x\in  X^j\cap P^\perp \cap B^n_{(0,r_P)}$,  $\T_xX^j$ est
transverse \`a $P^\perp$.
Dans ce cas :
\vskip2mm
{\bf -} Soit $i>\dim(X^j) $, et  $P^\perp$ et $X^j$ ne se coupent pas dans $B_{(0, r_P)}$, ce qui contredit
l'existence de $x_r$.

{\bf -} Soit
$i\le \dim(X^j) $, et tout un voisinage de $x $ dans $X^j$ se proj\`ete
sur tout un voisinage de $0$ dans $P$.
Si  $X$ est ferm\é, $X^j\subset \adh(X)=X$, donc tout un
voisinage de $x$ dans $X$  se proj\`ete  sur tout un voisinage de $0$ dans
$P$. Ce qui contredit l'existence de la
suite $(x^r_\ell)_{\ell \in \N}$.
\vskip2mm
\n
Si $X$ n'est pas ferm\'e, apr\`es avoir prouv\'e que le projet\'e du germe
de $F^0=\adh(X)$   d\'efinit bien un germe,  pour une projection g\'en\'erique,
on applique le m\ême raisonnement \`a
$F^1=\adh\big(\adh(X)\setminus X\big) $, puis \`a
$F^2=\adh \bigg( X^1 \setminus \big(\adh(X)\setminus X\big)  \bigg) $ etc...
Les projet\'es g\'en\'eriques de chacun des germes de ces ferm\'es d\'efinissent
bien des germes. Il en est alors alors de m\^eme du germe de $X$. \f

\ev
{\bf Proposition 2.5. --- }
{\sl  Soit $X$ un ensemble sous-analytique
de $\d R^n$ contenant l'origine, 
$i\in \{ 0, \ldots, n \}$ et $\r E_X^i\subset G(i,n) $
 l'ouvert  sous-analytique
dense introduit ci-dessus et relatif \`a une  stratification 
$(X^j)_{j\in\{0, \cdots, k\}}$ de  adh$(X)$. Quels que soient
$P \in  \r E_X^i $, $r \in ]0,r_P] $, il existe $s>0 $ v\'erifiant :
$$ \forall j\in \{0, \cdots, k\}, \ \ 
B^n_{(0,s)}\cap \r D_{X^j \cap B^n_{(0,r)}}(P)
=  B^n_{(0,s)} \cap \r D_{X^j \cap B^n_{(0,r_P)}}(P),$$
autrement dit pour des projections g\'en\'erales de $ G(i,n)$, les images
polaires du germe $X_0$ sont bien d\'efinies. }
\ev
{\bf Preuve.}
D'apr\`es la Proposition $2.4$, la condition
$ (B_{(0,r_P)}\cap \r P_{X^j}(P) \cap P^\perp)
\setminus \{ 0\}=\emptyset$ assure la conclusion
de la Proposition $2.5$.\f

\ev
{\bf Proposition 2.6.}
{\sl Soient $X$ un ensemble sous-analytique ferm\é de $\d R^n$  
contenant l'origine, $i \in \{ 0, \ldots, n\} $ et $P \in \r E_X^i $,
o\`u  $\r E_X^i $ est relatif \`a une  stratification $(a)$-r\'eguli\`ere
$(X^j)_{j\in\{0, \cdots, k\}}$ de $X$.
Il existe alors $r'_P >0$, tel que quel que soit
$r\in ]0, r'_P]$, quel que soit $j\in \{0, \cdots, k\}$~:
 $$0 \not \in  \r D_{L^r_{X^j}} (P),$$
o\`u
$L^r_{X^j}=S^{n-1}_{(0,r)} \cap X^j $ est le link de
$X^j$ dans $S^{n-1}_{(0,r)} $.
   }
\ev
{\bf Preuve.}
Soit $\rho>0$ tel que $S_{(0,\rho')}$
soit transverse aux strates de $(X^j)_{j\in \{0, \cdots, k\} }$, 
d\`es que $0<\rho'<\rho$. Avec les notations
  pr\'ec\'edentes, soit $r\in ]0,\min(\rho,r_P)[$.
  Commen\c cons par remarquer que quelle que soit la strate 
$X^j$, celle-ci est transverse \`a $\pi_P^{-1}(\pi_P(x))=\pi_P(x)+P^\perp$ en $x$, 
pour $x$ dans $S(0,r)$, pourvu que $\pi_P(x)$ soit 
suffisamment proche de $0$. Ceci r\'esulte de la $(a)$-r\'egularit\'e
de la stratification $(X^j)_{j\in\{0, \cdots, k\}}$, 
du fait que $X$ est ferm\'e
et du fait que $P^\perp$ ne coupe pas $\r P_{X^j}(P)$ dans $S(0,r)$.
Raisonnons par l'absurde pour prouver la Proposition  $2.6$ : si 
$0 \in \adh (\r D_{L^r_{X^j}} (P))$,
il existe une suite $(x_\ell^r)_{\ell\in \hbox{\bbbb N}}$ de 
$S_{(0,r)}\cap X^j$ telle que
 $\T_{x_\ell^r} L^{r}_{X^j}= \T_{x_\ell^r} S_{(0,r)}\cap
\T_{x_\ell^r}  X^j $ et $P^\perp$ sont non transverses
et $\D\lim_{\ell\to \infty}\pi_P(x_\ell^r)=0$.
Comme   $\pi_P(x_\ell^r)+P^\perp$ est transverse en 
$x_\ell^r$ \à $X^j$, on en d\'eduit que 
 : $\T_{x_\ell^r}  X^j\cap P^\perp \subset \T_{x_\ell^r} S_{(0,r)}$.
Une sous-suite de $(x_\ell^r)_{\ell\in \hbox{\bbbb N}}$ converge vers
un point $x^r$
de $S_{(0,r)}\cap X$, qui  est dans une state $X^m$ et 
\`a nouveau par la condition $(a)$, on a :  
$ \T_{x^r}X^m \cap P^\perp\subset  \T_{x^r} S_{(0, \Vert x^r\Vert)}$.
En faisant maintenant varier $r$,
on en d\'eduit une suite non constante
$(x_p)_{p\in \hbox{\bbbb N}}$ contenue dans une strate $X^m$, de 
limite $0$ et telle que :
$$ \T_{x_p}X^m \cap P^\perp\subset  \T_{x_p} S_{(0, \Vert x_p\Vert)}.$$ 
 Mais comme $P^\perp$ et $X^m$ sont transverses dans $B_{(0,r_P)}$, on a :
$$  \T_{x_p} (X^m \cap P^\perp)= \T_{x_p}X^m \cap P^\perp \subset  \T_{x_p} S_{(0, \Vert x_p\Vert)}.$$
Or $X^m \cap P^\perp$ est un sous-analytique de dimension $\not=0$
dont l'adh\'erence contient $0$, l'inclusion
$ \T_{x_p} (X^m \cap P^\perp)\subset  \T_{x_p} S_{(0, \Vert x_p\Vert)}$ 
est donc en contradiction avec le lemme de Whitney.\f
\ev
{\bf Remarque.} 
La preuve montre que l'on peut prendre 
$r'_P$ de sorte qu'en tout point $x$ de $B(0,r'_P)$,
$\T_x(X^m\cap P^\perp)$ et la direction $x$ ne sont pas orthogonaux, quelle
que soit la strate $X^m$.
\ev 
Nous rassemblons maintenant en un th\éor\ème les diverses propositions \établies,
et nous introduisons comme dans [Co2] le discriminant local $\Delta_{X_0}(P) $,
les profils polaires locaux $\r K_\ell^P $, et les caract\éristiques $\chi_\ell^P $
du germe $X_0$, associ\és \à une projection g\én\érale $\pi_P $, pour $P\in G(i,n) $
et $i\in \{ 0, \ldots, n\} $.
L'introduction de ces objets et la formule de repr\ésentation
int\égrale pour les
courbures de Lipschitz-Killing locales $ \L_i^{\ell oc}(X_0)$
que allons d\émontrer (Th\éor\ème $3.1$), nous permet de ramener
l'\étude du comportement des  courbures $ \L_i^{\ell oc}(X_0)$
\à la g\'eom\étrie
des discriminants $\Delta_{X_0}(P^i), \ldots,\Delta_{X_0}(P^d)  $,
(avec $P^i, \ldots, P^d$
g\'en\'eraux dans respectivement $G(i,n), \ldots, G(d,n) $)
qui est bien comprise le long de strates de Verdier (cf Th\'eor\`eme
$4.10$).

\ev
{\bf Th\éor\`eme  2.7. --- }
{\sl Soit  $X$ un ensemble sous-analytique ferm\é et born\é de $\d R^n $ de
dimension $d $
contenant l'origine et $(X^j)_{j \in \{0, \ldots , k \}}$
une stratification de Whitney de $X$. Soit  $ i \in  \{0, \ldots, n \}$.
\ev
\n
Avec les notations qui pr\'ec\`edent la Proposition $2.4$,
quel que soit $ P$ dans l'ouvert sous-analytique dense $\r E_X^i
=\D \bigcap_{j=1}^k \r E_{X^j}^i $ de  $G(i,n) $ :
\ev
\n
{(i)} \hskip2mm
 \vtop{\hsize12,95cm \n
Quel que soit $r>0$ suffisament petit, quel que soit
$j\in \{0, \ldots , k \} $,
le germe $[\r D_{X^j \cap B^n_{(0,r)}}(P)]_0 $ est ind\'ependant de $r$, 
on note $ \Delta_{X_0}(P) $
le germe  $[\D\bigcup_{j=0}^k\r D_{X^j \cap B^n_{(0,r)}}(P)]_0 $ et
on l'appelle le {\sl discriminant local de $X_0 $ associ\'e \`a $P$}. }
\vskip3mm
\n
{(ii)} \hskip2mm
 \vtop{\hsize12,95cm \n  Le germe $[\pi_P(X \cap B^n_{(0,r)})]_0\setminus \Delta_{X_0}(P) $
est bien d\'efini pour $r$ suffisamment petit. Il s'agit du germe d'un ouvert
sous-analytique de $P$. Les germes de ses composantes connexes sont not\'es $\r K_1^P,
\ldots, \r K_{n_P}^P $ et sont appel\'es les profils polaires locaux (\`a l'origine)
de $X_0 $ associ\'es \`a $P$.}

\vskip3mm
\n
{(iii)} \hskip2mm
 \vtop{\hsize12,95cm \n  Soit $\rho >0$ suffisamment petit pour 
que $S^{n-1}_{(0,\rho')}$
soit transverse \à chaque  $X^j$, lorsque
$\rho'\in ]0, \rho[ $.
Notons alors, pour $r\in ]0,\rho[ $ et pour $j \in \{0, \ldots , k \} $,
 $L^r_{X^j}=X^j\cap S^{n-1}_{(0,r)} $, puis  $\r D_{L^r}(P) =
\D \bigcup_{j=0}^k \pi_P(\r P_{L^r_{X^j}}(P))$ le discriminant du 
link $ L^r=X \cap S^{n-1}_{(0,r)}$.

\n Le compl\émentaire de  $\r D_{L_r}(P)  \cup
\D\bigcup_{j=0}^k\r D_{X^j \cap B^n_{(0,r)}}(P)$ dans $\pi_P(X \cap B^n_{(0,r)})$  est
un ouvert sous-analytique dense. On note chacune de ses composantes connexes
$\K_1^{P,r}, \ldots , \K_{N_P}^{P,r}$. Pour $r$ suffisamment petit, les germes des ouverts
$\K_j^{P,r}$ qui adh\èrent \à $0$ sont les profils polaires locaux  $\r K_1^P,
\ldots, \r K_{n_P}^P $ de $X_0$. De plus \à chaque $\r K_j^{P} $ on peut associer un
entier $\chi_j^{P} \in \d Z $ \égal \à $\chi (\pi_P^{-1}(y)\cap X\cap B^n_{(0,r)} ) $, cet
entier ne d\épendant ni de $r>0$, qui est suppos\é suffisamment
petit,  ni de $y\in\K_j^{P,r}  $, pourvu que $\Vert y\Vert <<r $.} }

\ev
{\bf Preuve.} (i) r\ésulte de la Proposition $2.5$ et  (ii) r\ésulte de la
Proposition $2.4 $.

\n Prouvons (iii) : par la Proposition $2.6$ les images polaires
des $L^r_{X^j} $ associ\ées \à $P$ n'adh\èrent pas \à l'origine, pour $r $ petit,
de sorte que les $\K_j^{P,r}$ qui adh\èrent  \à l'origine n'ont pas dans leur bord,
au voisinage de l'origine, de points communs avec les images polaires
des $L^r_{X^j}$; au voisinage  de l'origine le bord de  $\K_j^{P,r}$ est obtenu comme la
r\éunion des images polaires des  $X_j$. Les germes de tels  $\K_j^{P,r}$ donnent donc
les profils polaires locaux  $\r K_1^P, \ldots, \r K_{n_P}^P $ de $X_0$.

Enfin montrons que les entiers relatifs $\chi_j^{P,r}=
\chi (\pi_P^{-1}(y)\cap X\cap B^n_{(0,r)} )$ ne d\épendent ni de
$y\in \K_j^{P,r}  $, ni de $r$, pourvu que ce dernier soit choisi petit,
et que la partition $(X^0, \ldots , X^k) $ soit une stratification de Whitney.
On commence par remarquer que pour $r$ fix\'e, quels que soient $ y$ et $z$ dans $\K_j^{P,r} $,
les fibres $  \pi_P^{-1}( y)\cap X\cap B^n_{(0,r)} $   et
$ \pi_P^{-1}( z)\cap X\cap B^n_{(0,r)}  $ sont hom\'eomorphes. Il s'agit une fois de plus
du premier lemme d'isotopie de Thom-Mather; le morphisme sous-analytique propre
$ \pi_P : X\cap  B^n_{(0,r)} \to P\cap  B^n_{(0,r)}$ \'etant stratifi\'e par
$(X^j)_{j\in \{0, \ldots, k \}} $ et $(\K_j^{P,r})_{j\in \{1, \ldots, N_P \}}$.

On montre ensuite qu'\`a $P$ fix\'e dans $\r E^i_X $,  pour $r$ suffisamment proche de $0$,
l'entier $\chi_j^{P,r}$ associ\'e \`a $\K_j^{P,r} $,
o\`u $j\in \{1, \ldots, n_P \} $ ($\K_j^{P,r} $ adh\`ere \`a l'origine),
 est ind\'ependant de $r$.
Consid\'erons pour cela
$E=\{ (y,r); y\in  \adh(\K_j^{P,r}), r\in  [0,\u]    \}$  ($\u>0 $ suffisamment petit),
$F=\{ (x,r); \pi_P(x)\in \adh(\K_j^{P,r}), x\in  \cap X\cap B^n_{(0,r)}, r\in  [0,\u] \} $, et
$p: F \to E$  d\'efinie par $p(x,r)=(\pi_P(x), r) $. Il existe une stratification  $\Sigma'$
de $F$, une stratification finie $\Sigma $ de $E$, compatible avec
$E'=\{ (y,r); y\in   \K_j^{P,r}, r\in  [0,\u]    \} $, qui stratifient $p$. Soit $\sigma$ une strate
de $\Sigma $ qui soit contenue dans $ E'$ et qui adh\`ere au voisinage de l'origine \`a
tout un segment $ [0,\u']\subset [0,\u]$. Si $r,r'\le \u' $, et si $(y,r) $ et $(z,r') $ sont
suffisamment proche de $(0,r) $ et $(0,r')$ respectivement dans $(\K_j^{P,r}\times \{r\}) \cap
\sigma$ et $(\K_j^{P,r'}\times \{r\}) \cap \sigma$, les fibres $p^{-1}(y,r) $ et $p^{-1}(z,r')$
\'etant hom\'eomorphes et de caract\'eristique d'Euler-Poincar\'e respectivement $\chi_j^{P,r} $
et  $\chi_j^{P,r'} $, on en d\'eduit l'ind\'ependance de $\chi_j^{P,r} $ relativement \`a $r<\u' $.
\f
\ev
{\bf Remarque. }
Dans le Th\'eor\`eme 2.8.(iii), on peut choisir $r<r'_P$
($r'_P$ donn\'e par la Proposition $2.6$) et $y$ dans $P$
de sorte que $\Vert y \Vert < \delta $, o\`u $\delta>0$ minore
la distance  de $\r D_{L^r}(P)$ \`a $0$.
\ev 
L'introduction des profils polaires locaux $ (\r K_j^P)_{j\in \{1, \ldots, n_P \},
P\in {\cal E}_X^i}$ et des caract\éristiques
$(\chi_j^P)_{j \in \{ 1, \ldots, n_P  \},P\in {\cal E}_X^i}$, qui leur sont attach\'ees nous permet
de d\'efinir des invariants polaires~: les moyennes sur $G(i,n)$ des densit\'es des profils
polaires locaux de $X_0 $ associ\és \à $P$, affect\és des caract\éristiques locales
 $\chi_j^P$.

\ev
{\bf D\'efinition 2.8.}
Avec les notations du Th\'eor\`eme $2.8 $,
nous d\'efinissons, pour $i\in \{0, \ldots, n \}$, des
 {\sl invariants polaires}  $ \sigma_i(X_0)$ par :
$$  \sigma_i(X_0)=  \D\int_{P\in {\cal E}_X^i} \ \
\sum_{j=1}^{n_P} \ \chi_j^P \cdot \Theta_i(\r K_j^P) \ \ d \gamma_{i,n}(P) .$$
On note $\sigma_*(X_0)$ la suite $(  \sigma_i(X_0) )_{i\in \{0, \cdots,n\}}.$
\ev

Notons maintenant
$\r C(X)$ le groupe des fonctions constructibles
sur X, c'est-\`a-dire les fonctions
du type : $\varphi=\D\sum_{j=1}^Nn_j\cdot 1_{K^j}$, pour
 $n_j\in \d Z$,  $K^j$ des sous-analytiques de $X$, et
pour $Z\subset X$ et  $y\in Y$, $f_*( 1_Z)(y)
=\chi(f^{-1}(y)\cap Z)$.
Si on  consid\`ere le foncteur suivant de la cat\'egorie des ensembles
sous-analytiques compacts \`a
la cat\'egorie des groupes~:
\ev
$$\matrix{& X & -\!\!\!-\!\!\!\!-\!\!\!\!\longrightarrow &
\r C(X) &\cr
& f\downarrow \ \ \ && \downarrow f_*\cr
& Y &  -\!\!\!-\!\!\!\!-\!\!\!\!\longrightarrow & \r C(Y)} $$ \ev
\n
en vertu du Th\'eor\`eme $2.8$, pour $\pi_P$ une projection g\'en\'erale dans $G(i,n)$, ce diagramme admet
l'\'equivalent local : \ev
$$\matrix{& X_0 & -  \!\!\! -  \!\!\!  -  \!\!\!   \longrightarrow & \r C(X_0) & \cr
& \pi_{P_0}\downarrow \ \ \ && \downarrow \pi_{P_0*}\cr
& P_0  &  -\!\!\!-\!\!\!-\!\!\!\longrightarrow & \r C(P_0)  \cr }$$
\ev\n
o\`u pour $Z_0\subset X_0$ et $y\in P$, $\pi_{P_0*}( 1_{Z_0})(y)=\chi(\pi_P^{-1}(y)\cap Z\cap B(0,r))$,
$r$ \'etant suffisamment petit, et $0<\Vert y\Vert<<r$.
Si l'on note ensuite par $\theta(\varphi)$ l'int\'egrale relativement
\`a la densit\'e locale en $0$ d'un germe $\varphi:P_0 \to \d Z$ de
fonction constructible,
c'est-\`a-dire : $\theta(\varphi)= \D\sum n_j\cdot \Theta(K^j_0)$
 lorsque $\varphi=\D\sum_{j=1}^N n_j\cdot 1_{K^j_0}$, pour des germes
d'ensembles sous-analytiques $K^j_0\subset P_0$, on obtient la formule :
$\sigma_i(X_0)=\D\int_{P\in G(i,n)} \theta( \pi_{P_0*}(1_{X_0}) )\ d P.$
 \ev
 La figure ci-dessous compare, pour $X_0$ un germe de surface de
$\d R^3$ et
pour une projection $\pi_P$
particuli\`ere sur un $2$-plan de $\d R^3$, ce qui est respectivement
pris en compte dans le calcul de $\sigma_2(X_0)$ et $\Ll_2(X_0)$.

\Ev
$$\hskip-8,5cm  X\cap B^n_{(0,r)}$$
\vskip-1,5cm
\vglue0cm
 \vglue0cm
\vskip12cm
\hskip0cm\includegraphics{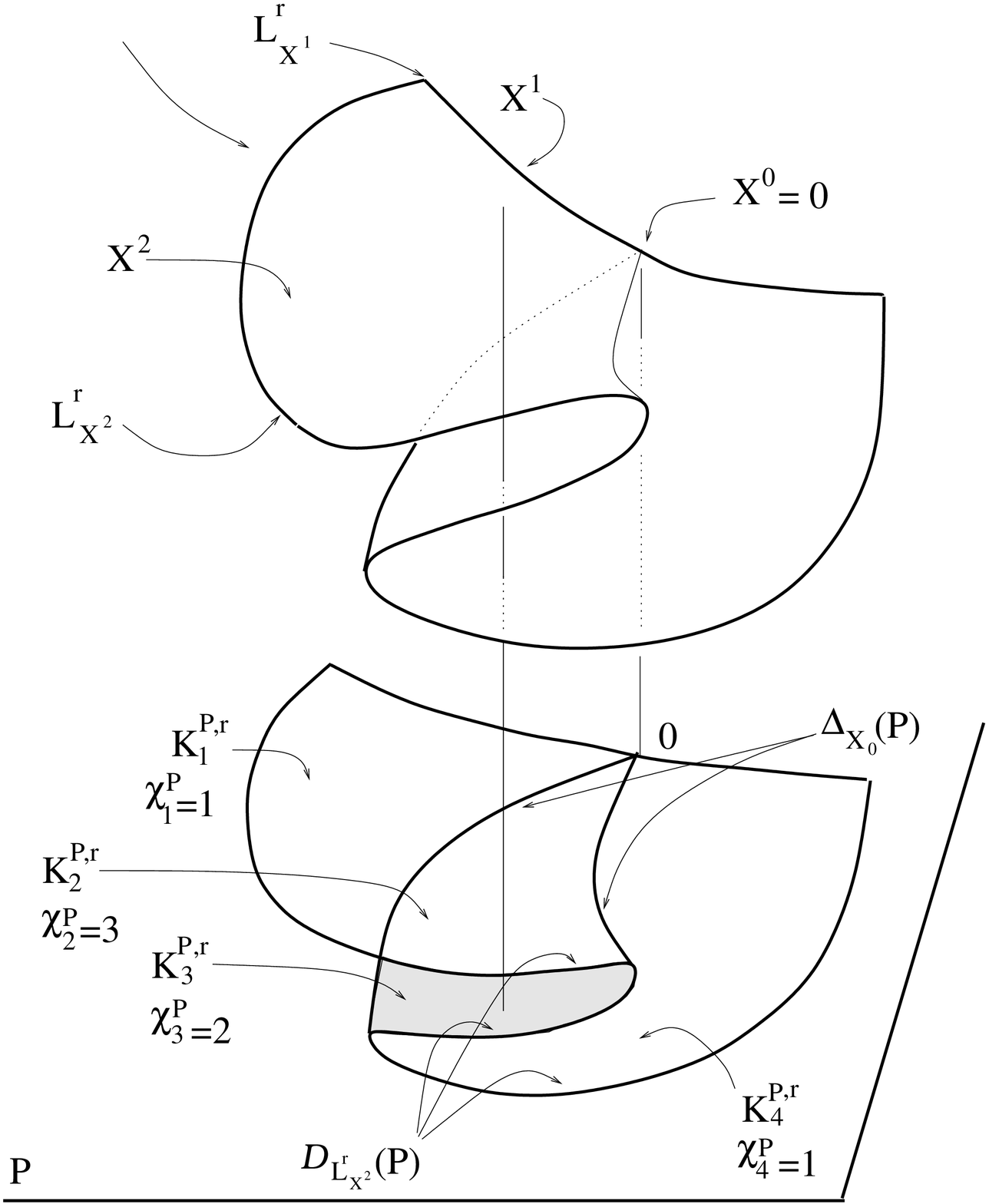}
\vglue0cm
 \vskip0,8cm

\n
 Sur cet
exemple, on a :

\n
- $\theta( \pi_{P_0*}(1_{X_0}))= \chi_1^P\cdot \Theta_2((K^P_1 )_0 )+
\chi_2^P\cdot \Theta_2((K^P_2)_0 ) +
\chi_4^P\cdot \Theta_2((K^P_4)_0 )  $ et
$\sigma_2(X_0)= \D\int_{Q\in G(2,3)} \theta( \pi_{Q_0*}(1_{X_0}) )\ d Q.$
Notons que dans le calcul de $\theta(\pi_{P_0*}(1_{X_0}))$,
pour la projection sur $P$,
le domaine $K^P_3$ (en gris sur le dessin) n'est pas pris en compte, puisque seuls les
profils polaires locaux adh\'erents
\`a l'origine comptent dans la d\'efinition de
$\sigma_i$. Ces profils sont ceux d\'etermin\'es par
$\Delta_{X_0}(P)$.
\ev
\n
- En revanche dans le calcul de $\Ll_2(X_0)$ le terme
$\chi_3^P\cdot \r H^2(K^P_3)$
 est pris en compte, car tous les
domaines $K^{P,r}_j$ interviennent dans ce calcul, y-compris ceux qui ne contiennent pas $0$ dans leur adh\'erence.
En notant $v(\pi_{P*}(1_{X\cap B_{0,r}}))= \D\sum_{j=1}^4
\chi_j^P\cdot \r H^2(K^j_0)$, on a :
$\Ll_2(X_0)=\D \lim_{r \to 0} {1\over \alpha_2r^2}
\int_{Q\in G(2,3)} v( \pi_{Q*}(1_{X\cap B^n_{(0,r)}}) )\ d Q.$

\ev
{\bf Remarques. }
Trivialement, lorsque $i=0$,  $\sigma_i(X_0)=1=\Ll_0(X_0)$ et
lorsque $i=n$,  $\sigma_i(X_0)=\Theta_n(X_0)=\Ll_n(X_0)$. Enfin
lorsque $i=d=\dim(X_0)$,  $\sigma_i(X_0)=\Theta_d(X_0)$, c'est-\à-dire
$\sigma_i(X_0)=\Ll_d(X_0)$. Cette derni\`ere \'egalit\'e, contrairement
aux deux autres
n'est pas imm\'ediate (cf [Co2]) :
\ev
{\bf Th\éor\ème 2.9. (Formule de Cauchy-Crofton locale)}
([Co2], Th\éor\èmes $1.10$ et $1.16$)  {\bf  --- }
{\sl Soit $X$ un sous-ensemble sous-analytique de $\d R^n$ de dimension $d$
et soient $\r G \subset G(d,n)$ un sous-ensemble sous-analytique de $G(d,n)$
sur lequel agit transitivement un sous-groupe $G$ de $\r O_n(\d R)$ et
$\mu_{d,n}$ une mesure $G$-invariante sur $\r G$, tels que : \par
{\bf -} les espaces tangents \à $\r C_0X_0$ sont dans $\r G$, \par
{\bf -} il existe $P^0\in \r G$ dont le fixateur $G_{P^0}$ agit transitivement
sur le $d$-espace vectoriel sous-jacent \à $P^0$ et $\mu_{d,n}(\r G)=
\mu_{d,n}(\r G \cap \r E^d_X)=1$. \par
L'\égalit\é suivante a alors lieu :
$$ \int_{P \in \r G \cap \r E^d_X} \sum_{j=1}^{n_P}\chi_j^P\cdot
\Theta_d(\r K_j^P) \ \ d \mu_{d,n}(P)=\Theta_d(X_0). $$}

\n
Dans le cas o\`u $\r G=G(d,n)$ et
 $G=\r O_n(\d R)$, la formule donne : $\sigma_d(X_0)=\Ll_d(X_0) $.

\n
Dans le cas o\`u $X$ est analytique complexe,
$\r G=\widetilde G(d/2,n)$ et
 $G=U_n(\d C)$, la formule donne~: $\Theta_{d}(X_0)=e(X,0)$.
\ev
{\bf Remarque 2.10.}
 Lorsque $(X^j)_{j\in \{0, \cdots, k\}}$ est une stratification de 
Whitney
de $X$,  $0\in X^0 $, $\sigma_i(X_0)=1$, lorsque $i<\dim( X^0 )$.
En effet, soit une direction de projection g\'en\'erale 
$P^\perp$ ($P\in \r E^i_X$), transverse \`a $X^0$ en $0$.
Comme la stratification   $(X^j)_{j\in \{0, \cdots, k\}}$ est
$(a)$-r\'eguli\`ere, $\pi_P$ est transverse
aux strates adjacentes \`a $X^0$ en $0$, et d'apr\`es la Proposition 
$2.6$, $\pi_P$ est aussi transverse
aux strates du link $X\cap S^n_{(0,r)}$ induites par
$(X^j)_{j\in \{0, \cdots, k\}}$, ceci
au-dessus d'un voisinage $P\cap B^n_{(0,\eta)}$ de $0$
dans $P$, pour $0<\eta<<r$. Il en r\'esulte que
la restriction $\pi^r_P$ de $\pi_P$
\`a $X\cap B^n_{(0,r)}$, pour $r$ suffisamment petit est une
submersion propre au-dessus de
$P\cap B^n_{(0,\eta)}$,
stratifi\'ee par
$ ( X^j\cap \hbox{int}( B^n_{(0,r)} ),X^j\cap S^n_{(0,r)} )_{j\in \{0, \cdots, k\}}  $ et les projections de
ces strates dans $P\cap B^n_{(0,\eta)}$.
Les fibres de $\pi_r$ au dessus de $P\cap B^n_{(0,\eta)}$ sont toutes hom\'eomorphes
et puisque $\pi_r^{-1}(\{0\})$ est contratile, on a bien~: $\sigma_i(X_0)=1$.
En conclusion, si $(X^j)_{j\in \{0, \cdots, k\}}$ est une stratification de Whitney
de $X$ et si $d_0$ est la dimension de la strate qui contient $0$, on a :
$$\sigma_*(X_0)=(1, \cdots, 1, \sigma_{d_0+1}(X_0),
\cdots, \sigma_{d-1}(X_0), \Ll_d(X_0), 0, \cdots, 0).$$

\ev
{\bf Remarque.}
Au m\ême titre que les invariants $\Lambda_i$, et donc au m\ême titre que les
courbures locales $\Ll_i$, les invariants polaires $\sigma_i$ sont
{\sl intrins\èques} : si le germe $X_0\subset \d R^n_0 \subset \d R^{n+1}_0$,
$\sigma_{i,n}(X_0)=\sigma_{i,n+1}(X_0) $, o\ù $\sigma_{i,n}(X_0)$ d\ésigne
l'invariant $\sigma_i(X_0)$, calcul\é gr\^ace \à $G(i,n)$. Il s'agit \à
nouveau d'une cons\'equence de la formule de Cauchy-Crofton pour la densit\é.
Comme c'est le cas pour les invariants $\Ll_i$, il est facile de
s'assurer que les invariants $\sigma_i$
ne d\épendent pas de la dimension de l'espace
euclidien dans lequel $X_0$ est consid\'er\'e.
De plus, pour $i\in \{0, d, d+1,\cdots,  n\}$,
les invariants polaires $\sigma_i$ sont
\égaux aux courbures locales~$\Ll_i$.
     Nous allons \'etablir dans la section suivante, au   Th\'eor\`eme $3.1$, que les
courbures $\Ll_i$ s'obtiennent en r\éalit\é comme des combinaisons lin\'eaires
(\à coefficients universels) des
invariants $\sigma_j$; les \égalit\és et les propri\'et\'es
communes que nous venons d'observer ne sont que des cas sp\éciaux de ce
th\'eor\`eme.

\Ev
\centerline{\bf 2-b.  Les invariants polaires complexes.}
\Ev
On peut bien s\^ur d\éfinir les invariants polaires  $\sigma_i(X_0)$
dans le cadre complexe, lorsque $X_0$ est un germe d'ensemble analytique
complexe \`a l'origine de $\d C^n$; pour cela il suffit de ne consid\'erer,
 dans la D\'efinition $2.8$, que les $i$-plans complexes de $\d C^n$
appartenant \`a l'\'equivalent complexe  $\D \widetilde{\r E^i_X}$ de
l'ensemble  $ \r E^i_X$ du Th\'eor\`eme $2.8$.
Pour de tels plans $P$ g\én\ériques, il existe un seul domaine
$\r K^P$ au-dessus duquel la caract\éristique
d'Euler-Poincar\é typique
$\chi^P=\chi (\pi_P^{-1}(y)\cap X\cap B^{2n}_{(0,r)})$ ne d\'epend ni
de $y$ g\én\érique dans $P\setminus \{0\}$ et suffisamment proche de $0$, ni du choix du repr\ésentant $X$ de  $X_0$,
c'est-\`a-dire du choix de  $r>0$, lorsque celui-ci est suffisamment petit. De plus $\chi^P$
ne d\'epend pas du choix de $P$, lorsque $P$ est g\'en\'eral. Ceci r\'esulte
des m\^emes arguments que dans le cas r\'eel et du fait que le
compl\'ementaire d'un ensemble analytique complexe de
codimension plus grande que $1$ est connexe.
Notons $\widetilde{\sigma_i}(X_0)$ l'invariant polaire complexe
d'ordre $i$.
Comme $\Theta_{2i}(\r K^P)=1$, on obtient~:
$$ \widetilde{\sigma_i}(X_0)=
\chi
(\pi_P^{-1}(y)\cap X\cap B^{2n}_{(0,r)}) $$
pour $y\in P$ g\én\érique
suffisamment proche de l'origine et $r>0$ suffisamment
petit.
\ev

Dans le cas particulier important o\`u  $X$  est l'hypersurface $f^{-1}(0)$, donn\ée par
une application analytique $f:\d C^n \to \d C$ :
$$ \chi (\pi_P^{-1}(y)\cap X\cap B^{2n}_{(0,r)}))= \chi
(\pi_P^{-1}(0)\cap f^{-1}(\epsilon)\cap B^{2n}_{(0,r)})),   $$
o\`u $\epsilon$ est g\én\érique dans $\d C$ et suffisamment proche de
$0$.
Ceci r\ésulte classiquement du premier Lemme d'isotopie de
Thom-Mather,
de la g\én\éricit\é de la transversalit\é en l'origine,
pour des
$(n-i)$-plans affines de $\d C^n$, aux strates d'une
stratification de Whitney du germe $X_0$ et de la connexit\é
du compl\émentaire dans $\d C^n$ d'un ensemble analytique
complexe n\égligeable.
Supposons enfin que $f$ admette $0$ pour singularit\'e
isol\ée.
Dans ce cas $ \chi(\pi_P^{-1}(0)\cap f^{-1}(\epsilon)\cap
B^{2n}_{(0,r)}))   $
est la caract\éristique d'Euler-Poincar\'e de la fibre de Milnor
de $f$ dans  $P^\perp$. C'est-\`a-dire : $1+(-1)^{n-i-1}\mu_{n-i} $,
o\`u $\mu_{n-i}$ est le nombre de Milnor de la section $(n-i)$-plane de
$X_0$ (introduit dans [Te1]). On a donc, dans le cas
o\`u $X_0$ est le germe d'une hypersurface complexe de $\d C^n$ :
$$\widetilde{\sigma_i}(X_0)=1+(-1)^{n-i-1}\mu_{n-i}.$$
 Dans [Te1] il est montr\'e que  l'ind\épendance  de la suite $(\mu_0(X_t), \cdots, \mu_n(X_t))$ relativement
aux para\-m\`e\-tres~$t$ (et donc, dans ce cadre, de la suite  $(\widetilde{\sigma}_0(X_t), \cdots,
 \widetilde{\sigma}_n(X_t))$), pour
une famille analytique $X=(X_t)_{t \in \hbox{\bbbb C}}
\subset \d C^{n+1}$  de germes d'hypersurfaces analytiques  de $\d C^n$
ayant $0$ pour singularit\é  isol\ée  implique la condition
de Whitney, au voisinage de $0$ dans $\d C$,
pour le couple $(X\setminus \d C, \d C)$.
Et d'apr\`es  [Br-Sp], l'implication  r\éciproque est vraie.

\ev
En toute g\'en\'eralit\'e
($X$ analytique complexe dans $\d C^n$ de dimension quelconque
et non plus seulement une hypersurface)
les invariants complexes
$\widetilde{\sigma}_i(X_t)$
ont \'et\'e consid\'er\'e pour la premi\`ere fois par M. Kashiwara
dans [Ka1]
(o\`u les boules sont ouvertes et non ferm\'ees comme
c'est le cas ici) : 
 un invariant $E_{X_0}^0$ y est d\'efini
par r\'ecurrence sur la dimension de $X_0$ \`a l'aide de $\ti \sigma_i$. L'\'etude de cet invariant est reprise par
A. Dubson ([Du1,2]) puis dans
[Br-Du-Ka] o\`u les auteurs en donnent une version
multidimensionnelle $E^k_{X_0}$. Leur d\'efinition est la suivante :

$$ E_{X_0}^k = \sum_{X^{j_0}\subset \bar X_j\setminus X^j, \ \dim(X^j)<\dim(X_0)} E^k_{\bar X^j}\cdot
\ti \sigma_{k+\dim(X^j)+1}(X_0),$$
o\`u $(X^j)_{j\in\{0, \cdots, k\}}$ est une stratification de Whitney de $X_0$, et $X^{j_0}$ la strate contenant $0$.
Les auteurs remarquent ensuite  (cf aussi [Du1,2]) que :
$$ E_{X_0}^k = Eu_{X_0},$$
o\`u $Eu_{X_0}$ est l'obstruction d'Euler locale de $X$ en $0$,
introduite par R. MacPherson dans [MacPh], et :
$$ (-1)^k(E_{X_0}^{\dim(X_0)-k-1}-E_{X_0}^{\dim(X_0)-k})= 
e(\r P^k(X_0),0), $$
o\`u $e(\r P^k(X_0),0)$ est la multiplicit\'e en $0$ de la
vari\'et\'e polaire $\r P^k(X_0)$ de
codimension $k$ de $X_0$ en $0$
(Voir aussi [Me], [L\^e-Te1,$\cdots$,3], [Du2]).
Il est annonc\'e sans preuve dans
[Du1] Proposition 1, [Du2] Th\'eor\`eme II.2.7, page 30,
 et  [Br-Du-Ka],
 que les invariants $\ti \sigma_i(X_y)$ sont constants lorsque
$y$ varie dans une strate d'une stratification de Whitney de $X_0$.
Mais dans [He-Me1], [Na2], [Te2]
 il est prouv\'e que la constance des multiplicit\'es
$e(\r P^k(X_y),y)$ lorsque $y$ varie dans une strate d'une stratification donn\'ee de $X_0$ \'equivaut
\`a la $(b)$-r\'egularit\'e de cette stratification, ce qui donne une preuve,
compte tenu de l'\'egalit\'e ci-dessus reliant les $e(\r P^k(X_y),y)$ et les $\ti \sigma_i(X_y)$,
de la constance des  $\ti \sigma_i(X_y)$ le long de strates de Whitney.

On rassemble ces r\'esultats dans le th\'eor\`eme suivant, o\`u
$e( \Delta^k(X_y),y)$ est la multiplicit\'e en $y$ du discriminant
$\Delta^k(X_y)$ associ\'e \`a $\r P^k(X_y)$.
 \ev
{\bf Th\'eor\`eme 2.12. }
([He-Me1], [Na2], [L\^e-Te3], [Te2]) --- {\sl
Soit $X_0$ un germe en $0$ d'ensemble analytique complexe de $\d C^n$
muni d'une stratification
$(X^j)_{j\in \{0, \cdots , k\}}$.
Les conditions suivantes sont \'equivalentes :
\ev

\n
{(i)} \hskip2mm
 \vtop{\hsize12,95cm \n  
La stratification $(X^j)_{j\in \{0, \cdots , k\}}$ est de Whitney.}

 \vskip1mm
\n
{(ii)} \hskip2mm
 \vtop{\hsize12,95cm \n 
 Les fonctions $y\mapsto e(\r P^k(X_y),y) $ en
restriction aux strates $X^j$ sont constantes.}

 \vskip1mm
\n
{(iii)} \hskip2mm
 \vtop{\hsize12,95cm \n  Les fonctions $y\mapsto e( \Delta^k(X_y),y) $ en
restriction aux strates $X^j$ sont constantes.}

 \vskip1mm
\n
{(iv)} \hskip2mm
 \vtop{\hsize12,95cm \n   Les fonctions $y\mapsto \ti \sigma_i(X_y) $ en restriction aux strates $X^j$ sont constantes.
 }}
\ev
\par\n
Nous donnons une preuve directe de $(i)\Rightarrow (iv)$ au
Corollaire $4.5$, c'est-\`a-dire sans
utiliser $(i)\Rightarrow (ii)$ et le lien entre
$e(\r P^k(X_y),y)$ et $\sigma_i(X_y)$
(la preuve de $(i)\Rightarrow (iv)$ est faite dans
[Br-Sp] dans le cas des hypersurfaces \`a singularit\'e isol\'ee).

\Ev
\centerline{\bf 3. Courbures de Lipschitz-Killing locales et invariants polaires}
\Ev
L'objet de cette section est de prouver le Th\'eor\`eme $3.1$, qui
relie via la suite $\sigma_*(X_0)$,  les
courbures $\Ll_i$
\`a la g\'eom\'etrie des discriminants g\'en\'eraux des projections du
germe $X_0$ sur des
plans de dimension $i, i+1, \ldots, n$. La preuve se fait en deux
\'etapes. Tout d'abord
pour les c\^ones  sous-analytiques  de sommet l'origine
(section 3-a), ce qui permet, par d\'eformation sur le c\^one
tangent, d'\'etablir le th\'eor\`eme  pour les germes
d'ensembles sous-analytiques ferm\'es quelconques  (section 3-b).
\ev
{\bf Th\éor\`eme  3.1. --- }
{\sl Pour tout $i\in\{1, \ldots, n \} $, il
existe des constantes r\'eelles $m^i_i, \ldots, m^n_i$
tels que pour tout
 germe $X_0$ d'ensemble sous-analytique ferm\é de $\d R^n $,
 on ait l'\égalit\é :
$$\L_i^{\ell oc}(X_0)=\D \sum_{j=i}^n m^j_i \cdot \sigma_j(X_0) .$$
Autrement dit il existe une matrice triangulaire sup\'erieure
$(m_i^j)_{1\le i\le n,1\le j\le n }$ telle que :
$$ \pmatrix{\Ll_1 \cr \vdots \cr \Ll_n}=
\pmatrix{m_1^1 & m_1^2  & \ldots & m_1^{n-1}& m_1^n \cr
           0  & m_2^2  & \ldots & m_2^{n-1 }&m_2^n  \cr
                \vdots &        &        & & \vdots \cr
                     0 & 0      & \ldots &0& m_n^n  \cr}\cdot
         \pmatrix{\sigma_1  \cr \vdots \cr \sigma_n} $$
De plus :
 $ m_i^i=1, \ \ m_i^j=
\D {\alpha_j  \over \alpha_{j-i}\cdot
\alpha_i}C_j^i -  \D {\alpha_{j-1}  \over \alpha_{j-1-i}\cdot
\alpha_i}C_{j-1}^i, \ \ \hbox{ si }\ i+1\le j \le n.$}

\Ev
\centerline{\bf 3-a. Le cas conique }
\Ev
Soit $X_0$ un c\^one sous-analytique ferm\'e de sommet l'origine dans $\d R^n$, c'est-\`a-dire qu'existe
$L$ un sous-analytique compact de la sph\`ere unit\'e de $\d R^n$ tel que : $X_0=\d R_+\cdot L$.
On suppose que $(X^j)_{j\in \{0,\cdots, 2k\}}$ est une stratification de Whithney
 de $ X=X_0\cap B_{(0,1)}$ provenant de $L$, c'est-\`a-dire que $(X^j)_{j\in \{k+1, \cdots, 2k\}}$ est une stratification de Whitney de $L$, et que pour 
$j\in \{1, \cdots, k\}$, $X^j=\d R_+^*\cdot X^{j+k}$, $X^0=\{0\}$.
Puisque $X_0$ est un c\^one :
$$ \lkl(X_0)={1\over \alpha_i} \L_i( X).$$
On commence par rappeler le calcul
 fait dans [Br-Ku], section $5$.
Pour cela on reprend les notations de [Br-Ku].
Soit $Y\in \d R^n$ un ensemble sous-analytique compact,
$(Y^j)_{j\in \{1, \cdots, k\}} $ une stratification de 
Whitney de $Y$ et $ \nu  \in \d R^n$. On note :
$$\gamma( Y,\nu  ,x)=1-\chi(B(x,\delta)\cap 
g_\nu^{-1}(g_\nu (x)-\epsilon),$$
pour $0<\epsilon< \!\! <\delta<\!\!<1$, $g_\nu  : Y\to \d R$ o\`u  
$g_\nu(x)=\Vert \nu  -x\Vert^2$
et $x\in Y^j$ est un point critique de  $g_{\nu \vert Y^j}$.
On note encore :
$$\gamma_{tan}(Y, \nu  ,x)=1-\chi(B(x,\delta)\cap
 g_{\nu  \vert Y^j}^{-1}(g_\nu  (x)-\epsilon),$$
pour $0<\epsilon< \!\! <\delta<\!\!<1$
et $x\in Y^j$ un point critique de   $g_{\nu   \vert Y^j}$ et enfin :
$$\gamma_{nor}(Y, \nu  ,x)=1-\chi(B(x,\delta)\cap 
g_{\nu  \vert Y\cap S}^{-1}(g_\nu  (x)-\epsilon),$$
pour $0<\epsilon< \!\! <\delta<\!\!<1$
et $x\in Y^j$ un point critique de   $g_{\nu \vert Y^j\cap S}$, o\`u $S$ est une vari\'et\'e
lisse, de dimension $n-\dim(Y^j)$, transverse \`a $Y^j$ en $x$.
On a alors ([Go-MacPh]) :
$$ \gamma(Y, \nu  ,x) = \gamma_{tan}(Y, \nu  ,x)\cdot \gamma_{nor}(Y, \nu  ,x),$$
et   pour presque tout $\nu  \in \d R^n$ ([Br-Ku], Lemme $3.5$) :
$$\chi(Y)=\sum_{x\in Y} \gamma(Y, \nu  ,x). $$

Notons $\pi:\r N\to X^j$ le fibr\'e normal d'une strate $X^j$ de $X$, $\r N_1$ son fibr\'e normal unitaire et
pour $r\in \d R$, $\r N_r=\{(t,\nu,x)\in \d R_+\times \r N_1\subset \d R_+\times \d R^n\times X^j; \Vert\nu\Vert\le r\}$, enfin soit $\Psi :\r N \to \d R^n$,
l'application d\'efinie par $\Psi ((t,\nu,x))=x+t\cdot\nu$.

On pose : $$ \L(X,X^j)(x,r)=
\int_{t\le r}\int_{\nu\in \pi^{-1}(x)\cap {\cal N}_1} 
\gamma_{nor}(X,\nu,x)\cdot \phi(t,\nu,x)\ \
dtd\nu,$$
o\`u $\phi$ est $Jac(\Psi)$.
On va calculer : $$ \sum_{j=0}^{2k}\int_{x\in X^j} \L(X,X^j)(x,r) \ \ dx.$$
Pour cela remarquons que si $d=\dim(X^j)$ :
$$\phi(t,\nu,x)=\det \pmatrix{ 1&0&0\cr *& t\cdot I_{n-d-1}&0\cr *& *& I_d+t\cdot \hbox{II}(\nu,x)\cr },$$
o\`u $I_d$ est la matrice unit\'e d'ordre $d$ et $ \hbox{II} (\nu,x)=\hbox{II} _{X^j}(\nu,x)$ la matrice
dans une base orthonorm\'ee de $T_xX^j$ de la seconde forme
fondamentale de $X^j$ en $x$ suivant
la direction normale $\nu$. La forme $\hbox{II} (\nu,x)$ est d\'efinie de la fa\c con suivante~:
si $\mu(s)$ et $\tau(s)$ sont deux chemins diff\'erentiables trac\'es sur $X^j$, tels que $\mu(0)=\tau(0)=x$ et si
$\nu(y)$ est un champ de vecteurs normaux \`a $X^j$, tel que $\nu(x)=\nu$ et si d'autre part
$x=x(u_1, \cdots, u_d)$ sont des coordonn\'ees sur $X^j$ telles que $x=x(0, \cdots, 0)$,
$\D {\partial x\over \partial u_i}(0)=\mu'(0)$,  $\D {\partial x\over \partial u_j}(0)=\tau'(0)$ :
$$\hbox{II}(\nu,x)(\mu'(0),\tau'(0))= <D\nu_{(x)}\big(\mu'(0)\big)\vert \tau'(0)>=-< \nu \vert {\partial^2 x\over \partial u_i\partial u_j}(0)>.$$
Les propri\'et\'es classiques de finitude en g\'eom\'etrie
mod\'er\'ee montrent alors que
$ \L(X,X^j)(x,r)$ est un polyn\^ome en $r$ que l'on \'ecrit :
$$ \L(X,X^j)(x,r)=  \sum_{i=n-d}^n \lambda_{n-i}(X,X^j)(x)\cdot \alpha_i  \cdot r^i.$$
En observant que le signe de $\phi(t,\nu,x)$ est $\gamma_{tan }(X,t.\nu,x)$ (cf [Bro-Kup], (5.1.4)) et
$\gamma_{nor}(X,\nu,x)=\gamma_{nor}(X,t.\nu, x)$,
on d\'eduit que :
$$ \int_{x\in X^j} \L(X,X^j)(x,r)\ \ dx
= \sum_{i=n-d}^n \alpha_i\cdot  \int_{x\in X^j} \lambda_{n-i}(X,X^j)(x)\ \ dx \cdot r^i=$$
$$\int_{y \in \R^n} \ \sum_{x\in X^j, \Vert x-y\Vert\le r } \gamma(X,y,x)\ \ dy.$$
Puis en sommant sur toutes les strates et en convenant que  $\lambda_\ell(X,X^j)=0$ si $ \ell>\dim(X^j)$~:
 $$ \sum_{j=0}^{2k}\int_{x\in X^j} \L(X,X^j)(x,r) \ \ dx=
\sum_{i=0}^n \alpha_i \Big(\sum_{j=0}^{2k} \int_{x\in X^j} \lambda_{n-i}(X,X^j)(x)\ \ dx\Big) \cdot r^i=$$
$$ \int_{y \in \R^n} \ \sum_{x\in X, \Vert x-y\Vert\le r } \gamma(X,y,x) \ \ dy=
\int_{y \in \R^n} \  \chi(X\cap B(y,r))   \ \ dy.$$

On conclut de ce calcul que : $$ \Ll_i(X_0)= \D {1\over \alpha_i} \L_i(X) =
\D {1\over \alpha_i} \sum_{j=0}^{2k} \int_{x\in X^j} \lambda_{i}(X,X^j)(x)\ \ dx, $$
Dans cette somme seules les strates $X^j$ de dimension plus grandes que $i$ interviennent effectivement.
De la m\^eme fa\c con, en faisant jouer \`a $L=X\cap S^{n-1}_{(0,1)}$ et \`a sa stratification
$(X^j)_{j\in \{k+1, \cdots, 2k\}}$  le r\^ole que viennent de jouer  $X$ et $(X^j)_{j\in \{0, \cdots, 2k\}}$,
on obtient :
$$\L_i(L) =
\D {1\over \alpha_i} \sum_{j=k+1}^{2k} \int_{x\in X^j} \lambda_{i}(L,X^j)(x)\ \ dx, $$
l\`a encore, seules les strates de dimension plus grandes que $i$ ont une contribution non nulle.
\ev
On   montre maintenant que tout comme $ \L_i(L)$, on peut exprimer
$ \Ll_i(X_0)$ comme une combinaison lin\'eaire \`a coefficients 
universels
des seules quantit\'es $ \D \int_{x\in X^j} \lambda_{i}(L,X^j)(x)$, pour $j\in\{k+1, \cdots, 2k\}$. Ce qui
permet d'exprimer $ \Ll_i(X_0)$ comme combinaison lin\'eaire des $ \L_j(L)$. Ceci est
possible gr\^ace au caract\`ere conique de $X$ et montre
que les invariants introduits dans [Be-Br2] sont des combinaisons
lin\'eaires des $\Ll_i$ (cf [Be-Br2], Remarque $5.4$), dans un premier temps pour les
c\^ones sous-analytiques, puis d'apr\`es la section 3-b
qui suit, dans le cas le plus g\'en\'eral.
\ev
Soit $j\in \{1, \cdots, k\}$, $X^j$ une strate conique de dimension $d$ et $x\in X^j$.
Alors  $\Vert x \Vert\in ]0,1[$ et $r(x)=x/\Vert x\Vert \in X^{j+k}$.
Pour $\nu $ est un vecteur  unitaire normal \`a $X^j$ en $x$, et pour $\theta\in [0,\pi]$,
on note $ \nu_\theta =\sin(\theta).\nu+\cos(\theta).r(x)$,  un vecteur unitaire
normal \`a $X^{j+k}$ en $r(x)$.
Lorsque $\nu $ d\'ecrit $\pi^{-1}(\{x\})\cap \r N_1$ et $\theta $ d\'ecrit $[0,\pi]$,
$\nu_\theta$ d\'ecrit $\pi^{-1}(\{r(x)\})\cap \r N_1 $.
On remarque alors que :
$$ \gamma_{nor }(X,\nu,x) = \gamma_{nor }(X,\nu,u\cdot r(x))= \gamma_{nor }(L,\nu_\theta,r(x)), \ \ \forall
u\in ]0,1[, \ \forall \theta\in ]0,\pi[, )$$
puis que :
$$ \gamma_{nor }(X,\nu,x) = \gamma_{nor }(X,\nu_\theta,r(x)), 
\ \ \forall \theta\in ]0,{\pi\over 2}[, \eqno(a)$$
tandis que du fait de la structure conique de $X$ :
$$\gamma_{nor }(X,\nu_\theta,r(x))=0, 
\ \ \forall \theta\in ]{\pi\over 2}, \pi[. \eqno(b)$$
Fixons $\theta\in [0,\pi]$ et comparons maintenant
$\hbox{II} _{X^j}(\nu, x)$ et $\hbox{II}_{X^{j+k}}(\nu_\theta, r(x))$.
Suivant la direction $r(x)$, $D\nu_{(x)}$ est nulle du fait de la structure conique de $X^j$,
d'autre part si $\nu_\theta$ est un champ de vecteurs unitaires  normaux \`a
$X^{j+k}$, \'ecrit sous la forme :
$$X^{j}\ni y\mapsto \nu_\theta(r(y)) =\sin(\theta).\nu(y)+\cos(\theta).r(y),$$
 $X^j\ni y\mapsto \nu(\Vert x\Vert.y)$ est un champ de vecteurs unitaires normaux \`a
$X^j$ en $\Vert x\Vert.y$, et si $\hat\mu$ est un chemin trac\'e dans
$X^{j+k}$ passant par $r(x)$ en $0$ et tel que $\hat\mu'(0)\in T_{(r(x))}X^{j+k}$ est unitaire,
$\mu=\D\Vert x\Vert \hat\mu$ est un chemin trac\'e dans $X^j$, passant
par $x$ en $0$ et $\mu'(0)\in T_{(x)}X^{j}$ est de norme $\Vert x\Vert$ . On a alors, puisque  :
$$ \nu_\theta\circ \hat \mu=\sin(\theta).\nu\circ\mu+\cos(\theta).r\circ \mu $$
$$ D\nu_{\theta[r(x)]}( \hat\mu'(0)) = \sin(\theta)\cdot D\nu_{(x)}(\mu'(0))+\cos(\theta)\cdot Dr_{(x)}(\mu'(0)) $$
$$ \hskip-3mm^t[\hat\mu'(0)]\cdot\hbox{II}_{X^{j+k}}(\nu_\theta,r(x))=
\sin(\theta)\cdot \Vert x\Vert \cdot^t [\hat\mu'(0)]\cdot\hbox{II}_{X^j}(\nu,x) +\cos(\theta)\cdot
\Vert x\Vert \cdot ^t [\hat\mu'(0)]\cdot\hbox{II}_{S^{n-1}_{(0,\Vert x \Vert)}}(\nu,x)_{\vert T_{x}X^j}  $$
$$ \hbox{II}_{X^{j+k}}(\nu_\theta,r(x))=
\sin(\theta)\cdot \Vert x\Vert
\cdot\hbox{II}_{X^j}(\nu,x) +\cos(\theta)\cdot
I_{d-1}.  $$

On note $A$, $B$ et $C$ les
matrices d'ordre $d-1$   d\'efinie par les relations suivantes :
$$\hbox{II}_{X^j}(\nu, x) =\pmatrix{0 & 0\cr 0 & A \cr},
\hskip0,5cm  C=\Vert x\Vert \cdot A, $$
$$ B=\sin(\theta)\cdot C+
\cos(\theta)\cdot I_{d-1}.$$

On a alors d'une part :
$$ \L(X,X^j)(x,r)=\int_{t\le r}\int_{\nu\in \pi^{-1}(x)\cap {\cal N}_1} \gamma_{nor}(X,\nu,x)\cdot \phi(t,\nu,x)\ \
dtd\nu$$
$$\hskip-0cm = \int_{t\le r}\int_{\nu\in \pi^{-1}(x)\cap {\cal N}_1} \gamma_{nor}(X,\nu,x)  \cdot
\det \pmatrix{ 1&0&0  &0 \cr *& t\cdot I_{n-d-1}&0&0 \cr
*&*&1&0  \cr *& *& *& I_{d-1}+t\cdot A\cr }
dtd\nu$$
$$\hskip0cm = \int_{t\le r}\int_{\nu\in \pi^{-1}(x)\cap {\cal N}_1} \gamma_{nor}(X,\nu,x)  \cdot
\det \pmatrix{ 1&0&0&0 \cr *& t\cdot I_{n-d-1}&0&0 \cr
*&*&1 &0 \cr *& *& *& I_{d-1}+t\cdot \D{C\over \Vert x\Vert}\cr }
dtd\nu. \eqno(1)$$

Et d'autre part, d'apr\`es $(a) $ et $(b) $:
$$ \L(X,X^{j+k})(r(x),r)=
\int_{t\le r}\int_{\nu_\theta\in \pi^{-1}(r(x))\cap {\cal N}_1} \gamma_{nor}(X,\nu_\theta,r(x))
\cdot \phi(t,\nu_\theta,r(x))\ \
dtd\nu$$
$$= \int_{t\le r}\int_{\nu\in \pi^{-1}(x)\cap {\cal N}_1} \gamma_{nor}(X,\nu,x)$$

$$\int_{\theta\in [0,{\pi\over 2}]}  f(\theta)\cdot
\det \pmatrix{ 1&0&0\cr *& t\cdot I_{n-d}&0\cr *& *& I_{d-1}+t\cdot (\sin(\theta)\cdot C+
\cos(\theta)\cdot I_{d-1})\cr }
dtd\nu d\theta, \eqno(2)$$
avec $f(\theta)=\sin^{n-d-1}(\theta)$.
Notons : $$ \det(C+T\cdot I_{d-1})=\sum_{p=0}^{d-1}a_pT^p.$$
D'apr\`es $(1)$ :
$$  \L(X,X^j)(x,r)=\int_{t\le r}\int_{\nu\in \pi^{-1}(x)\cap {\cal N}_1} \gamma_{nor}(X,\nu,x)\cdot
t^{n-d-1}\cdot \det(I_{d-1}+\D{t\over \Vert x \Vert }\cdot C)\ \
dtd\nu $$
$$ \L(X,X^j)(x,r)= \sum_{p=0}^{d-1}\ \
\int_{t\le r}\int_{\nu\in \pi^{-1}(x)\cap {\cal N}_1} {\gamma_{nor}(X,\nu,x)\over \Vert x\Vert^{d-1-p}}\cdot a_p\cdot
t^{n-2-p} \ \ dtd\nu  $$
$$ \L(X,X^j)(x,r)= \sum_{p=0}^{d-1}\ \ {a_p\over n-p-1}\cdot
r^{n-1-p}
\int_{\nu\in \pi^{-1}(x)\cap {\cal N}_1} {\gamma_{nor}(X,\nu,x)\over \Vert x\Vert^{d-1-p}} \ \ d\nu  $$
$$ \int_{x\in X^j} \L(X,X^j)(x,r) \ \ dx
= \Gamma(X)\cdot \sum_{i=n-d}^{n-1}\ \ a_{n-i-1} \cdot  { K_{n,i} \over i} \cdot r^i
 \eqno(3)$$
o\`u
$$\Gamma(X)=\int_{r(x)\in X^{j+k}} \int_{\nu\in \pi^{-1}(x)\cap {\cal N}_1} \gamma_{nor}(X,\nu,r(x)) \ \
d\nu  dr(x) $$ et o\`u
$K_{n,i}$ est une constante ne d\'ependant que de $n$ et $i$ et qui provient de :
$$ \int_{x\in X^j} \int_{\nu\in \pi^{-1}(x)\cap {\cal N}_1} {\gamma_{nor}(X,\nu,x)\over \Vert x\Vert^{d+i-n}} \ \ d\nu dx
=K_{n,i}\cdot \Gamma(X).$$
Notons que l'ind\'ependance de $K_{n,i}$ relativement \`a $d$ vient du fait que $X^j$ est de dimension $d$ et que
l'on calcule $ \D\int_{x\in X^j} \int_{\nu\in \pi^{-1}(x)\cap {\cal N}_1} {\gamma_{nor}(X,\nu,x)\over \Vert x\Vert^{d+i-n}} \ \ d\nu dx $  gr\^ace au changement de variables $X^{j+k}\times ]0,1[\ni (r(x),t)\mapsto t\cdot r(x)\in X^j $.
L'\'egalit\'e $(3)$ ensuite donne :
$$ \int_{x\in X^j} \L(X,X^j)(x,r) \ \ dx
= \Gamma(X)\cdot \sum_{i=n-d}^{n-1}\ \ a_{n-i-1} \cdot  { K_{n,i} \over i} \cdot r^i
=\sum_{i=n-d}^{n} \ \int_{x\in X^j} \ \lambda_{n-i}(X,X^j)\cdot\alpha_i\cdot  r^i,$$
et donc :
$$\int_{x\in X^j} \ \lambda_0(X,X^j)\ \  dx=0,$$
 $$ \int_{x\in X^j} \ \lambda_{p+1}(X,X^j)\ \  dx= \Gamma(X)\cdot  { K_{n,n-p-1} \over (n-p-1)\cdot \alpha_{n-p-1}} \cdot  a_p  ,
\ \ p\in\{ 0, \cdots, d-1\}. \eqno (4)$$

Maintenant, d'apr\`es $(2)$ :

$$ \L(X,X^{j+k})(r(x),r)= \int_{t\le r}\int_{\nu\in \pi^{-1}(x)\cap {\cal N}_1} \gamma_{nor}(X,\nu,x)$$
$$\int_{\theta\in [0,{\pi\over 2}]}  f(\theta)\cdot \sin^{d-1}(\theta)\cdot t^{n-1}\cdot
\det \Big( C+({1\over t \sin(\theta)}+{\cos(\theta)\over \sin(\theta)})\cdot I_{d-1} \Big) d\theta d\nu dt.$$

Posons : $$ \alpha=\alpha(\theta)=1/\sin(\theta), \ \ \beta=\beta(\theta)=\cos(\theta)/\sin(\theta) \ \ \hbox{ et }
g(\theta)=\sin^{d-1}(\theta)\cdot f(\theta)= \sin^{n-2}(\theta).$$
On obtient :

$$ \L(X,X^{j+k})(r(x),r)= \int_{t\le r}\int_{\nu\in \pi^{-1}(x)\cap {\cal N}_1}
\gamma_{nor}(X,\nu,x)$$
$$ \int_{\theta\in [0,{\pi\over 2}]}  g(\theta)\cdot t^{n-1}\cdot
\sum_{p=0}^{d-1} a_p\cdot \big({\alpha \over t}+\beta\big)^p d\theta d\nu dt$$
$$ \L(X,X^{j+k})(r(x),r)= \int_{t\le r}\int_{\nu\in \pi^{-1}(x)\cap {\cal N}_1}
\gamma_{nor}(X,\nu,x)\  d\nu$$
$$ \sum_{q=0}^{d-1} t^{n-q-1} \ \int_{\theta\in [0,{\pi\over 2}]}  g(\theta) \cdot \alpha^q \sum_{p=q}^{d-1}
a_p\cdot C_p^q\cdot \beta^{p-q} \ d\theta  dt$$
$$ \int_{r(x)\in X^j} \L(X,X^{j+k})(r(x),r) \ \ dr(x)= \Gamma(X)\
\sum_{i=n-d+1}^{n} {r^i\over i} \    \sum_{p=n-i}^{d-1}
a_p\cdot \delta(i,p),\eqno(5)$$
o\`u :
$$ \delta(i,p)=\int_{\theta\in [0,{\pi\over 2}]}  g(\theta)\cdot \alpha^{n-i}\cdot C_p^{n-i}\cdot \beta^{p-n+i}.$$

L'\'egalit\'e $(5)$ donne :
$$\sum_{i=n-d+1}^{n} \ \int_{r(x)\in X^j} \ \lambda_{n-i}(X,X^{j+k}) \ \ dr(x)\cdot\alpha_i\cdot  r^i=
 \Gamma(X)\
\sum_{i=n-d+1}^{n} {r^i\over i} \    \sum_{p=n-i}^{d-1}
a_p\cdot \delta(i,p)$$
soit  :
$$ \int_{r(x)\in X^j} \ \lambda_{n-i}(X,X^{j+k})\ \ dr(x)=  { \Gamma(X)\over i\cdot \alpha_i} \    \sum_{p=n-i}^{d-1}
a_p\cdot \delta(i,p), $$
et d'apr\`es $(4)$ :

$$\hskip0mm \int_{r(x)\in X^{j+k}} \ \lambda_{n-i}(X,X^{j+k})\ \ dr(x)$$
$$={1\over i\cdot \alpha_i}
 \    \sum_{p=n-i}^{d-1} {(n-p-1)\cdot \alpha_{n-p-1}\over K_{n,n-p-1}}\cdot \delta(i,p)\int_{x\in X^j} \ \lambda_{p+1}(X,X^j)\ \  dx $$
$$\hskip0mm \int_{r(x)\in X^{j+k}} \ \lambda_{\ell}(X,X^{j+k})\ \ dr(x)$$
$$=\ \sum_{p=\ell}^{d-1} {(n-p-1)\cdot \alpha_{n-p-1}\over (n-\ell)\cdot \alpha_{n-\ell}\cdot K_{n,n-p-1}}\cdot
\delta(n-\ell,p)\int_{x\in X^j} \
\lambda_{p+1}(X,X^j)\ \  dx $$
Soit en notant $\D M_n(\ell, p)={(n-p-1)\cdot \alpha_{n-p-1}\over (n-\ell)\cdot \alpha_{n-\ell}\cdot K_{n,n-p-1}}\cdot
\delta(n-\ell,p)$ :
$$\int_{r(x)\in X^{j+k}} \ \lambda_{\ell}(X,X^{j+k})\ \  dr(x)=
 \    \sum_{p=\ell}^{d-1} M_n(\ell, p)\cdot \int_{x\in X^j} \
\lambda_{p+1}(X,X^j)\ \ dx , \ \  \ell\in \{0, \cdots, d-1 \} . \eqno (6)$$
ou encore si :
$$Q=(Q_n(i,j))_{i,j\in \{0, \cdots, d-1\}}  =              \pmatrix{M_n(0,0)& M_n(0,1)& \cdots &M_n(0,d-1)\cr
                                  0&M_n(1,1) & \cdots & M_n(1,d-1)\cr
                                  \vdots & & & \vdots \cr
                                 0& \cdots  & 0 & M_n(d-1, d-1)}^{-1}$$

$$\int_{r(x)\in X^{j+k}} \ \lambda_{\ell+1}(X,X^{j}) \ \ dr(x)$$
$$=
 \    \sum_{p=\ell}^{d-1} Q_n(\ell, p)\cdot \int_{x\in X^j} \
\lambda_{p}(X,X^{j+k})\ \ dx , \ \  \ell\in \{0, \cdots, d-1 \} . \eqno (7)$$

Enfin, le m\^eme calcul montre que :
$$ \int_{r(x)\in X^{j+k}} \ \lambda_{n-i}(L,X^{j+k})\ \ dr(x)=  { \Gamma(X)\over i\cdot \alpha_i} \
\sum_{p=n-i}^{d-1} a_p \cdot \ti\delta(i,p), $$

o\`u :
$$ \ti \delta(i,p)=\int_{\theta\in [0,\pi]}  g(\theta)\cdot \alpha^{n-i}\cdot C_p^{n-i}\cdot \beta^{p-n+i}.$$
Si on d\'efinit $\ti M_n(i,j)$ de la m\^eme fa\c con que $M_n(i,j)$, mais en rempla\c cant
$\delta(i,p)$ par $\ti \delta(i,p)$, puis si $\ti Q$ d\'esigne l'inverse de la matrice triangulaire
$(\ti M_n(i,j))_{j\ge i}$, on obtient finalement :
$$ \pmatrix{\D \int_{x\in X^j} \lambda_1(X, X^j) \ dx \cr \vdots \cr \D\int_{x\in X^j}  \lambda_d(X, X^j) \ dx }=Q\cdot
\pmatrix{ \D\int_{x\in X^{j+k}}   \lambda_0(X, X^{j+k})\ dx \cr \vdots \cr
\D\int_{x\in X^{j+k}}    \lambda_{d-1}(X, X^{j+k})\ dx \cr }$$
$$=\ti Q\cdot
\pmatrix{ \D\int_{x\in X^{j+k}}   \lambda_0(L, X^{j+k}) \ dx \cr \vdots \cr
\D\int_{x\in X^{j+k}}   \lambda_{d-1}(L, X^{j+k})\ dx \cr }$$

Notons que la matrice $\ti Q=\ti Q(X^j)$, dont l'ordre est la dimension de $X^j$,
 ne d\'epend de $X^j$, ou plut\^ot de $\dim(X^j)$, que par sa taille;
les coefficients de $\ti Q$ ne d\'ependent en effet pas de la dimension de $X^j$.
Pour \^etre plus pr\'ecis,  si pour
$j,j'\le k$, $\dim(X^j)\le \dim(X^{j'})$, on a :
$$\ti Q(X^{j'})=\pmatrix{\ti Q(X^j)& *\cr * &* \cr}.$$

De plus  comme $\lambda_i(X,X^0)=0$, pour $i\not=0$, on a
 pour $i\ge 1$ :
 $$ \Ll_i(X_0)= \D{1\over \alpha_i}\L_i(X)= {1\over \alpha_i}\sum_{j=0}^{2k} \int_{x\in X^j}\lambda_i(X,X^j)\ dx=
{1\over \alpha_i}\sum_{j=1}^{2k} \int_{x\in X^j}\lambda_i(X,X^j)\ dx,$$
 et  en notant
$\ti Q_i(X^j)$ la $i$\`eme ligne de $\ti Q(X^j)$, on obtient  :

$$ \Ll_i(X_0)= {1\over \alpha_i}\sum_{j=1}^{2k}\int_{x\in X^j} \lambda_i(X,X^j) \ dx  $$
$$={1\over \alpha_i}\Big(\sum_{j=1}^{k}\int_{x\in X^j} \lambda_i(X,X^j) \ dx+
\sum_{j=1}^{k}\int_{x\in X^{j+k}} \lambda_i(X,X^{j+k}) \ dx \Big)$$
$$={1\over \alpha_i} \sum_{j=1}^{2k} \Big(\ti Q_i(X^j)+ [Q^{-1}(X^j)\ti Q(X^j)]_i \Big)
\pmatrix{ \D\int_{x\in X^{j+k}} \lambda_0(L,X^{j+k}) \ dx \cr \vdots\cr
\D\int_{x\in X^{j+k}} \lambda_{\dim(X^j)-1}(L,X^{j+k}) \ dx \cr }$$
Soit  :
$$\Ll_i(X_0)=   \hat Q_i
\pmatrix{  \L_0(L)  \cr \vdots\cr
 \L_{\ell -1}(L) \cr }
, $$
avec $ \dim(X)=\ell=\dim(X^{\hat j})$, $\hat Q$ la matrice d'ordre $ \ell$ dont la $i$\`eme ligne $\hat Q_i$ est
$\D{1\over \alpha_i}[\ti Q(X^{\hat j})+Q^{-1}(X^{\hat j}) \ti Q(X^{\hat j})]_i $. Notons que $\hat Q$ est
triangulaire sup\'erieure et que :
$$\pmatrix{ \Ll_1(X_0)\cr \vdots \cr
\Ll_\ell (X_0)\cr }
= \hat Q
\pmatrix{  \L_0(L)  \cr \vdots\cr
 \L_{\ell-1}(L) \cr }. \eqno(8)$$
Cette \'egalit\'e est \`a comparer avec celle donn\'ee
dans [Be-Br2], Remarque $5.4$.

Exprimons maintenant $ \sigma_j(X_0) $ comme combinaison lin\'eaire \`a coefficients universels des invariants
$\L_i(L)  $.
Par d\'efinition :
$$ \sigma_j(X_0) = \int_{P\in {\cal E}^j_X } {1\over j\cdot \alpha_j}
\int_{\ell\in S^P_{(0,1)}}
{1\over 2}\cdot (\chi_-^\ell+\chi_+^\ell) \ \  d\ell\  d\gamma_{j,n}(P),$$
o\`u $S^P_{(0,1)}$ est la sph\`ere unit\'e de $P$,
$\pi_{P,\ell}$ la projection
orthogonale de $\ell\oplus P^\perp$ sur $\ell$ (on identifie la direction $\ell$ et la droite qu'elle supporte)
et  $0<\epsilon=\epsilon(X,\ell)<<1$,
$\chi^\ell_-=\chi(X\cap \pi^{-1}_{P,\ell}(-\epsilon) )$, $\chi^\ell_+=\chi(X\cap \pi^{-1}_{P,\ell}(\epsilon))$.
Remarquons que du fait que $X_0$ est un c\^one, $\chi^\ell_-=
\chi(L \cap \pi^{-1}_{P,\ell} (]-\infty,-\epsilon])) $ et  $\chi^\ell_+= \chi(L \cap \pi^{-1}_{P,\ell} ([\epsilon,+\infty [)) $.
Mais d'autre part, le premier  lemme d'isotopie de Thom-Mather montre que :
$$ \chi(L \cap \pi^{-1}_{P,\ell} (]-\infty,-\epsilon]))= \chi(L \cap \pi^{-1}_{P,\ell} (]-\infty,0])),$$
$$\chi(L \cap \pi^{-1}_{P,\ell} ([\epsilon,+\infty [))= \chi(L \cap \pi^{-1}_{P,\ell} ([0,+\infty [)),  $$
puisque pour $\epsilon $ suffisamment proche de $0$, $\pi_{P,\ell}^{-1}(\epsilon)$ est transverse aux strates de
$L$ (Proposition $2.6$).
On obtient alors :
$$ \sigma_j(X_0) = {1\over j\cdot \alpha_j} \int_{P\in {\cal E}^j_X }
\int_{\ell\in S^P_{(0,1)}}
{1\over 2}\cdot (\chi(L\cap P^\perp)+\chi(L\cap \ell\oplus P^\perp)  \ \  d\ell\  d\gamma_{j,n}(P).$$
Il existe par cons\'equent deux  constantes $a(j,n)$ et $a'(j,n)$
ne d\'ependant que de $j$ et $n$ telles que :
 $$ \sigma_j(X_0) = a(j,n) \int_{P\in {\cal E}^j_X }  \chi(L\cap P^\perp) \ \ d\gamma_{j,n}(P) +
a'(j,n)\int_{H\in {\cal E}'^{n-j+1}_X }
  \chi(L\cap H)  \ \  d\gamma_{j,n}(P),$$
avec ${\cal E}'^{n-j+1}_X $ un sous-analytique dense de $G(n-j+1,n) $.
En notant $\r O(S^{n-1})$ le groupe des isom\'etries de $S^{n-1}_{(0,1)}$, il existe un sous-analytique dense
$O_X^j $ dans $\r O(S^{n-1}) $ tel que :
$$ \sigma_j(X_0) = \ti a(j,n) \int_{\sigma \in O^j_X }  \chi(L\cap \sigma\cdot K) \ \ d\sigma +
\ti a'(j,n)\int_{\sigma \in O^j_X }
  \chi(L\cap \sigma\cdot H)  \ \  d\sigma,$$
o\`u $K$ est la trace dans $S^{n-1}_{(0,1)}$ d'un $(n-j)$ plan vectoriel fix\'e et $H$ celle d'un $(n-j+1)$ plan
vectoriel fix\'e.
En utilisant $\chi=\Lambda_0$, puis le fait que $\Lambda_0$ (comme les $\Lambda_p$)
est combinaison lin\'eaire des
courbures    sph\'eriques relatives $\ti \Lambda_i$ ([Be-Br2],
Th\'eor\`eme $1.1$, ou formule de Gauss-Bonnet
sph\'erique, [Be-Br2], Th\'eor\`eme $1.2$)
introduites dans [Be-Br2], on obtient en notant $\ell=\dim(X)$ :

$$ \sigma_j(X_0) = \ti a(j,n) \sum_{i=0}^{\ell-1} \int_{\sigma \in O^j_X }  b(i)\cdot \ti \L_i(L\cap \sigma\cdot K) \ \ d\sigma $$
$$+  \ti a'(j,n)\sum_{i=0}^{\ell-1} \int_{\sigma \in O^j_X }
 b(i)\cdot  \ti\L_i (L\cap \sigma\cdot H)  \ \  d\sigma.$$
La formule cin\'ematique ([Be-Br2], Th\'eor\`eme $4.4$) pour les courbures sph\'eriques relatives donne ensuite :
$$ \sigma_j(X_0) =  \sum_{i=j-1}^{\ell-1} c(i,j,n)\cdot \ti \L_i(L).$$
En utilisant encore que $\ti \L_i$ est combinaison lin\'eaire des $\L_p$, pour $p\ge i$, on en d\'eduit :
$$ \sigma_j(X_0) =  \sum_{i=j-1}^{\ell-1} \hat c(i,j,n)\cdot \L_i(L).$$
En notant $R$ la matrice triangulaire sup\'erieure dont les coefficients sont les
$\hat c(i,j,n), j-1\le i\le \ell-1 $ :
$$ \pmatrix{ \sigma_1(X_0) \cr \vdots \cr \sigma_\ell(X_0)} =R \pmatrix{ \L_0(L) \cr \vdots \cr \L_{\ell-1}(L)}.
\eqno(9)$$
Les \'egalit\'es $(8)$ et $(9)$ donne pour finir le Th\'eor\`eme $3.1$ pour les c\^ones sous-analytiques
 (la r\'egularit\'e des matrices $\hat Q$ et $R$ est \'etablie en appendice par un calcul explicite pour $X$ un c\^one
poly\'edral). \f

\Ev
\centerline{\bf 3-b. Le cas g\'en\'eral }
\Ev
Dans cette section $X_0$ est un germe d'ensemble
sous-analytique ferm\'e, contenant l'origine. Contrairement
\`a la section pr\'ec\'edente, on ne suppose plus ici que $X_0$
est un c\^one. On  montre le Th\'eor\`eme
$3.1$ en toute g\'en\'eralit\'e, en d\'eformant $X_0$ sur son
c\^one tangent. On se donne $(X^j)$ une stratification de
Whitney de $X$.
Soit $X$ un repr\'esentant de $X_0$, notons $X_r$
le sous-analytique compact $ \D{1\over r}\cdot(X\cap B_{(0,r)})$,
$L_r=\D{1\over r}\cdot(X\cap S_{(0,r)})$,
$C_rX=\d R_+\cdot L_r$. La famille $(C_rX)_{r\in ]0,1]}$ est la d\'eformation de $X_0$ sur son c\^one tangent
$C_0X$.
\ev
{\bf Lemme 3.2. --- } {\sl Avec les notations ci-dessus, on a, pour $P$ dans un ouvert sous-analytique dense de
$\bar G(n-i,n) $ :}
$$ \lim_{r\to 0}\chi(X_r\cap P)=\lim_{r\to 0}\chi(C_rX\cap P).$$
\ev
{\bf Preuve. } Commen\c cons par prouver ce lemme pour $P$ un hyperplan de $\d R^n$, ie pour
$i=1$.
Soit donc $P$ un  hyperplan affine de $\bar G(n-1,n)$ ne passant pas par $0$, et $E$ le demi-espace ne contenant pas $0$
qu'il d\'efinit. L'application $R:E\to S^{n-1}_{(0,1)}$, d\'efinie par $R(x)=x/\Vert x\Vert $ est une hom\'eomorphisme
de $P\cap B_{(0,1)}$ sur $S^{n-1}_{(0,1)}\cap E$ qui envoie $ C_rX\cap P$ sur $L_r\cap E$. Il suffit alors de prouver
que $\D \lim_{r\to 0}\chi(X_r\cap P)=\lim_{r\to 0}\chi(L_r \cap E)$.
Soit $x\in P\cap S^{n-1}_{(0,1)}$.
L'angle entre $P$ et $T_x S^{n-1}_{(0,1)}$,  est plus petit que $\D{\pi\over 2}-C$, pour
$C>0$ une certaine constante, puisque $0\not\in P$.
On peut d\'eformer $F_0=S^{n-1}_{(0,1)}\cap E$ sur $F_1=P\cap B_{(0,1)}$ dans une famille telle que
l'angle entre $T_aF_t $ et $T_{R(a)}F_0$ soit major\'e par $\D {\pi\over 2}-C$, quel que soit $a\in E\cap B_{(0,1)}$.
On note $X_r^j$ la strate de $1/r\cdot X$ provenant de $X^j$.
Soit alors $r$ suffisamment petit pour que~:
$\forall x\in X_r, \ \exists \ell_x \ C_0 X\cap S^{n-1}_{(0,1)},
\ \Vert R(x)-\ell_x\Vert \le C/3$ et l'angle entre
$\ell_x $ et $T_xX^j_r$ est plus petit que $C/3$ (lorsque $x\in x_r^j$).
Quelle que soit la strate $X_r^j$ et $x\in X_r^j$, quel que soit $t\in [0,1]$,
l'angle $ \angle(T_xF_t,T_xX_r^j)$  v\'erifie :
 $$\hskip-1mm \angle(T_xF_t,T_xX_r^j)\ge \angle(T_{R(x)}F_0,\ell_x)-\angle(T_{R(x)}F_0, T_xF_t)-\angle(\ell_x,T_xX_r^j)$$
$$\ge (\D{\pi\over 2}- C/3)- (\D{\pi\over 2}- C)- C/3=C/3>0. $$
De sorte que la d\'eformation $(F_t)_{t\in [0,1]}$ est transverse aux strates de Whitney $(X_r^j)$ de $X_r$.
Le premier lemme d'isotopie de Thom-Mather assure alors que 
$F_0\cap X_r$ et $F_1\cap X_r$ sont hom\'eomorphes.
Lorsque $P$ n'est pas un hyperplan, on place $P$ dans 
$\pi_V^{-1}(\d R\cdot \pi_V(P))$, o\`u $V$ est le $i$-plan 
vectoriel de $\d R^n$ orthogonal \`a $P$. Comme $\dim(X)\ge i$, 
le raisonnement p\'ec\'edent s'applique pour presque toutes les droites
$\d R\cdot \pi_V(P)$ de $V$.
  \f
\ev
{\bf Lemme 3.3. --- }
 {\sl Avec les notations pr\'ec\'edentes et celles du Th\'eor\`eme $2.8$ :
$$ \lim_{r\to 0} \sigma_j (C_r X)=\sigma_j (X_0).$$}

{\bf Preuve.}
 Soit $P$ un $j$ plan vectoriel de $\d R^n$ qui coupe
$B_{(0,1)}$ et $y\in P$, par le Lemme $3.2$, pour $r$ suffisamment petit on a : $\chi(C_rX\cap (y+P^\perp))=\chi(X_r\cap (y+P^\perp))$. Par d\'efinition des constantes
$ \chi_1^P(X_0), \cdots , \chi_{n_P}^P(X_0)$ et des domaines
$ \r K_1^P(X_0), \cdots , \r K_{n_P}^P(X_0)$ (Th\'eor\`eme $2.8$),
on a :
$$\lim_{r\to 0}  \Theta_j(\{ y\in P\cap B_{(0,1)};
\chi(X_r\cap (y+P^\perp))=\chi_k^P(X_0) \})
=\Theta_j(\r K_k^P(X_0)),$$
 Comme les caract\'etristiques d'Euler-Poincar\'e de toutes les sections planes de $X_r$ et de $C_rX$ sont uniform\'ement
born\'ees relativement \`a $r$, le th\'eor\`eme de convergence
domin\'ee et le Lemme $3.2$ donnent :
$$ \lim_{r\to 0} \sigma_j (C_rX)=
\lim_{r\to 0} \int_{P\in G(j,n)}\ \sum_{k=1}^{n_P^r}\chi_k^P(C_rX)\cdot
\Theta_j \big(\r K(C_rX)\big) \ d\gamma_{j,n}(P)$$
$$=\int_{P\in G(j,n)}\ \sum_{k=1}^{n_P}\chi_k^P(X_0)\cdot
\Theta_j \big(\r K(X_0)\big) \ d\gamma_{j,n}(P)=\sigma_j(X_0).\f$$

\ev
{\bf Preuve du Th\'eor\`eme 3.1.}
 Comme ci-dessus on utilise le fait que
les familles $(C_rX)_{r\in [0,1]}$ et $(X_r)_{r\in [0,1]}$ sont
des sous-analytiques compacts et donc
que les carac\-t\'eristiques d'Euler-Poincar\'e de toutes les
sections planes de $X_r$ et de $C_rX$ sont uniform\'ement
born\'ees relativement \`a $r$ pour appliquer le th\'eor\`eme de
convergence domin\'ee.
 On a ainsi, par le Lemme $3.2$ :
$$ \Ll_i(X_0)=
{1\over \alpha_i}\lim_{r\to 0}\int_{\bar P\in \bar G(n-i,n)} \ 
\chi(X_r\cap \bar P)\ d\bar\gamma_{n-i,n}(\bar P)$$
$$=
{1\over \alpha_i}\lim_{r\to 0}\int_{\bar P\in \bar G(n-i,n)} \ 
\chi(C_rX\cap \bar P)\ d\bar \gamma_{n-i,n}(\bar P) = \lim_{r\to 0} \Ll_i(C_rX).$$
Mais le Th\'eor\`eme $3.1$ \'etant vrai pour les c\^ones :
$$\Ll_i(X_0)= \lim_{r\to 0} \sum_{j=i}^n m_i^j\cdot \sigma_j(C_rX).$$
Ce qui donne par le Lemme $3.3$ l'\'egalit\'e annonc\'ee du Th\`eor\`eme $3.1$ :
$$\Ll_i(X_0)=  \lim_{r\to 0} \sum_{j=i}^n m_i^j\cdot \sigma_j(C_rX) = \sum_{j=i}^n m_i^j\cdot \sigma_j(X_0).\f $$

\Ev
\centerline{\bf 4. Conditions de r\'egularit\'e et invariants locaux}
\Ev

Lorsque $X$ est un ensemble sous-analytique, les fonctions
$\sigma_j$ et $\Ll_j$ sont des fonctions Log-analytiques (voir
par exemple [Li-Ro], [Co-Li-Ro]). Nous nous int\'eressons ici
\`a la r\'egularit\'e
de ces fonctions.
Nous prouvons dans cette section la continuit\'e des invariants
$\sigma_j$ le long de strates d'une stratification de Verdier
d'un ensemble sous-analytique ferm\'e. Il s'agit de montrer que
le long (de la projection) d'une telle strate, les discriminants
des projections
sur des $i$-plans ont un bon comportement du point de vue de la
densit\'e des domaines qu'ils bordent~: l'\'equis\'ecabilit\'e.
Dans le cadre analytique complexe, l'\'equis\'ecabilit\'e des
discriminants le long d'une strate est la constance de leur
multiplicit\'e locale, condition que l'on sait \^etre \'equivalente
\`a la $(w)$, la $(b^*)$, la $(b)$-r\'egularit\'e et \`a la constance des
caract\'eristiques \'evanescentes (cf Section $2$-$b$).
Nous montrons au passage
(Corollaire $4.5$) l'implication $(i)\Rightarrow (iv)$ du
Th\éor\`eme $2.12$, c'est-\`a-dire
la constance des caract\'eristiques \évanescentes le long de
strates de Whitney d'un ensemble analytique complexe
([Br-Sp] dans le cas des familles analytiques
d'hypersurfaces \à singularit\'e isol\ée,
c'est-\à-dire pour les nombres de Milnor des sections planes).
Bien qu'annonc\'ee par ailleurs, il s'agit \`a notre connaissance 
de la seule preuve 
dans le cas g\'en\'eral de cette implication 
qui n'utilise pas le calcul qui lie la multiplicit\'e 
des vari\'et\'es polaires aux caract\'eristiques \'evanescentes.

\ev 

Rappelons que si $E$ est une condition de  r\'egularit\'e portant 
sur un couple
$(Y,Z)$ de sous-vari\'et\'es de $\d R^n$, dans  [L\^e-Tei 2]
et [Tr2] est d\'efinie la 
$E_{\hbox{\t cod} (\ell)}$-r\'egularit\'e ($0 \leq \ell \leq
\hbox{cod}(Y)$) de la fa\c con suivante :  le couple $(Y,Z)$
est dit $E_{\hbox{\t cod} (\ell)}$-r\'egulier en $y \in Y$,
s'il existe  un ouvert dense ${\cal M}_y$ dans $\{ \Pi \in
G(n-\ell,n), \T_yY \subset
\Pi \}$ pour lequel on ait, quelle que soit la 
sous-vari\'et\'e $\cal W$ de $\d R^n$ telle que  $Y \subset
{\cal W} $ dans un voisinage de $y$, l'implication :
 $\T_y{\cal W} \in \r M_y$  $ \Longrightarrow$ $Z$ et ${\cal W}$ 
sont transverses dans un
voisinage de $ y$ et $(Y, {\cal W} \cap Z) $ est
$E$-r\'egulier en $y$. Lorsque $(Y,Z)$
 est $E_{\hbox{\t cod} (\ell)}$-r\'egulier en $y$ pour tout
$\ell$, $0 \leq \ell \leq$ cod $(Y)$, on dit que $(Y,Z)$ est
$(E^*)$-r\'egulier en $y$. On se reportera par exemple \`a 
[Na1], [Na-Tr], [Or]
et [Or-Tr] pour une \'etude compar\'ee 
de la $(r^*)$, de la $(b^*)$ et de la
$(w^*)$-r\'egularit\'e. 
Notons que la condition $(b^*)$ est strictement plus faible que
 la condition $(w)$, m\^eme en  alg\'ebrique.
 En effet,
d'apr\`es [Na-Tr], pour la cat\'egorie 
sous-analytique 
$(w) \Longrightarrow (b^*)$, puis 
 d'une part d'apr\`es [Na-Tr] Corollaire 3.13, si le couple de 
strates $(Y,Z)$ est $(b)$-r\'egulier et si $\dim(Y)=1$, ce couple
  est aussi $(b^*)$-r\'egulier, et  d'autre part dans [Br-Tr] et [Tr1]
sont donn\'es des exemples de couples $(Y,Z)$ qui sont $(b)$-r\'eguliers
mais non $(w)$-r\'eguliers et o\`u $\dim(Y)=1$.
 \ev
 Nous utilisons dans cette partie 
la condition $(*)$ suivante qui porte sur une projection 
$\pi_P$, avec $P\in G(i,n)$ et la stratification $(X^j)_{j\in 
\{0, \cdots, k\}}$ de $X$ que l'on se donne :

\ev
\hskip5mm\vtop{\hsize11cm \n \sl La projection
$\pi_P:\d R^n\to P$ est une projection pour laquelle existe
 un voisinage ou\-vert $\r U(=\r U_P)$ de l'origine
tel que quel que soit la strate $ X^j$ adh\'erente \`a $X^0$:
$$\r U  \cap \r P_{X^j}(P)\cap  (\pi_P^{-1}(\pi_P(X^0))\setminus X^0)
 = \emptyset \eqno(*)$$}

 En nous appuyant sur le Th\'eor\`eme 4.10, nous 
remarquons dans la proposition qui suit 
que la condition $(*)$ est une condition qui a lieu 
pour des projections $\pi_P$ g\'en\'eriques, d\`es lors
que
la stratification $(X^j)_{j\in \{0, \cdots, k\}}$ est 
$(w)$-r\'eguli\`ere.

\ev
{\bf Proposition 4.1. ---}
{\sl Soit $X$ un ensemble sous-analytique
ferm\'e de $\d R^n$, contenant l'origine, muni d'une stratification
$(w)$-r\'eguli\`ere $(X^j)_{j\in \{0, \cdots, k\}}$ et soit 
$X^0$ la strate contenant $0$. Soit $i$
un entier.
 Il existe un ensemble sous-analytique dense $\r G^i_X$ dans $G(i,n)$, 
tel que pour tout $P\in \r G^i_X$,
 $\pi_P:\d R^n\to P\in G(i,n)$ est une projection pour laquelle existe
 un voisinage ouvert $\r U(=\r U_P)$ de l'origine
tel que quel que soit la strate $ X^j$ adh\'erente \`a $X^0$:
$$\r U  \cap \r P_{X^j}(P)\cap  (\pi_P^{-1}(\pi_P(X^0))\setminus X^0) = \emptyset \eqno(*)$$
}
\ev
{\bf Preuve. } Il s'agit d'une cons\'equence imm\'ediate du 
Th\'eor\`eme 4.10.(iii).\f
\ev 
{\bf Remarque. }
Dans la Proposition 4.1, 
on peut se contenter de supposer que la stratification 
$(X^j)_{j\in \{0, \cdots, k\}}$ est seulement 
$(a^*)$-r\'eguli\`ere,  condition ad hoc pour 
que $(*)$ ait g\'en\'eriquement lieu. En effet,
si $P\in G(i,n)$ est tel que la propri\'et\'e $(*)$ n'a
pas lieu, il existe
une strate $X^j$ contenant dans son adh\'erence $X^0$, une suite
$(x_n)_{n\in \hbox{\bbbb N}}$ de limite $0$ et contenue dans $ P_{X^j}(P)\cap  (\pi_P^{-1}(X^0)\setminus X^0)$.
Soit $q=\dim(X^j)$. En chaque point $x_n$, $P^\perp$ et $\T_{x_n}X^j$ ont une intersection exc\'edentaire en
dimension $\max(0,q-i+1)$. Supposons $i\in [\dim(X^0), \dim(X)]$, qui est le seul cas int\'eressant.
Supposons  que dans un voisinage de $0$ la strate $X^0$ soit 
un sous-espace vectoriel de $\d R^n$ de dimension
$d$. Si $P$ est g\'en\'erique, $X^0\cap P^\perp=\{0\}$ et
si la stratification $(X^j)_{0\in \{0, \cdots, k\}}$ est
$(a^*)$-r\'eguli\`ere, par d\'efinition m\^eme
il existe un sous-analytique dense $ \r W_{n-i+d} $ de l'ensemble des
$(d+n-i)$-plans contenant  $X^0$ tel que si 
$X^0\oplus P^\perp=W\in \r W_{n-i+d} $,
$W$ v\'erifie :

$$\pi_P^{-1}(\pi_P(X^0))=W    \hbox{ coupe
transversalement}  \ X^j \eqno(1)$$
et 
$$(X^0, X^j\cap W)   \hbox{ soit un couple de strates } (a)
\hbox{-r\'egulier en }  0\eqno(2)$$
On a, d'apr\`es $(1)$ : $\dim(W\cap X^j)=q-i+d$ et d'apr\`es ce qui pr\'ec\`ede,
 $P^\perp$ et $\T_{x_n}(X^j\cap W)$
 ont une intersection en dimension au moins $q-i+1$. Comme $P^\perp$ est transverse \`a $X^0$,
$\D \lim_{n\to \infty}\T_{x_n}(X^j\cap W)$ ne peut contenir
l'espace $X^0$ qui est de dimension $d$, ce qui con\-tre\-dit $(2)$. 

 Comme la $(b^*)$-r\'egularit\'e, la
 $(a^*)$-r\'egularit\'e est une  
condition strictement plus faible
que la $(w)$-r\'egularit\'e et on ne sait pas, par exemple, si
en sous-analytique $(b)\Longrightarrow (a^*)$,
 comme c'est le cas en complexe.

\ev
Rappelons maintenant que les caract\'eristiques $\chi_j^P(0)$ 
sont d\'efinies de la fa\c con suivante
(cf Th\'eor\`eme 2.8.(iii) pour les notations) :
si $X$ est un ensemble sous-analytique ferm\'e de $\d R^n$ contenant $0$ et si
$P$ est un $i$-plan vectoriel (g\'en\'eral) de $\d R^n$, 
il existe un r\'eel non nul
$R_0=R_0(P)$, tel que pour tout $r\le R_0$, les ensembles 
sous-analytiques ouverts $ K_1^{P,r},
\cdots,  K_{N_P}^{P,r}$, qui sont les composantes connexes 
du compl\'ementaire dans $P$
de $\r D_{L_r}(P)\cup \D \bigcup_{j=0}^k {\r D}_{X^j\cap B^n_{(0,r)}}(P)$, ont une
 r\'eunion dense dans $P$ au voisinage de $0$,  et v\'erifient~:

- les germes $( K_j^{P,r})_0$ et $( K_j^{P,s})_0$ co\"\i ncident
pour $r,s\le R_0$,

- pour tout $y$ suffisamment proche de $\pi_P(0)=0$ dans un domaine  $ K_j^{P,r}$,
(ie $0<\Vert y\Vert  <<r$),
l'entier $\chi_j^P=\chi(\pi_P^{-1}(y)\cap X\cap B^n_{(0,r)})$
ne d\'epend pas de $r$ (lorsque $r\le R_0$).

\ev
{\bf Notation 4.2.} On notera ${\cal R}_0={\r R}_0(P)$  
le plus grand rayon $R_0$ qui v\'erifie cette~proposition.
\ev
{\bf Remarque.} 
D'apr\`es les Propositions 2.4, 2.5 et le lemme d'isotopie de 
Thom-Mather, 
l'existence du rayon $R_0(P)$ est garantie d\`es que $P\in \r E_X^i$
(notation qui pr\'ec\`ede la Proposition 2.4), d\`es que pour tout
$j\in \{0,\cdots, k\}$,  
$P^\perp$ ne rencontre pas  $\r P_{X^j}(P)$ 
dans $B(0, R_0(P))$ et lorsque pour 
tout $r\le R_0(P)$, $S(0,r)$ est transverse \`a toutes les strates
$X^j$ et  $0\not\in  \r D_{L_r}(P)$.
De sorte que si $R>{\cal R}_0$,   pour tout
$\eta>0$, il existe $y_\eta\in P$, $j\in \{0, \cdots, k\}$
et $x_\eta\in X$ tels que :
\ev 
-  $0<\Vert y_\eta\Vert \le \eta$,
${\cal R}_0\le  \lim_{\eta\to 0} \Vert x_\eta \Vert\le R$,
$x_\eta\in X^j\cap \pi_P^{-1}(y_\eta)$,

-  $\pi_P^{-1}(y_\eta)\cap T_{x_\eta}X^j$  n'est
pas transverse en  $x_\eta $ \à $S_{(0,\Vert x_\eta\Vert)}$.

\ev 

Nous montrons maintenant dans la Proposition 
$4.3$ que la Proposition $2.6$, qui est locale en un 
point, admet une version le long d'une strate d'une 
stratification suffisamment 
r\'eguli\`ere. En particulier~: 
si $Y$ est une strate d'une stratification $(b^*)$-r\'eguli\`ere
de $X$, si $P$ est g\'en\'erique dans $ G(i,n)$, 
il existe un voisinage
$\r U$ de $0$ dans $Y$, deux r\'eels $0<r<<R$, tels
pour tout $y\in Y\cap \r U$,   pour tout 
$z$ et $z'$ dans le m\^eme profil polaire $\K_j^{P,R}(y)$
avec $\Vert z-y \Vert <r$ et $\Vert z'-y \Vert <r$, 
$\pi_P^{-1}(z)\cap X\cap B(y,R)$ et   $\pi_P^{-1}(z')\cap 
X\cap B(y,R)$
sont hom\'eomorphes. 
\ev 
{\bf Proposition 4.3. --- }
{\sl Soit $X$ un ensemble sous-analytique
ferm\'e de $\d R^n$  muni d'une stratification $(b)$-r\'eguli\`ere
  $(X^j)_{j\in \{0, \cdots, k\}}$. Soient
   $X^0$ la strate contenant $0$, 
$i \in \{ 0, \ldots, n\} $ et 
 $\pi_P:\d R^n\to P\in G(i,n)$  une projection telle qu'existe
un voisinage ouvert $\r U$ de l'origine
pour lequel, quelle que soit la strate $ X^j$ adh\'erente \`a $X^0$~:
 
$$\r U \cap  \r P_{X^j}(P) \cap (\pi_P^{-1}(\pi_P(X^0))\setminus X^0)
=\emptyset \eqno(*)$$

\n
Si de plus  les intersections 
$\pi_P^{-1}(X^0)\cap X^j$ sont $(b)$-r\'eguli\`eres le long
de $X^0$, on a : 
\ev 
 \vskip3mm
\n
{(i)} \hskip2mm
 \vtop{\hsize12,95cm \n Le rayon $r'_P(y)$ de la Proposition $2.6$
est minor\'e par un r\'eel $R>0$, ind\'ependant de $y\in X^0$
(dans un voisinage de $0$).}

 \vskip3mm
\n
{(ii)} \hskip2mm
 \vtop{\hsize12,95cm \n  Quel que soit $j\in \{0, \cdots, k\}$, 
la distance $\delta(y)$ de $\r D_{L^R_{X^j}(y)}$ \`a $y$ 
ne tend pas vers $0$ lorsque $y$ tend vers $0$ dans $X^0$.
  } }
\ev
{\bf Preuve.}
(i)- 
Notons que pour chaque $y\in X^0$, $P\in \r E^i_X(y)$, par la condition 
$(*)$. Soit $y\in X^0$.
 La preuve de la Proposition $2.6$ (cf la Remarque qui la suit)
montre que $r'_P(y)$ est r\'ealis\'e d\`es que :

- pour $r\le r'_P(y)$,
les sph\`eres $S(y,r)$ sont transverses aux strates $X^j$,

- en tout point $x$ de $B(y,r'_P)$,
$\T_x(X^j\cap (y+P^\perp))$ et la direction $y-x$ ne sont pas orthogonaux. 
 
Pour le premier point :
 on peut choisir un r\'eel $R>0$ ind\'ependant de $y$
 tel que pour tout
$y\in X^0$, $S(y,r)$ et  $X^j$ sont transverses pour $r\le R$. 
Ceci r\'esulte de la condition $(b)$. Le second point r\'esulte
de la condition $(b)$ pour  le couple 
$(X^0,\pi_P^{-1}(X^0)\cap X^j )$.
\ev
(ii)- On choisit $R$ comme dans $(i)$ et on regarde les links
$L^R_{X^j}(y)$ des strates $X^j$
en les points $y$ de $X^0$. D'apr\`es $(i)$, 
la distance de leur discriminant \`a $y$, not\'ee $\delta(y)$
n'est pas nulle. On veut montrer que $\delta(y)$ est minor\'ee par 
$\delta>0$ ind\'ependamment de $y$. Autrement dit qu'il y existe
un voisinage de $0$ dans $P$   dans lequel n'entrent pas les 
discrimants de la famille de links 
$(L^R_{X^j}(y))_{y\in X^0, j\in \{0, \cdots, k\}}$. Si  tel 
n'est pas le cas,  
on trouve une suite de points $(x_\ell^R)_{\ell\in \N}$ avec 
$x_\ell^R\in S(y_\ell, R)$, $y_\ell\in X_0$ telle que 
$\T_{x_\ell^R} [ L^{R}_{X^j}(y_\ell)]= \T_{x_\ell^R} S_{(y_\ell,R)}\cap
\T_{x_\ell^R}  X^j $ et $P^\perp$ sont non transverses
et $\D\lim_{\ell\to \infty}\pi_P(x_\ell^R)=0$. Comme on l'a remarqu\'e
dans la preuve de la Proposition $2.6$, pour $\pi_P(x_\ell^R)$
suffisamment proche de $0$,  
 $\pi_P(x_\ell^R)+P^\perp$ et  $X^j$ sont transverses en $x_\ell^R$
 (par l'hypoth\`ese $(*)$ et la condition $(a)$)
et de ce fait  on obtient l'inclusion~: $\T_{x_\ell^R}X^j\cap P^\perp
\subset \T_{x_\ell^R}S(y_\ell,R)$. En faisant $\ell\to \infty$, on en 
d\'eduit un point $x^R\in S(0,R)\cap X^m$, pour certaine strate
$X^m$ et par la condition $(a)$ au point $x^R$ pour le couple
$(X^m,X^j)$, l'inclusion 
$\T_{x_\ell^R}X^j\cap P^\perp
\subset \T_{x_\ell^R}S(y_\ell,R)$ donne par passage \`a la limite~:
$\T_{x^R}X^m\cap P^\perp
\subset \T_{x^R}S(0,R)$, ce qui contredit  $\delta(0)\not=0$.\f

 \ev
{\bf Proposition 4.4. ---}
 {\sl Soit $X$ un ensemble sous-analytique
ferm\'e de $\d R^n$, contenant l'origine, muni d'une stratification
$(b^*)$-r\'eguli\`ere $(X^j)_{j\in \{0, \cdots, k\}}$ et soit $X^0$ la strate contenant $0$.
 Soit $\pi_P:\d R^n\to P\in G(i,n)$ une projection telle qu'existe un voisinage ouvert $\r U$ de l'origine
pour lequel quel que soit la strate $ X^j$ adh\'erente \`a $X^0$ :
$$\r U  \cap \r P_{X^j}(P)\cap  (\pi_P^{-1}(\pi_P(X^0))\setminus X^0) 
= \emptyset \eqno(*)$$
 \ev
\n
{(i)}\hskip4mm
 \vtop{\hsize12,95cm \n La projection $\pi_P$ permet de d\'efinir 
les caract\'eristiques $\chi_j^P(y)$
(\`a la fa\c con du Th\'eor\`eme 2.7.(iii)), quels que soient 
$y\in X^0\cap \r U$.}
\ev
\n
{(ii)} \hskip1mm
 \vtop{\hsize12,95cm \n   Il existe un ensemble sous-analytique 
$\r L^i_X$ dense dans $G(i,n)$ tel que pour $P\in \r L^i_X$
 existent un voisinage ouvert $\r U'=\r U'(P)$ de $0$ dans $\d R^n$ et
 un r\'eel ${\cal R}={\cal R}(P)>0$,
tels que pour tout 
$j\in \{0, \cdots, k\}$ :
\ev 
- $\D\bigcup_{y\in \cal U'}B(y,\r R)\cap \r P_{X^j}(P)
\cap  (\pi_P^{-1}(\pi_P(X^0))\setminus X^0) 
= \emptyset$,
 
-  pour tout $y\in X^0\cap \r U'$, ${\cal R} \le {\cal R}_y(P)$. 
\ev
\n 
En cons\'equence chaque caract\'eristique $\chi^P_\ell(y)$ en 
$y$ est donn\'ee par $\chi(\pi_P^{-1}(z)\cap X\cap B(y,\r R))$, pour 
$z$ suffisamment proche de $y$ dans $P$.}

\ev 
\n 
Pour tout $P$ dans le sous-analytique $\r L^i_X\subset G(i,n)$,
pour le r\'eel $\r R=\r R(P)$  donn\'e en (ii) et avec
les notations de la D\'efinition $2.1$ et du Th\'eor\`eme $2.7$, 
il existe un voisinage de $0$ dans $X^0$, tel que pour tout $y$ 
dans ce voisinage :
\ev
\n
{(iii)}\hskip3,5mm
 \vtop{\hsize12,95cm \n
L'ensemble $\D\bigcup_{j=0}^k \r D_{X^j\cap B_{(0,\cal R)}}(P)$ est 
un repr\'esentant du discri\-mi\-nant local
$\Delta_{X_y}(P)$ de $X$ en~$y$ et
chaque profil polaire local
$\r K_\ell^P(y)$ de $X$ en $y$ admet pour repr\'esentant 
un domaine choisi parmi les
 $K_1^{P,\cal R}(0), \cdots, K_{N_P}^{P,\cal R}(0) $ qui 
  contiennent $0$ dans leur adh\'erence.}
\ev
\n
{(iv)}\hskip3mm
 \vtop{\hsize12,95cm \n
Chaque entier $\chi_\ell^P(y)$ associ\'e au profil polaire local 
$\r K_\ell^P(y)$ s'obtient comme
$\chi (\pi_P^{-1}(z)\cap X\cap B^n_{(0,\cal R)})$, pour $z$ suffisamment
proche de $y$ dans le domaine $K_{m_\ell}^{P,\cal R}(0)$
qui repr\'esente $\r K_\ell^P(y)$, et en particulier 
$\chi^P_\ell(y)=\chi^P_{m_\ell}(0)$.}}

\ev
{\bf Preuve.} 
(i)-
 La stratification $(X^j)_{j\in \{0, \cdots, k\}}$ 
\'etant une stratification
  $(b^*)$-r\'eguli\`ere, il s'agit aussi d'une
stratification $(b)$-r\'eguli\ère.
D'apr\`es les hypoth\`eses, la projection $\pi_P$ est dans l'ouvert 
$\r E_X^i=\r E_X^i(y)$ d\éfini dans la section $2$ et attach\'e 
\`a chaque point
$y$ de $X^0\cap \r U$. 
On peut alors appliquer le Th\'eor\`eme $2.7$; 
$\pi_P$ est bien une projection qui permet de
d\'efinir les profils polaires locaux $\r K_j^P$ et 
les caract\'eristiques $\chi_j^P$ qui leur correspondent.

(ii)-    
 La stratification $(X^j)_{j\in \{0, \cdots, k\}}$ 
\'etant par hypoth\`ese une stratification
  $(b^*)$-r\'eguli\`ere, en particulier 
l'ensemble des $P\in G(i,n)$ v\'erifiant la condition $(*)$ est un 
sous-analytique dense dans $G(i,n)$ (cf la remarque qui suit la Proposition 
4.1). Si de plus $P$ est choisi de sorte
que les couples $(X^0,\pi_P^{-1}(X^0)\cap X^j)$ soient $(b)$-r\'eguliers, 
d'apr\`es la Proposition 4.2.(i) l'existence de $\r U'$ et de $\r R$
est garantie. Or un tel choix donne encore un ensemble sous-analytique
dense de $G(i,n)$, puisque  la stratification 
$(X^j)_{j\in \{0, \cdots, k\}}$ est $(b^*)$-r\'eguli\`ere.

(iii)-
Pour tout $y$ dans un voisinage de $0$ dans $X^0$,
$(B(0,\r R)\cap \pi_P^{-1}(y)\cap \r P_{X^j})\setminus\{y\}=\emptyset$, 
puisque $P \in \r L_X^i $.
La remarque qui pr\'ec\`ede la preuve de la Proposition $2.4$ 
montre que cette condition suffit pour que
$\Delta_{X^j_y}(P)=\big(\r D_{X^j\cap B_{(y,\cal R)}}(P)\big)_y  =
 \big(\r D_{X^j\cap B_{(0,\cal R)}}(P)\big)_y  $.
On remarque ensuite que 
 $\Delta_{X_y}(P)$ est le bord de la r\'eunion des profils
 polaires locaux $\r K^P_\ell(y)$.

(iv)- 
Soit $\delta>0$ un r\'eel minorant la distance
$\delta(y)$ de $\r D_{L^{\cal R}_{X^j}}(y)$ \`a $y$
(Proposition 4.2.(ii)), pour 
$j\in \{0, \cdots, k\}$ et pour tout $y$ voisin de $0$ dans $X^0$.
Soit $y$ proche de $0$ dans $X^0$ (et $\Vert y\Vert << \min(\delta, \r R)$)
 de sorte que les profils polaires 
locaux de $y$ admettent pour repr\'esentant certains des domaines 
$K^{P,\cal R}_1(0), \cdots, K^{P,\cal R}_{N_P}(0)$ (d'apr\`es le point 
(iii))
 et soit  $z\in P$ tel que $\Vert y -z\Vert << \delta$. Alors
 d'une part,  pour un certain~$\ell$,
 $\chi(\pi_P^{-1}(z)\cap B(y,\r R)\cap X)=\chi_\ell(y)$ et
 le profil polaire local $\r K_\ell^P(y)$
admet le domaine $K^{P,\cal R}_{m_\ell}(0)$ pour repr\'esentant, 
et d'autre part 
 $\chi(\pi_P^{-1}(z)\cap B(y,\r R)\cap X)=
\chi(\pi_P^{-1}(z)\cap B(0,\r R)\cap X)$, puisque 
$\Vert y\Vert << \delta$ et $\Vert y -z\Vert << \delta$. Enfin 
$\chi(\pi_P^{-1}(z)\cap B(0,\r R)\cap X)=\chi_{m_\ell}(0)$. \f

\ev
On en vient maintenant \`a l'implication $(i)\Rightarrow (iv)$ du Th\'eor\`eme $2.12$ que l'on s'est
propos\'ee de d\'emontrer.
\ev
{\bf Corollaire 4.5. ---}
{\sl Soit $X$ un ensemble analytique complexe de $\d C^n$ muni d'une
stra\-ti\-fi\-cation analytique complexe
$(b)$-r\'eguli\`ere $(X^j)_{j\in \{0, \cdots, k\}}$.
  La suite des caract\'eristiques \'evanescentes 
$(\ti \sigma_i(X_y))_{i\in \{1, \cdots, n\}}$
est constante lorsque $y$ varie le long des strates de  
$(X^j)_{j\in \{0, \cdots, k\}}$.}

\ev
{\bf Preuve.} Voyons $X$ comme un sous-analytique de $\d R^{2n}$ 
et montrons la constance des caract\'eristiques
\'evanescentes le long de la strate $X^0$ au voisinage de $0\in X^0$.
Notons que la stratification $(X^j)_{j\in \{0, \cdots, k\}}$ du 
sous-analytique $X$ est $(b^*)$-r\'eguli\`ere, puisque
$(b)$-r\'eguli\`ere et analytique complexe.
 Soit $i\in \{1, \cdots, n\}$.
On remarque que l'ensemble $\r L_X^i$ de la Proposition  4.4.(ii)
peut \^etre choisi dense dans $\ti G(i,n)$  aussi bien
que dans $G(i,n)$, avec les m\^emes arguments.
Dans ces conditions, soit $P\in \r L_X^i$ et $\r K^{P}(0)$ 
l'unique profil polaire local de
$X$ en $0$ (puisque le compl\'ementaire de 
$ \Delta_{X_0}(P)$ est connexe). La Proposition  4.4.(iv) 
montre que  si $y$ est dans 
un voisinage de $0$ dans $X^0$,
$\r K^P(y)$ et $\r K^{P}(0)$ admettent un m\^eme repr\'esentant 
$K^{P,r}(0)$ et que :
$\ti \sigma_i(X_0)=\ti \sigma_i(X_y).\f$
\ev
{\bf Remarque. } L'implication utile dans la preuve du Corollaire
4.5 est $(b)\Longrightarrow (b^*)$. Implication dont on r\'ep\`ete
qu'on ne sait pas si elle vraie en sous-analytique r\'eel.

\ev
Pour $X$ un sous-ensemble de $\d R^n$ et $Y=\d
R^d\times\{0\}^{n-d}$, on note $\r D(X,Y)$ la d\'eformation de
$X$ sur son c\^one normal \`a $Y$. Il s'agit du
 sous-ensemble suivant de
$\d R^n\times \d R$ :
$$
\r D(X,Y)=\adh\Big\{(y,\ell,u)\in \d R^d\times \d R^{n-d}\times \d
R^*_+ \;;\; (y,u\ell)\in X\setminus Y \Big\}.
$$
L'application $u:\r D(X,Y)\longrightarrow \d R$ induite par la
deuxi\ème projection est appel\ée {\it d\éformation de $X$ au
c\^one normal \`a $Y$ dans $X$}. Cette application est compatible
avec l'action de $\d R^*_+$ d\éfinie par $(\lambda ,(y,\ell,u))
\longmapsto (y,\lambda \ell, \lambda^{-1} u)$ pour $\lambda$ dans
$\d R_+^*$. Sa fibre en $0\in \d R$ est le {\it c\^one normal \`a
$Y$ dans $X$}, not\é $C_YX$ et stable par l'action induite de $\d
R^*_+$ d\écrite par $(\lambda , (y,\ell)) \longmapsto (y,\lambda
\ell)$ pour $\lambda$ dans $\d R_+^*$. Si $Y=\{0\}$, $C_YX$ est le
{\it c\^one tangent \`a $X$ en $0$}. La premi\ère projection de
$\d R^n\times \d R$ induit une application $p: C_YX\longrightarrow
Y$. Nous notons $(C_YX)_y$ la fibre $p^{-1}(y)$. Une telle fibre
est non vide si et seulement si $y\in \adh(X\setminus Y)$.
L'application $(y,\ell,u)\longmapsto (y,u\ell,u)$ induit un
isomorphisme de $\r D(X,Y)\setminus u^{-1}(0)$ sur $\adh(X\setminus
Y)\times \d R^*_+$.

\ev
{\bf D\éfinition 4.6. } On dit que $X$ est {\sl normalement
pseudo-plat le long de $Y$} si la projection $p:C_YX\to Y $ est
ouverte (cf [Hi1]). \ev La dimension des fibres de $p$, la
pseudo-platitude normale, la condition $(b)$ de Whitney et
l'\'equimultiplicit\'e de $X$ le long de $Y$ lorsque $X$ et $Y$
sont analytiques complexes sont li\'es par la proposition suivante
: \ev
{\bf Proposition 4.7. --- } {\sl Soit $X$ et $Y$ deux ensembles sous-analytiques de
$\d R^n$, $X$ \'etant  de dimension $d$ et $Y$ lisse de dimension
$k$.
\ev
\n
{(i)}\hskip3,5mm
 \vtop{\hsize12,95cm \n
Si $Y$ est une strate d'une stratification de Whitney de
$\adh(X)$, $\adh(X)$ est normalement pseudo-plat le long de
$Y$, et si $y\in Y$, $C_{\{y\}}(X\cap (y+Y^\perp))=(C_YX)_y$.
En parti\-cu\-lier le long de strates de Whitney :
$\dim((C_YX)_y)\le d-k$ (cf \hbox{\rm [Hi1], [He-Me2], [Or-Tr]}).}
\ev
\n
{(ii)}\hskip2,6mm
 \vtop{\hsize12,95cm \n Si $X$ est normalement pseudo-plat le long de $Y$, quel que soit  $y\in Y$,
$\dim((C_YX)_y)\le d-k$ (cf (\hbox{\rm [Hi1], [He-Me2], [Or-Tr],
[Co2], Lemme $2.4$})).
 }
\ev
\n
{(iii)}\hskip2mm
 \vtop{\hsize12,95cm \n
Soit $y\in Y$. Notons   $s=\dim((C_YX)_y)$. Il existe un ensemble sous-analytique $\r Q_y^q$
dense dans $G(q,n)$, tel que pour tout $P\in \r Q_y^q$ on ait l'existence d'un voisinage $\r U_{y,P}$ de $y$ dans $\d R^n$
pour lequel :
$$ \big(\adh(X\cap \r U_{y,P}) \setminus Y \big)\cap
\pi_P^{-1}(\pi_P(Y))=\emptyset $$
si et seulement si $s+k\le q$ (cf [Co2], Lemme $2.3$). Lorsque
$\r Q_y^d$ existe, on dit que $X$ est
{\sl \'equis\écable le long de $Y$ en $y$ },
cette condition \'equivaut \`a $\dim((C_YX)_y)\le d-k$
(cf \hbox{\rm [He-Me2]}).
}
\ev
\n
{(iv)}\hskip2mm
 \vtop{\hsize12,95cm \n  Si $X$ et $Y$ sont des ensembles analytiques complexes de $\d C^n$,
la pseudo-platitude normale de $X$ le long de $Y$, l'\'equimultiplicit\'e de $X$ le long de $Y$, l'\'equis\'ecabilit\'e de
$X$ le long de $Y$ et la condition  $\dim((C_YX)_y)= d-k$
pour tout $y\in Y$, sont \'equivalentes (cf \hbox{\rm [Hi1],
[Sc], [He-Me2]}).  }
} \ev
Dans [Co2], Proposition $2.6$, il est montr\'e :
\ev
{\bf Proposition 4.8. --- } {\sl Soient $\r K$ et $Y$ deux
ensembles sous-analytiques de $\d R^i$, $\r K$ \'etant
 de dimension $i$ et $Y$ lisse. Notons $\inte(\r K)$ l'int\'erieur de $\r K$ et
 $\fr(\r K)=\adh(\r K)\setminus \inte(\r K)$ la fronti\`ere de $\r K$.
Si pour tout $y\in Y$,  $\dim((C_Y\fr(\r K))_y)\le i-1-\dim(Y)$, c'est-\`a-dire si $\fr(\r K)$ est \'equis\écable
le long de $Y$,
la fonction :
$Y\ni y\mapsto \Theta_i(\r K_y)$ est continue.}
\ev
Dans le th\'eor\`eme qui suit on identifie une fois de plus une 
strate $Y$ de $(X^j)_{\{0, \cdots, k\}}$ et son projet\'e sur
un plan vectoriel.
Ce th\'eor\`eme montre que l'\'equis\'ecabilit\'e des discriminants g\'en\'eraux en toute
dimension le long des strates d'une stratification  $(b^*)$-r\'eguli\`ere
implique la continuit\'e des invariants $\sigma_*$ et $\Ll_* $ le long des strates.
 \ev

{\bf Th\'eor\ème 4.9. --- } {\sl Soit $X$ un ensemble sous-analytique
ferm\é de $\d R^n$ de dimension $d$,  muni d'une stratification
$(b^*)$-r\'eguli\`ere $(X^j)_{j\in \{0, \cdots, k\}}$ dont  $Y$ soit
une strate. Soient $i\in \{\dim(Y)+1, \cdots, d\} $ et $\r L^i_X $
le sous-analytique dense de $G(i,n)$ de la Proposition  4.4.(ii)
(relatif \`a $Y$). \ev \n
{(i)}\hskip3,5mm
 \vtop{\hsize12,95cm \n
Si pour un ensemble sous-analytique dense $\r D^i\subset \r L_X^i$
de $G(i,n)$, on a :
$$ P\in \r D^i \Longrightarrow \dim\big((C_{Y}
\Delta_{X_y}(P))_y\big)\le i-1-\dim(Y), \ \forall y\in Y,$$
alors la restriction \`a $Y$ de l'invariant  $\sigma_i(X_y)$ est continue en tout point de $Y$.}
\ev
\n
{(ii)}\hskip2,6mm
 \vtop{\hsize12,95cm \n
Si le sous-analytique dense $\r D^i$ existe pour tout $i\in
\{\dim(Y)+1, \cdots, d\} $, les restrictions \`a $Y$ des
invariants  $ \sigma_1(X_y), \cdots, \sigma_n(X_y)$ et
$\Ll_1(X_y), \cdots, \Ll_n(X_y)$ sont continues en tout point de
$Y$. }}
\ev
{\bf Preuve. } D'apr\ès le Th\éor\ème $3.1$, le point (i) entra\^ \i 
ne le point (ii).
Montrons alors~(i). Soient $y\in Y$ et $(y_n)_{n\in \N}$ une suite de $Y$ de limite $y$. Comme la famille
$(\chi_j^P)_{P\in {\cal E }_X^i, j\in \{1, \cdots, n_P\}}$
est uniform\'ement born\'ee,
il suffit de prouver que pour $P$ g\'en\'eral dans $\r L_X^i$ a 
lieu l'\'egalit\'e :
$$ \lim_{n\to \infty} \sum_{j=1}^{n_P(y_n)} \chi_j^P(y_n)\cdot \Theta_i(\r K_j^P(y_n))= \sum_{j=1}^{n_P(y)}
\chi_j^P(y)\cdot \Theta_i(\r K_j^P(y)).$$
Soit alors $P\in \r L_X^i\cap \r D^i$. 
D'apr\`es la Proposition 4.4.(iii)
et (iv), 
il existe un voisinage de $y$ dans $Y$,
tel que pour tout $y_n$ dans ce voisinage, pour tout 
$j\in \{1, \cdots, n_P(y_n)\} $, il existe $\ell_j\in \{1, \cdots,
n_P(y)\}$ tel que $\chi_j^P(y_n)=\chi_{\ell_j}^P(y)$ et 
$\r K_j^P(y_n)= K^{P,r}_{\ell_j}(y)$. Avec la convention
$\Theta_i(K^{P,r}_j(y), y_n)=0$ si $y_n\not\in 
\adh(K^{P,r}_j(y))$,
on obtient pour $n$ suffisamment grand :
$$   \sum_{j=1}^{n_P(y_n)} \chi_j^P(y_n)\cdot
 \Theta_i(\r K_j^P(y_n))= \sum_{j=1}^{n_P(y)}
\chi_j^P(y)\cdot \Theta_i(K^{P,r}_j(y), y_n).$$
Or  puisque
$\fr(\K_j^{P,r}(y)) \subset \D\bigcup_{j=1}^k \r D_{X^j\cap B(y,r)}(P)$
et que $\D\bigcup_{j=1}^k \r D_{X^j\cap B(y,r)}(P)$
est un repr\ésentant du germe $\Delta_{X_{y_n}}(P)$ par la 
Proposition 4.5.(iii), l'hypoth\`ese
assure que la majoration suivante a lieu :
$\dim\big( C_Y(\fr(\K_j^{P,r}(y)))_{y_n}\big) 
\le i-1-\dim(Y)$, 
ce qui par la Proposition $4.8$ donne bien :
$\D\lim_{n\to \infty} \Theta_i(\K_j^{P,r}(y), y_n)=
\Theta_i(\K_j^{P,r}(y), y)=\Theta_i(\r K_j^{P,r}(y))   $.\f

\ev
Nous allons maintenant montrer que l'hypoth\`ese du Th\'eor\`eme
4.9.(ii) a lieu, c'est-\`a-dire que les sous-analytiques 
denses $\r D^{\dim(Y)+1}, \cdots, \r D^d$ existent bien, lorsque la
stratification 
$(X^j)_{j\in \{0, \cdots, k\}}$ est  $(w)$-r\'eguli\`ere.

On consid\ère pour cela la situation suivante : $X$ est un ensemble
sous-analytique ferm\é de $\d R^n$ de dimension $d$,  muni d'une
stratification $(w)$-r\'eguli\`ere $(X^j)_{j\in \{0, \ldots, k\}}$
dont  $Y$ est une strate. On s'int\éresse aux propri\ét\és locales
de la stratification au voisinage d'un point de $Y$. Au voisinage
de ce point, l'adh\'erence de $X\setminus Y$ est $X$. Quitte \à faire
un changement de coordonn\ées locales, on peut supposer que ce
point est l'origine et que $Y$ est un sous-espace vectoriel de $\d
R^n$. Pour $i$ de $\dim(Y) +1$ \à $d$, et un plan $P$ de dimension
$i$ contenant $Y$, on consid\ère la projection orthogonale de
noyau $P^\perp$ sur $P$. On note N le sous-espace vectoriel de
dimension $n-\dim(Y)$ normal \à $Y$ et on identifie $\d R^n$ avec
le produit $Y\times \hbox{\tm N}$.

\`A un morphisme sous-analytique lisse $f : Z \longrightarrow \d
R$ d\éfini sur un ensemble sous-analytique de dimension $d$ dans
$\d R^n$, on associe le fibr\é conormal relatif $T^*_f \d R^n$,
ensemble sous-analytique lisse de dimension $n$ du fibr\é
cotangent $T^*\d R^n$. Pour un morphisme sous-analytique $f : Z
\longrightarrow \d R$ lisse sur un ouvert partout dense $U$
de $Z$, le {\it conormal relatif} $T^*_f \d R^n$ est, par
d\éfinition, l'adh\érence de $T^*_{f|_U}$ dans $T^*\d
R^n$. Cet ensemble est stable par les homoth\éties r\éelles du
fibr\é vectoriel cotangent $T^*\d R^n$. Lorsque $f$ est constante
sur $Z$, on note $T^*_Z \d R^n$ au lieu de $T^*_f \d R^n$.
Lorsqu'il n'y a pas d'ambigu\"{\i}t\é, on abr\ège $T^*_f \d R^n$
en $T^*_f$.

\ev
{\bf Th\'eor\ème 4.10. --- } {\sl  Soit $X$ un ensemble
sous-analytique ferm\é de $\d R^n$ de dimension $d$,  muni d'une
stratification $(w)$-r\'eguli\`ere $(X^j)_{j\in \{0, \ldots, k\}}$
dont  $Y$ est une strate. Pour $i$ de $\dim(Y)+1$ \à $d$, les deux
propri\ét\és suivantes sont vraies :

\ev \n
{(i)}\hskip3,5mm
 \vtop{\hsize12,95cm \n \sl
Il existe un ouvert sous-analytique partout dense $\r D^i$ dans
$G(i-\dim(Y), \hbox{\tm N})$ tel que, pour $P$ dans $\r D^i$, le plan
$P^\perp$ ne rencontre le c{\^o}ne normal $C_Y(\r P_{X}(P))$ qu'\à
l'origine.}

\ev \n
{(ii)}\hskip2,6mm
 \vtop{\hsize12,95cm \n \sl
Il existe un ouvert sous-analytique partout dense $\r D^i$ dans
$G(i-\dim(Y), \hbox{\tm N})$ et un ouvert sous-analytique partout
dense $\r D^{n-i+1}$ dans $G(n-i+1,\hbox{\tm N})$, tels que, pour $P$
(contenant $Y$) dans $\r D^i$ et $W$ (orthogonal \à $Y$) dans $\r
D^{n-i+1}$, le plan $W$ ne rencontre le c{\^o}ne normal $C_Y(\r
P_{X}(P))$ qu'\à l'origine.}

\ev \n
{(iii)}\hskip1,7mm
 \vtop{\hsize12,95cm \n \sl
Il existe un ouvert sous-analytique partout dense $\r A^i$ dans la
vari\ét\é des drapeaux $F(1, i-\dim (Y), \hbox{\tm N})$ tel que pour $(L\subset
P)$ \él\ément de $\r A^i$, $L$ normal \à $Y$ dans $P$, le plan
$\pi_P^{-1}(L)$ ne rencontre le c{\^o}ne normal $C_Y(\r P_{X}(P))$
qu'\à l'origine.}
}

\ev
{\bf Preuve. }On consid{\`e}re la d{\'e}formation $u: \r
D(X,Y)\longrightarrow \d R$ de $X$ au c{\^o}ne normal \à $Y$,
naturellement plong\ée dans le produit $T_Y\d R^n \times \d R$. Le
conormal relatif de cette d{\'e}formation s'identifie {\`a} la
d{\'e}formation $\r D(T^*_X\d R^n,T^*_Y\d R^n)$ de $T^*_X\d R^n$
au c{\^o}ne normal {\`a} $T^*_Y\d R^n$, elle-m{\^e}me sous espace
du produit $T^*(T_Y\d R^n) \times \d R$. Nous noterons $\hbox{\ce
X}$ (resp. $\hbox{\ce T}$) l'espace total de la d\éformation de
$X$ au c{\^o}ne normal \à $Y$, (resp. de la d\éformation de
$T^*_X\d R^n$ au c\^one normal \à $T^*_Y\d R^n$) et sa fibre
sp\'eciale $\hbox{\ce X}_0$ (resp. $\hbox{\ce T}_0$). La fibre
${\hbox{\ce T}}_{\!0}$ est un sous-ensemble de $T^*(T_Y\d R^n)$.
On d{\'e}signe par $\pi$ la projection naturelle du fibr{\'e}
normal $T_Y\d R^n$ sur $Y$. Elle induit une projection ${\rm pr}$
de $T^*(T_Y\d R^n)$ sur le fibr\é cotangent relatif $T^*{\pi}$.
D'autre part, l'isomorphisme hamiltonien identifie le c{\^o}ne
normal {\`a} $T^*_Y\d R^n$ dans $T^*\d R^n$ au fibr\é cotangent
$T^*(T^*_Y\d R^n)$. Si on d{\'e}signe par $\pi'$ la projection
naturelle du fibr{\'e} conormal $T^*_Y\d R^n$ sur $Y$, elle induit
\à son tour une projection ${\rm pr}'$ de $T^*(T^*_Y\d R^n)$ sur
$T^*{\pi'}$. On remarque que ${\rm pr}'$ s'identifie \à ${\rm pr}$
et $T^*{\pi'}$ \à $T^*{\pi}$ via l'isomorphisme hamiltonien.

L'{\'e}nonc{\'e} qui suit a {\'e}t{\'e} prouv{\'e}, lorsque $Y$
est un point par Kashiwara dans [Ka2] dans le cadre analytique complexe.

\ev
{\bf Lemme 4.11. --- } {\sl L'ensemble ${\rm pr}({\hbox{\ce
T}}_{\!0})$ est un sous-ensemble {\sl relativement isotrope} et
homog\ène du cotangent relatif $T^*\pi$.}

\ev
{\bf Preuve du lemme 4.11. }Le r{\'e}sultat se d{\'e}montre
localement en un point lisse de la fibre ${\hbox{\ce T}}_{\!0}$ de
$\hbox{\ce T}$ dont la projection dans $T_Y\d R^n$ (resp. $T^*_Y\d
R^n$) n'est pas dans la section nulle de $T_Y\d R^n$ (resp.
$T^*_Y\d R^n$). Gr{\^a}ce au lemme d'aile ([Dr-Mi], [Lo]) on se
ram{\`e}ne au cas o{\`u} ${\hbox{\ce T}}_{\!0}$ est de codimension
1 dans $\hbox{\ce T}$. Dans ce cas, l'aile est
une param{\'e}trisation locale $\varphi$ de $\hbox{\ce T}$,
$\varphi : \hbox{\ce T}_0 \times (0,\alpha) \longrightarrow
\hbox{\ce T}$ telle que que $\varphi(s,t)=(\varphi(s,0),
\rho(s,t))$.

On choisit des coordonn{\'e}es locales $(y,x)$ dans  $\d R^n$ de
telle sorte que $Y$ est d{\'e}fini par l'id{\'e}al $(x)$. On note
$(\eta,\xi)$ les coordonn{\'e}es correspondantes sur les fibres du
fibr{\'e} cotangent $T^*\d R^n$. Sur $T^*(T_Y\d R^n)$ on en
d{\'e}duit des coordonn{\'e}es $(y,\ell,\nu,\xi)$ qui
s'{\'e}tendent en des coordonn{\'e}es $(y,\ell,\nu,\xi,u)$ sur le
produit $T^*(T_Y\d R^n)\times \d R$. L'application :
 $$
(y,\ell,u)\longmapsto(y,u\ell,u)
 $$
induit un isomorphisme de $\hbox{\ce X}\setminus \hbox{\ce X}_0$
sur $X\times \d R_+^*$. Elle se rel{\`e}ve en une application :
 $$
(y,\ell,\nu,\xi,u)\longmapsto(y,u\ell,u\nu,\xi,u)
 $$
qui induit un isomorphisme de $\hbox{\ce T}\setminus {\hbox{\ce
T}}_{\!0}$ sur $T^*_X \d R^n\times \d R_+^*$. Par d{\'e}finition
de $T^*_X\d R^n$, la forme diff{\'e}rentielle $\eta dy+ \xi dx$
est nulle en restriction {\`a} la partie lisse de $T^*_X\d R^n$.
On en d{\'e}duit l'{\'e}galit{\'e} suivante, vraie en restriction
aux points lisses de $\hbox{\ce T}\setminus {\hbox{\ce T}}_{\!0}$
:
 $$
{{\xi}\over{\|\xi\|\,\|\ell\|}} d \ell =
-{{\nu}\over{\|\xi\|\,\|\ell\|}} dy - {{\langle\xi\cdot
\ell\rangle}\over{\|\xi\|\,\|\ell\|}} {{du}\over{u}}.
 $$
La condition $(w)$ se traduit par le fait que la quantit{\'e} :
 $$
{{\|\nu\|}\over{\| \xi\|\,\|\ell\|}}
 $$
est localement born{\'e}e. La condition $(b')$, cons{\'e}quence de
$(w)$, se traduit par le fait que la quantit{\'e} :
 $$
{{\langle\xi\cdot \ell\rangle}\over{\|\xi\|\,\|\ell\|}}
 $$
tend vers 0 avec $t$. La forme diff{\'e}rentielle $\D
{{\langle\xi\cdot \ell\rangle}\over{\|\xi\|\,\|\ell\|}}
{{du}\over{u}}$ s'{\'e}tend donc en une forme sur $\hbox{\ce T}$
au voisinage des points consid{\'e}r{\'e}s. On en conclut que la
forme diff{\'e}rentielle $\D{{\xi}\over{\|\xi\|\,\|\ell\|}} d
\ell$ induit une forme diff{\'e}rentielle relative au dessus de
$Y$, nulle sur la partie lisse de ${\rm pr}({\hbox{\ce
T}}_{\!0})$.\f

\ev {\bf Fin de la preuve du Th\éor\ème 4.10. } Consid\érons le
fibr\é cotangent relatif $T^*\pi$ et sa projection canonique $p$
sur $T_Y\d R^n$. L'isomorphisme hamiltonien identifie $T^*\pi$ au
fibr\é cotangent relatif $T^*{\pi'}$ de la projection $\pi':
T^*_Y\d R^n \longrightarrow Y$. On a donc \également une
projection canonique $p'$ de $T^*{\pi}$ sur $T^*_Y\d R^n$.
Un plan $P$ de dimension $i$ contenant $Y$ \étant donn\é dans $\d
R^n$, il induit, gr{\^a}ce \à la structure euclidienne de $\d
R^n$, un sous-fibr\é de rang $i-\dim(Y)$ de $T^*_Y\d R^n$, not\é
encore $P$. D'autre part, l'orthogonal $P^\perp$ induit un
sous-fibr\é de rang $n-i$ de $T_Y\d R^n$, not\é encore $P^\perp$.
Consid\érons le produit fibr\é $P^\perp\times_Y P$. Le Lemme 4.11
a les cons\équences suivantes : il existe un ouvert $\r D^i$
partout dense dans la grassmannienne $G(i-\dim(Y),\hbox{\tm N})$ tel
que pour $V$ de dimension $i$ contenant $Y$,
 \ev \n
 - \hskip1,5mm
 \vtop{\hsize13,3cm \n 
L'intersection de $P^\perp\times_Y P$ avec tout sous-ensemble
relativement isotrope de $T^*\pi$, en particulier avec ${\rm
pr}({\hbox{\ce T}}_{\!0})$ est r\éduite \à la section nulle de
$T^*\pi$.}
 \ev \n
-  \hskip1,5mm
 \vtop{\hsize13,3cm \n 
 La projection $p(p'^{-1}(P)\cap {\rm pr}({\r
 T}_{\!0}))$ est le c{\^o}ne normal \à la polaire
 $\r P_X(P)$.}
\ev \n Nous venons de montrer la premi\ère assertion du
th\éor\ème. Pour terminer la preuve, il suffit de v\érifier que
les trois assertions sont en fait \équivalentes. On consid{\`e}re
un espace vectoriel $E$ de dimension $k$ et la vari{\'e}t{\'e} des
drapeaux $\ell\subset H$ form{\'e}s d'une droite $\ell$ et d'un
hyperplan $H$ de $E$. On la note $\r F$. Elle est contenue dans le
produit $\d P\times \d P^\vee$ du projectif des droites par le
projectif des hyperplans de $E$. Si on se donne un drapeau :
 $$
{\bf D} : \{0\}=V_0\subset V_{1}\subset \ldots \subset
V_{k-1}\subset V_k=E,
 $$
les cycles $\sigma_i({\bf D}) : \ell\subset V_i \subset H$ pour
$i$ de 1 {\`a} $k-1$ engendrent la cohomologie de $\r F$ en
codimension $k-1$. Il en est de m\ême des cycles :
 $$
\varsigma_i({\bf D}) : \ell \subset V_{i+1}, \; V_{i} \subset H.
 $$
pour $i$ de $1$ \`a $k-1$, puisque dans la cohomologie de $\r F$ on
a : $\varsigma_i({\bf D})=\sigma_i({\bf D})+\sigma_{i+1}({\bf D})$
pour $i$ de $1$ \`a $k-1$.

Donnons-nous un autre drapeau ${\bf D}': W_{1}\subset \ldots
\subset W_k=E$ de sous-espaces vectoriels de $E$. La
vari{\'e}t{\'e} de drapeaux $\r F$ est une sous-vari{\'e}t{\'e} de
codimension 1 dans le produit $\d P\times \d P^\vee$ o{\`u} la
cohomologie en codimension $k-1$ est engendr{\'e}e par les cycles
 $$
\varsigma_i({\bf D},{\bf D}') : \ell \subset W_{i+1}, \; V_{i}
\subset H.
 $$
pour $i$ de 0 \`a $k-1$. Pour un sous-ensemble $\Sigma$ de $\r
F$, les quatre propri{\'e}t{\'e}s suivantes sont donc
{\'e}quivalentes :
 \ev \n
{1 - }\hskip1mm
 \vtop{\hsize12,95cm \n $\Sigma$ est de dimension strictement inf\'erieure {\`a}
    $k-1$.}
 \ev \n
{2 - }\hskip1mm
 \vtop{\hsize12,95cm \n Il existe un ouvert partout dense de drapeaux tel que
    pour tout drapeau ${\bf D}$ de cet ouvert et $i$ de 1 {\`a} $k-1$,
    l'intersection de $\Sigma$ avec les $\sigma_i({\bf D})$ est
    vide.}
 \ev \n
{3 - }\hskip1mm
 \vtop{\hsize12,95cm \n Il existe un ouvert partout dense des couples de drapeaux tel que
pour tout couple de drapeaux $({\bf D},{\bf D}')$ de cet ouvert et
$i$ de 0 {\`a} $k-1$, l'intersection de $\Sigma$ avec les
$\varsigma_i({\bf D},{\bf D}')$ est vide.}
 \ev \n
{4 - }\hskip1mm
 \vtop{\hsize12,95cm \n Il existe un ouvert partout dense de drapeaux tel que
    pour tout drapeau ${\bf D}$ de cet ouvert et $i$ de 1 {\`a} $k-1$,
    l'intersection de $\Sigma$ avec les $\varsigma_i({\bf D})$ est
    vide. \f}

\ev
{\bf Remarque.} On note que la condition $4.10.(i)$ 
est \'equivalente \`a la condition $(b^*)$
le long de $Y$ et \`a  la suivante :
\ev
\n
$4.10.(iv)$  : \hskip3,5mm
 \vtop{\hsize11,6cm \n \sl
Il existe un ouvert sous-analytique partout dense $\r G^i$ dans
$G(i-\dim(Y), \hbox{\tm N})$ tel que, pour $P$ dans $\r G^i$, 
$P$ v\'erifie la condition $(*)$ introduite au d\'ebut de la section 
$4$.}
\ev
\n
D'autre part tout sous-ensemble isotrope de $\d P\times \d P^\vee$
\'etant de dimension strictement inf\'erieure \`a $k-1$, la conclusion 
du Lemme 4.11 implique chacune des conditions \'equivalentes 4.10.(i) 
\`a 
4.10.(iv). Nous ne savons pas si la r\'eciproque est vraie.

\Ev
\centerline{\bf  Appendice.   Poly\`edres et valuations sph\'eriques -
Calcul  des coefficients $m_i^j$.}
\Ev

Nous traitons ici le cas o\`u l'ensemble $X_0$ (not\'e $V$) est
un c\^one poly\édral  de $\d R^n$ de
sommet l'origine, de dimension  maximale et
 nous expliquons comment le Th\'eor\`eme $3.1$ se rattache dans ce
contexte \`a
des questions de g\'eom\'etrie convexe sph\'erique. Le cas conique
poly\'edral permet en outre un calcul ais\'e des coefficients
$m_i^j$ de la matrice $\r M$ du  Th\'eor\`eme~$3.1$.

\ev
Par {\sl poly\èdre de $\d R^n$ } nous entendons une partie de $\d R^n$
qui est une intersection finie de demi-espaces ferm\és. Un  poly\èdre de $\d R^n$
est ainsi une partie convexe de  $\d R^n$. Un {\sl c\^one   poly\édral de
$\d R^n$ de sommet l'origine } est un poly\èdre de $\d R^n$ obtenu comme
intersection de demi-espaces ferm\és contenant tous dans leur bord l'origine.
Par {\sl polytope de  $\d R^n$ } nous entendons un poly\èdre compact de
$\d R^n$, c'est-\à-dire l'enveloppe convexe d'un nombre fini de points de
$\d R^n$.
Enfin un {\sl polytope sph\érique (de $S^{n-1}$) }  est l'intersection
d'un c\^one poly\édral (de $\d R^n$) de sommet l'origine et de la sph\ère
$S^{n-1}$.
Nous noterons $\r K^n$ l'ensemble des convexes compacts de $\d R^n$ et
$\r KS^{n-1}$ l'ensemble des convexes compacts de $S^{n-1}$.
\ev

Nous introduisons maintenant bri\èvement la notion de valuations sur les
convexes. Pour plus de d\étails on pourra se reporter \à [Sc3,4] ou
[McMu-Sc].

Une application $v: \r K^n \to \d R$ (resp.
$v: \r KS^{n-1}\to \d R $) est une {valuation r\éelle} (resp. une
{\sl valuation sph\érique r\éelle}) si $v(\emptyset)=0$ et
si quels que soient $K, L \in \r K^n$
(resp. $K, L \in \r KS^{n-1}$) tels que $K\cup L\in \r K^n$
(resp. $K\cup L \in \r KS^{n-1}$), on a :
$$ v(K \cup L)=v(K) \cup v(L) - v(K\cap L).  $$
On dit qu'une valuation $v$ sur $\r K^n$ (resp. $\r KS^{n-1}$) est {\sl
continue} si elle est continue pour la m\étrique de Hausdorff sur
$\r K^n$ (resp. sur $\r KS^{n-1}$).
On dit qu'une valuation  $v$ sur $\r K^n$ (resp. $\r KS^{n-1}$) est {\sl
simple } si la restriction de $v$ aux convexes de dimension non maximale
est nulle.
Soit $\r Iso^n $ le groupe des isom\étries de $\d R^n$.
Ce groupe agit sur $\r K^n$, de m\ême que le groupe $\r O^n$
des isom\'etries de la sph\`ere agit sur $\r KS^{n-1}$.
 On dit qu'une valuation  $v$ sur $\r K^n$ (resp. $\r KS^{n-1}$) est
{\sl invariante (sous l'action de $\r Iso^n $, resp. sous l'action de
$\r O^n$ )} si quel que soient  $K\in \r K^n$ (resp. $K\in \r KS^{n-1})$,
  $g\in \r Iso^n $  (resp.  $g\in\r O^n$),
on a : $v(g.K)=v(K)$.

\ev
Soient $i,j \in \{ 0, \cdots, n\}$.
Si $C$ est un c\^one convexe sous-analytique de $\d R^n$ de sommet
l'origine, par $i$-homog\én\éit\é de $\L_i$ et par convexit\'e, on a :
$$\L_i^{\ell oc}(C_0)= \D{1\over \alpha_i}\cdot \L_i(C \cap B^n_{(0,1)})=
\D{1\over \alpha_i}\cdot \bar\gamma_{n-i,n}(\{ \bar P \in \bar G(n-i,n);
\ \bar P \cap
C \cap  B^n_{(0,1)})   \not=  \emptyset\}) \ \hbox{ et }$$

$$ \sigma_j(C_0 )
=\int_{P\in  G(j,n) } \ \ \Theta_j\big((\pi_P(C))_0\big)   \ \
d \gamma_{j,n}(P). $$
\ev
\n
Notons $\widehat\L_i^{\ell oc}$ et $\widehat\sigma_j$ les
valuations sph\ériques  sur $\r KS^{n-1}$,
auxquelles donnent lieu $\L_i^{\ell oc}$ et $\sigma_j$~:
$$\forall K \in \r KS^{n-1}, \ \ \  \widehat{\L_i}(K)=\L_i^{\ell oc}
(\widehat K_0)\ \
\hbox{ et } \ \ \    \widehat{\sigma_j}(K)=  \sigma_j(\widehat K_0),$$
o\ù $\widehat K$ est le c\^one dans $\d R^n$ sur $K$ de sommet l'origine.
Ces valuations sont
continues et invariantes sous l'action de $\r O^n$.

\ev
En tant qu'objet d'\étude syst\ématique
l'introduction des valuations vues comme
invariants de d\'ecoupage,
remonte au troisi\ème probl\ème de Hilbert (voir [Bo], [Sah]). Cette
\'etude
culmine dans la caract\érisation, due \à Hadwiger ([Had], [Kl])
des valuations invariantes et continues :
\ev
{\bf Th\éor\`eme }  ([Had], [Kl]) {\bf . --- }
{\sl Si $v$ est une valuation sur $\r K^n$ continue et invariante sous
l'action de $\r Iso^n$, il existe des constantes r\éelles
$\alpha_0, \cdots, \alpha_n$
telles que : $v = \D \sum_{i=0}^n \alpha_i \cdot \L_i$. }
\ev
\n
Une fa\c con \équivalente d'\énoncer  le th\éor\ème  d'Hadwiger
est la suivante :
\ev
{\sl Si $v$ est une valuation simple sur $\r K^n$, continue et invariante
sous l'action de $\r Iso^n$, il existe une constante r\éelle $\alpha$
telle que : $v=\alpha\cdot \L_n=\alpha\cdot \r H^n $.  }
\ev
\n
Un probl\ème apparemment d\'elicat, consiste
\à savoir si le th\éor\ème d'Hadwiger admet un \équivalent sph\érique :
\ev
{\bf Question}  ([Gr-Sc] Probl\`eme 74, [Sc-McMu] Probl\`eme 14.3) {\bf . --- }
{\sl Si $v$ est une valuation simple sur $\r KS^{n-1}$,
continue et invariante sous l'action de $\r O^n$, est-il vrai que $v$
est proportionnelle au volume $\r H^{n-1}$ de $S^{n-1}_{(0,1)}$ ?
}
\ev
Une r\éponse positive \à cette question  aurait pour corollaire
imm\'ediat le Th\éor\ème $3.1$ et inversement la preuve du Th\'eor\`eme
$3.1$ est une r\'eponse positive \`a la question de la validit\'e
de la version sph\'erique du th\'eor\`eme de Hadwiger, dans le cas
particulier des valuations $\widehat \sigma_*$ et $\widehat\L_*$.
Notons que l'on
 dispose d'une r\éponse positive \à cette question, dans le cas
o\ù $n\le 3$ (cf [McMu-Sc] Th\'eor\`eme 14.4) et dans le cas o\`u
la valuation est de signe constant. Dans ce dernier
cas la continuit\é de la valuation n'est pas requise et
la valuation est d\éfinie a priori sur les polytopes convexes
et non pas n\écessairement sur
tous les convexes de $S^{n-1}$~:
\ev
{\bf Th\éor\ème A.1.} ([Sc1] Th\'eor\`eme 6.2, [Sc2]) {\bf   ---}
{\sl Soit $v$ une valuation simple d\éfinie sur les polytopes  de $S^{n-1}$,
invariante sous l'action de $\r O^n$ et \à valeurs dans $\d R_+$. Il existe
alors $c\in \d R_+$ tel que  $v=c\cdot \r H^{n-1}_{|S^{n-1}}$.    }

\ev
Venons-en maintenant au calcul des coefficients de la matrice $\r M$
du Th\'eor\`eme 3.1.
\ev
Si $V$ est  un poly\èdre de $\d R^n$ de dimension $n$, on dit  qu'un
hyperplan qui le borde
est une {\sl facette de $V$}. Le vecteur normal \à une facette $F$ de $V$
est le vecteur unitaire orthogonal \à $F$ situ\é dans le
demi-espace d\éfini
par $F$ qui ne contient pas $V$. Pour $i\in \{0, \cdots , n-1 \}$,
une {\sl $i$-face } de $V$ est l'intersection de
$n-i$ facettes distinctes de $V$ et de $V$.
On note $\r F_{i}(V) $ l'ensemble des $i$-faces de $V$. Par convention
$\r F_n(V)=\{  V \}$.
Si $V$ est un polytope sph\érique, pour $i\in \{0, \cdots, n-1 \}$,
 une $i$-face de $V$ est d\éfinie comme
\étant l'intersection de $S^{n-1}$ et d'une  $(i+1)$-face de $\widehat V$.
\`A tout point $x\in V$ on peut associer $F_x$, l'unique
face de $V$ de dimension minimale contenant $x$.
Si $x\in \partial V$ (le bord de $V$), on d\éfinit $C(x,V)$,
le {\sl c\^one conormal} \à $V$ en $x$,
de la fa\c con suivante :  $C(x,V)$ est le c\^one positif de
$\T_x\d R^n$ engendr\é par les vecteurs normaux aux
facettes de $V$ contenant $x$. On convient que :
$C(x,V)=\{ 0 \}$, lorsque $x\in V \setminus \partial V$.
\ev
{\bf Remarque. }
 Si $V$ est  un poly\èdre de $\d R^n$,  $x\in V$, et
si $F_x$ est de dimension $i\in \{0, \cdots, n-1 \}$,
$C(x,V)$ est un c\^one
de dimension $n-i$. De plus quel que soit $y \in F_x$,
$C(x,V)=C(y,V)$. On d\éfinit donc le {\sl
c\^one conormal  de $V$ le long d'une face $F$ de $V$} par~:
$C(F,V)=C(x,V)$, o\`u $x$ est quelconque dans $F$. On a :
$C(V,V)=\{ 0 \}$.
\ev
Si $V$ est  un poly\èdre
d\ég\én\ér\é de $\d R^n$, c'est-\à-dire si le
sous-espace vectoriel $[V]$ de $\d R^n$ engendr\é par $V$
est de dimension
$<n$, on note $C_{[V]}(x,V)$ le c\^one conormal de $V$ en $x$
dans $[V]$, au sens de la d\'efinition qui pr\'ec\`ede, puisque
$V$ est de codimension nulle dans $[V]$.
Avec cette notation,
le c\^one conormal de $V$ en $x$ dans $\d R^n$,
not\é $C_{\R^n}(x,V)$
est d\éfini par $C_{[V]}(x,V)\times [V]^\perp $.
Le c\^one
conormal est ainsi relatif \à l'espace dans lequel il est obtenu,
mais nous donnons
une d\éfinition intrins\èque attach\ée au c\^one conormal :
l'angle ext\érieur.

\ev
{\bf D\éfinition A.2. }
Soit $V$  un poly\`edre $\d R^n$ et $F\in \r F_i(V)$. On d\éfinit
$\gamma(F,V)$, l'{\sl angle ext\érieur de $V$ en $F$} par :
$$\gamma(F,V)=  {1\over \alpha_{n-i}}(\r H^{n-i}(C(F,V)\cap B^n_{(0,1)}))=
\Theta_{n-i}(C(F,V)_0).  $$
Par convention  $\gamma(V,V)=1$ (ce qui, compte tenu de $C(V,V)=\{ 0 \}$,
signifie que l'on opte pour la convention $\Theta_0(\{ 0  \})=1$).
\ev
{\bf Remarque.}
Comme la densit\é d'un produit  est le produit
des densit\és des facteurs (dans les  dimensions ad\équates),
lorsque $V$ est un poly\èdre  d\ég\én\ér\é de $\d R^n$, l'\égalit\é
 $C_{\R^n}(x,V)= C_{[V]}(x,V)\times [V]^\perp $ montre que l'angle
exterieur
de $V$ en l'une de ses faces ne d\épend pas de l'espace ambiant dans
lequel
on le calcule.
\ev
La formule cin\ématique principale donne l'expression de $\L_i$
\`a l'aide des angles ext\'erieurs attach\'es aux $i$-faces :
\ev
{\bf Th\éor\ème (formule cin\ématique principale)} [Sc3] 4.5.2,
[Sc4] 7.2]) {\bf . --- }
{\sl  Soit $V$ un polytope de $\d R^n$, pour tout $i \in \{0, \cdots, n \}$,
on a :
$$ \L_i(V)= \D \sum_{F\in {\cal F}_i}^n\gamma(F,V)\cdot \r H^i(F). $$}
\ev
Le Th\éor\ème d'Hadwiger implique  alors  que  toute valuation
$v$ sur $\r K^n$,
continue et invariante sous l'action de $\r Iso^n$ v\érifie :
$$\exists \alpha_1, \cdots, \alpha_n \in \d R, \ \
\forall V \ \hbox{ polytope de } \d R^n, \ \ \
 v=\sum_{i=0}^n\alpha_i\sum_{F \in {\cal F}_i} \gamma(F,V)\cdot \r H^i(F)
\eqno{(*)}$$
\n
Bien que nous ne sachions pas si l'\équivalent sph\érique du Th\éor\ème
d'Hadwiger est
vrai, nous allons montrer, pour prouver le Th\éor\ème $3.1$
sur les c\^ones
poly\édraux, que les valuations $\widehat \L_i$ et
$\widehat \sigma_j$ sur les polytopes
de $S^{n-1}$
sont des combinaisons lin\éaires des \équivalents sph\ériques des
$ \gamma(F,V)\cdot \r H^i(F)  $.

Soit $V$ un polytope de
$S^{n-1}$ et $\widehat V$ le c\^one poly\édral de sommet l'origine de
$\d R^n $ qui lui est associ\é. Si $F$ est une $k$-face de $V$, on note
$\widehat F$ la $(k+1)$-face de $\widehat V$ qui lui est associ\ée.
Nous allons
prouver qu'existent  des constantes $a_0, \cdots, a_{n-1}, b_0,
\cdots, b_{n-1}$
ind\épendantes de $V$ telles que :
$$\widehat \L_i(V)=\sum_{k=0}^{n-1} a_k\sum_{ F\in {\cal F}_k( V)}
\gamma(\widehat F, \widehat V)\cdot \r H^{k}( F )  \ \
\hbox{ et }    \ \widehat\sigma_j(V)=
\sum_{k=0}^{n-1} b_k\sum_{ F\in {\cal F}_k( V)}
\gamma(\widehat F, \widehat V)\cdot \r H^{k}( F )\eqno(**) $$
En prouvant $(**)$ nous allons ainsi prouver l'\équivalent
sph\érique de $(*)$ pour
les valuations sph\ériques particuli\ères que sont
$\widehat \L_i$ et $\widehat \sigma_j$,
$ i,j\in \{ 0, \cdots, n \}$, c'est-\`a-dire r\'epondre positivement
au Probl\`eme 15.5 de [McMu-Sc] pour les valuations $\widehat \L_i$.
Remarquons que dans les \égalit\és $(**)$, $\r H^k(F)=
\r H^{k}(S^k)\cdot \Theta_{k+1}(\widehat F_0)$.
\ev
{\bf Lemme A.3. --- }
{\sl Soit $V$ un c\^one poly\édral de $\d R^n$  de sommet l'origine. On a :
$$ \sum_{k=0}^n \ \ \sum_{F \in {\cal F}_k(V) }
 \gamma(F,V)\cdot\Theta_k(F_0) =1 $$    }
\ev
{\bf Preuve.} Soit $F\in \r F_k(V)$, le c\^one $F+C(F,V)$ est de dimension $n$
et sa densit\é en l'origine est $\Theta_n((F+C(F,V)_0) )=\Theta_k(F_0)\cdot
\Theta_{n-k}((C(F,V))_0)=\Theta_k(F_0)\cdot\gamma(F,V) $.
Or comme $V$ est convexe,
$\D \bigcup_{k=0}^n \bigcup_{F\in {\cal F}_k(V)} F+C(F,V)=\d R^n$ et comme
lorsque $F\not= G$, $ \big( F+C(F,V)\big) \cap \big( G+C(G,V)\big) $
est un c\^one de dimension $<n$, on en d\éduit que :
$$\Theta_n\D \big((\bigcup_{k=0}^n \bigcup_{F\in {\cal F}_k(V)}
F+C(F,V))_0\big)= \sum_{k=0}^n \ \ \sum_{F \in {\cal F}_k(V)}
\Theta_n \big( (F+C(F,V))_0 \big)=1.\f $$
\ev
{\bf Th\éor\ème A.4. --- }
{\sl  Soit $j\in \{0, \cdots, n \}$ et
soit $V$ un c\^one poly\édral de $\d R^n$ de sommet l'origine. On a alors :
$$\sigma_j(V_0)= \sum_{k=j}^n \ \ \sum_{F\in {\cal F}_k(V)}
\gamma(F,V)\cdot\Theta_k(F_0).  $$}
\ev
{\bf Preuve.}
Par d\éfinition $\sigma_j(V_0)=\D \int_{P\in G(j,n)}
\Theta_j\big((\pi_P(V))_0\big) \ d \gamma_{j,n}(P)$. Par le Lemme~$A.3$
on obtient :
$$\sigma_j(V_0)= 1- \D \int_{P\in G(j,n)}\
\sum_{k=0}^{j-1} \ \ \sum_{G \in {\cal F}_k(\pi_P(V)) }
 \gamma(G,\pi_P(V))\cdot\Theta_k(G_0)\  d \gamma_{j,n}(P). $$
Soit $P\in G(j,n)$ g\én\érique et $k\in \{ 0, \cdots, j-1\}$.
Pour toute $k$-face $F$ de $V$, $\pi_P(F)$ est une $k$-face de
$\pi_P(V)$ si et seulement si $P\cap C(F,V)\not =\{ 0\}$ et dans ce cas :
$$C(\pi_P(F), \pi_P(V))=P\cap C(F,V). $$
R\éciproquement \à toute $k$-face $G$ de $\pi_P(V)$ on peut associer une
unique $k$-face $F$ de $V$ telle que $\pi_P(F)=G$. On peut donc \écrire :
$$\sigma_j(V_0)= 1- \D \int_{P\in G(j,n)}\
\sum_{k=0}^{j-1} \ \ \sum_{F \in {\cal F}_k(V) }
\Theta_{j-k}\big( (C(F,V) \cap P)_0 \big)\cdot
\Theta_k(\big(\pi_P(F))_0\big)\  d \gamma_{j,n}(P).   $$
Soit alors $[F]$ l'espace vectoriel engendr\é par $F$. On d\éfinit
une application $v$ sur les polytope sph\ériques de $[F]$ de la fa\c con
suivante. Pour tout polytope
sph\érique $ W $ de $[F]$ :
$$v(W)=\D \int_{P\in G(j,n)}\ \Theta_{j-k}\big( (C(F,V) \cap P)_0 \big)\cdot
\Theta_k(\big(\pi_P(\widehat W))_0\big) \  d \gamma_{j,n}(P). $$
L'application $v$ est une valuation simple sur les polytopes sph\ériques
de $[F]$, invariante sous l'action des rotations de $[F]$. Comme de plus
$v$ est positive, par le Th\éor\ème $A.1$,~$v$ est proportionnelle
au volume $\r H^{k-1}$ sur la sph\ère unit\é de $[F]$. Il existe
$c\in \d R_+$, tel que pour tout polytope sph\érique $W$ de $[F]$,
$v(W)=c\cdot \r H^{k-1}(W )$. En \égalant $W$ \à la sph\ère unit\é de
$[F]$, on obtient, puisqu'alors $ \Theta_k(\big(\pi_P(\widehat W))_0\big)=1$ :
$$ \D \int_{P\in G(j,n)}\ \Theta_{j-k}\big( (C(F,V) \cap P)_0 \big)
  \  d \gamma_{j,n}(P)=c \cdot \r H^{k-1}(S^{k-1}_{(0,1)}). $$
On en d\éduit que :
$$v(F)= \D {\r H^{k-1}(F)\over \r H^{k-1}(S^{k-1}_{(0,1)}) }\cdot
\D \int_{P\in G(j,n)}\ \Theta_{j-k}\big( (C(F,V) \cap P)_0 \big)
  \  d \gamma_{j,n}(P)$$

$$ = \Theta_k(F)\cdot\D \int_{P\in G(j,n)}\ \Theta_{j-k}\big( (C(F,V)
\cap P)_0 \big)   \  d \gamma_{j,n}(P). $$
Maintenant l'application $u$ d\éfinie sur les polytopes sph\ériques
$W$ de $[F]^\perp$ par :
$$u(W)=\D \int_{P\in G(j,n)}\ \Theta_{j-k}\big( (\widehat W
\cap P)_0 \big)   \  d \gamma_{j,n}(P)  $$
est une valuation simple, positive, invariante sous l'action des rotations
de $[F]^\perp$. \`A nouveau, par le Th\éor\ème $A.1$, $u$ est proportionnelle
au volume sur la sph\ère unit\é de $[F]^\perp$ (on peut aussi
invoquer la formule de Cauchy-Crofton sph\érique cf [Fe2] 3.2.48),
ce qui donne :
$$v(F)= \Theta_k(F)\cdot \Theta_{j-k}\big((C(F,V))_0\big). $$
On a ainsi prouv\é l'\égalit\é :
$$\sigma_j(V_0)= 1-   \sum_{k=0}^{j-1} \ \ \sum_{F \in {\cal F}_k(V) }
\Theta_{j-k}\big( (C(F,V))_0 \big)\cdot
\Theta_k(F_0) = 1- \sum_{k=0}^{j-1} \ \ \sum_{F \in {\cal F}_k(V) }
\gamma(F,V)\cdot \Theta_k(F_0).    $$
Le Th\éor\ème $A.4$ r\ésulte de cette \égalit\é et du Lemme $A.3$. \f
\ev
Nous exprimons maintenant \à leur tour  les invariants
$\L^{\ell oc}_i$, dans le cas
poly\édral,
comme combinaisons lin\éaires
des produits des densit\és des $k$-faces par leur angle ext\érieur.

\ev
{\bf Th\éor\ème A.5. --- }
 {\sl  Quels que soient $i\in \{0, \cdots, n \}$ et le
 c\^one poly\'edral $V$ de $\d R^n$  de sommet l'origine :
$$\L_i^{\ell oc}(V_0)=\D \sum_{k=i}^n\ \ a_i^k \sum_{F\in {\cal F}_k(V)}
\gamma(F,V)\cdot \Theta_k(F_0) , $$
o\`u quel que soit $k \in \{ i, \cdots, n\}$,
$a_i^k=  \D {\alpha_k  \over \alpha_{k-i}\cdot \alpha_i}C_k^i $.}

\ev
{\bf Preuve.}
On a d\éj\`a remarqu\é que $\L_i^{\ell oc}(V_0)=\D{1\over\alpha_i }\cdot
\L_i(V\cap B^n_{(0,1)})$.
Calculons le volume du voisinage tubulaire de rayon $r$
autour de $V\cap B^n_{(0,1)}$ afin  d'obtenir $\L^{\ell oc}_i(V_0)$ :
 $$\r H^n(\r T_r(V\cap  B^n_{(0,1)}))=
\D \sum_{i=0}^{n}\alpha_i \cdot \L^{\ell oc}_i(V_0)\cdot
\alpha_{n-i} \cdot r^{n-i}. \eqno(1)$$

\n
Rappelons que $C(x,V)$ d\'esigne le c\^one conormal de $V$ en $x$.
Si  $x\in S^n_{(0,1)}\cap V$ et $F(=F_x)$ est
la face de $V$ contenant $x$,
nous d\'efinissons le c\^one conormal de $V\cap B^n_{(0,1)}$ en $x$
par $\d R_+\cdot x+  C(F,V) $ et nous le notons aussi $C(x,V)$.
Avec cette notation :
$$\r T_r(V\cap  B^n_{(0,1)})= \coprod_{x \in V\cap  B^n_{(0,1)} }
x+\big[ C(x, V) \cap B^n_{(0,r)} \big].$$
Soit $F \in {\cal F}_0(V) \cup \cdots \cup {\cal F}_n(V) $,
$x \in F \cap S^n_{(0,1)}$ et $y \in C(F,V)\cap  S^n_{(0,1)}$,
on note (pour $r\le 1$)~:
$$\r A_{x,y}= [0,1]\cdot x+ \bigg( ([0,1]\cdot y + [0,1]\cdot x)
\cap B^n_{(0,r)}) \bigg).$$
On a alors :
$$\r T_r(V\cap  B^n_{(0,1)})=
\coprod_{F \in {\cal F}_0(V) \cup \cdots \cup {\cal F}_n(V)  }\ \ \
\coprod_{(x,y) \in (F \cap S^n_{(0,1)}) \times
(C(F,V)\cap  S^n_{(0,1) })} \r A_{x,y}. \eqno(2) $$
Si $F$ est une face de $V$ de dimension $k$,
le th\éor\ème de changement de variables donne l'existence d'une
constante $c_{k,n}(r)$ ind\'ependante de $F$, telle que :
$$\r H^n\big(\coprod_{(x,y) \in (F \cap S^n_{(0,1)}) \times
(C(F,V)\cap  S^n_{(0,1) })} \r A_{x,y}\big)= c_{k,n}(r) \cdot
\gamma(F,V)\cdot \Theta_k(F_0). \eqno(3) $$
Si $F=\d R^k\times \{ 0\}^{n-k}$, comme $ \D
\r H^n\big(\coprod_{(x,y) \in (F \cap S^n_{(0,1)} )\times (
C(F,V)\cap  S^n_{(0,1) })} \r A_{x,y}\big)=
\r H^n(\r T_r(B^k_{(0,1)})) $, on   d\éduit de $(3)$ que :
$$ c_k(r)= \r H^n(\r T_r(B^k_{(0,1)}))=
\sum_{j=0}^n\L_j(B^k_{(0,1)})\cdot \alpha_{n-j}\cdot r^{n-j}
=\sum_{j=0}^k\L_j(B^k_{(0,1)})\cdot \alpha_{n-j}\cdot r^{n-j}. $$
Or les $\L_j$ ne d\épendent pas de la dimension de l'espace
euclidien dans lequel on les calcule, en voyant
$\r B^k_{(0,1)}$ dans $\d R^k$ au lieu de $\d R^n$, on obtient :
$$\alpha_k\cdot (1+r)^k=\r H^k(B_{(0,r+1)})
=\r H^k(\r T_r(B^k_{(0,1)})) =
\sum_{j=0}^k\L_j(B^k_{(0,1)})\cdot \alpha_{k-j}\cdot r^{k-j}.  $$
Cette derni\ère \égalit\é donne, pour tout $j\in \{ 0, \cdots, k\}$ :
$$\L_j(B^k_{(0,1)}) = \D {\alpha_k \over \alpha_{k-j}}\cdot C_k^j,$$
et donc :
$$ c_{k,n}(r)=\sum_{j=0}^k  \ \
 \D {\alpha_k \cdot \alpha_{n-j} \over \alpha_{k-j}}\cdot C_k^j\cdot r^{n-j} .
\eqno(4)$$
Les \égalit\és $(1), (2), (3)$ et $(4)$ donnent enfin :
$$ \r V_{V\cap B^n_{(0,1)}}(r)= \D \sum_{k=0}^n
\bigg( \sum_{j=0}^k   \D {\alpha_k \cdot \alpha_{n-j} \over
 \alpha_{k-j}}\cdot C_k^j\cdot r^{n-j}  \bigg)\ \
\sum_{F \in {\cal F}_k}\gamma(F,V)\cdot \Theta_k(F_0)
  $$
$$= \D \sum_{i=0}^n\bigg( \sum_{k=i}^n \D {\alpha_k
\over \alpha_{k-i}\cdot \alpha_i}C_k^i
\D\sum_{F\in {\cal F}_k}\gamma(F, V)\cdot \Theta_k(F_0)
\bigg)\alpha_{n-i}\cdot\alpha_i\cdot  r^{n-i}, $$
c'est-\à-dire :
$$ \L_i^{\ell oc }(V_0)=  \sum_{k=i}^n \D {\alpha_k  \over \alpha_{k-i}\cdot
\alpha_i}C_k^i
\D\sum_{F\in {\cal F}_k}\gamma(F, V)\cdot \Theta_k(F_0) .\f  $$
\ev
Des Th\éor\èmes $A.4$ et $A.5$ on d\éduit imm\édiatement les valeurs des coefficients $m_i^j$
du Th\éor\ème
$3.1$ :
\ev
{\bf Calcul des coefficients $m_i^j$ du  Th\éor\ème 3.1.  }
Soit $V$ un c\^one poly\édral de $\d R^n$ de sommet l'origine et soit
$k \in \{ 0, \cdots, n \}$. Le Th\éor\ème $A.4$ donne :
$$\sum_{F \in {\cal F}_k}\gamma(F,V)\cdot \Theta_k(F_0)=
\sigma_k(V_0)-\sigma_{k+1}(V_0), \ \ \hbox{si } k\not= n$$
$$ \ \hbox{ et }
\sum_{F \in {\cal F}_n}\gamma(F,V)\cdot \Theta_n(F_0)=
\sigma_n(V)= \Theta_n(V_0).$$
Par le Th\éor\ème $A.5$, on en d\éduit :

$$ \pmatrix{\Ll_1(V) \cr \vdots \cr \Ll_n(V)}=
\pmatrix{m_1^1 & m_1^2  & \ldots & m_1^{n-1}& m_1^n \cr
           0  & m_2^2  & \ldots & m_2^{n-1 }&m_2^n  \cr
                \vdots &        &        & & \vdots \cr
                     0 & 0      & \ldots &0& m_n^n  \cr}\cdot
         \pmatrix{\sigma_1(V)  \cr \vdots \cr \sigma_n(V)},  $$
avec :  $$ m_i^i= a_i^i=1, \ \ m_i^j= a_i^j-a_i^{j-1}=
\D \D {\alpha_j  \over \alpha_{j-i}\cdot
\alpha_i}C_j^i  - {\alpha_{j-1}  \over \alpha_{j-1-i}\cdot
\alpha_i}C_{j-1}^i , \ \ \hbox{ si }\ i+1\le j \le n. \f$$

\vfill\eject
\centerline{\bf Remerciements}
\ev
Le premier auteur remercie entre autres A. Bernig, L. Br\"oker et
N. Dutertre pour les nombreuses discussions que nous avons eues
lors des rencontres organis\'ees dans le cadre du r\'eseau europ\'een
RAAG. Cet article a b\'en\'efici\'e de s\'ejours \`a  l'Institut
Max-Planck de Bonn et \`a l'Institut Weizmann.
\Ev
\Ev
\centerline{\Tmind R\'EF\'ERENCES}
\ev

[Be-Br1] A. Bernig, L. Br\"ocker,  Lipschitz-Killing invariants. {\sl Math. Nachr.} {\bf 245} (2002), 5-25
 \vskip1mm

 [Be-Br2] A. Bernig, L. Br\"ocker, Courbures intrins\`eques dans les cat\'egories analytico-g\'eom\'etriques. {\sl Ann. Inst. Fourier}
{\bf 53} (2003), no. 6, 1897-1924
 \vskip1mm
[Bl]  W. Blaschke, Vorlesungen \"uber Integralgeometrie.
{3rd ed., VEB Deutsch. Verl. d. Wiss.,  Berlin, } (1955)
(first edition 1937)
 \vskip1mm

[Bo]  V. G. Boltianskii, Hilbert's third problem.
{John Wiley and Sons, New-York}, (1978)
\vskip1mm

[Br-Sp] J. Brian\c con, J. P. Speder, Les conditions
de Whitney impliquent
$\mu^*$ constant.
 {\it Ann. Inst. Fourier, Grenoble} {\bf 26}
(1976), 153-163
\vskip1mm

[Br-Tr] H. Brodersen, D. Trotman, Whitney
(b)-regularity is strictly weaker than Kuo's ratio test
 for real algebraic stratifications. {\it Math. Scand.} {\bf
45}, (1979), 27-34.

 \vskip1mm
[Br-Du-Ka] J. L. Brylinski, A. S.  Dubson, M. Kashiwara,
Formule de l'indice pour modules holonomes
et obstruction d'Euler locale.
{\sl  C. R. Acad. Sci. Paris S\'er. I Math.} {\bf 293} (1981),
no. 12, 573-576
\vskip1mm

 \vskip1mm
 [Br-Ku] L. Br\"ocker, M. Kuppe,
Integral geometry of tame sets.
{\sl Geom. Dedicata} {\bf 82} (2000), no. 1-3, 285-323

 \vskip1mm
 [Ch-M\"u-Sc]
J. Cheeger, W. M\"uller, R. Schrader,
Kinematic and tube formulas for piecewise
linear spaces. {\sl Indiana Univ. Math. J.}
{\bf  35} (1986)  no. 4, 737-754

\vskip1mm
[Co1] G. Comte,
Formule de Cauchy-Crofton pour la densit\'e des ensembles
sous-analytiques.
{\sl C. R. Acad. Sci. Paris, t}. {\bf 328} (1999), {\sl S\'erie I},
505-508

\vskip1mm
[Co2] G. Comte,
\'Equisingularit\'e r\'eelle : nombre de Lelong et images polaires.
{\sl Ann. Scient. \'Ec. Norm. Sup.} {\bf 33} (2000), 757-788
\vskip1mm
[Co-Li-Ro] G. Comte, J. -M. Lion, J. -Ph. Rollin,
Nature Log-analytique du volume des sous-analytiques. {\it Illinois J. Math}
{\bf 44}, 4, (2000), 884-888.

\vskip1mm
[Dra] R. N. Draper, Intersection theory in analytic
geometry. {\it Math. Ann.} {\bf 180} (1969), 175-204
\vskip1mm

[Dri] L. van den Dries,
Tame topology and o-minimal structures. {\sl London Mathematical Society
 Lecture Note Series,} {\bf 248} (1998), Cambridge University Press,
Cambridge
\vskip1mm

[Dr-Mi] L. van den Dries, C. Miller, Geometric
categories and o-minimal structures.  {\it Duke Math. J.}
{\bf 84} (1996), 497-540
\vskip1mm

[Du1] A. S. Dubson,  Classes caract\'eristiques des vari\'et\'es
singuli\`eres.
{\sl C. R. Acad. Sci. Paris S\'er. A-B} {\bf 287}
(1978), no. 4, 237-240

\vskip1mm
[Du2] A. S. Dubson,  Calcul des invariants num\'eriques des
singularit\'es et des applications.
Th\`ese, Bonn University, (1981)

\vskip1mm
[Fe1] H. Federer, The $(\Phi,k)$ rectifiable subsets
of n space. {\it Trans. Amer. Math. Soc.} {\bf 62} (1947), 114-192

\vskip1mm
 [Fe2] H. Federer,
Geometric measure theory.
{\sl Grundlehren Math. Wiss.}, Vol. {\bf 153} (1969)
{\sl Springer-Verlag}, 418-491
\vskip1mm
[Fe3] H. Federer,
Curvature measures.
{\sl Trans. Amer. Math. Soc}, Vol. {\bf 93} (1969)

\vskip1mm
[Fu1] H. G. J. Fu, Tubular neighborhoods in Euclidean spaces.
{\sl Duke Math. J.} {\bf 52} (1985), no. 4, 1025-1046

\vskip1mm
[Fu2] H. G. J. Fu, Curvature measures and generalized Morse theory.
{\sl J. Differential Geom.} {\bf 30} (1989), no. 3,
619-642

\vskip1mm
[Fu3] H. G. J. Fu, Curvature of singular spaces via the normal cycle.
Differential geometry: geometry in mathematical physics and related
topics (Los Angeles, CA, 1990), 211-221,
{\sl Proc. Sympos. Pure Math.,} {\bf 54} (1993), Part 2,
{\sl Amer. Math. Soc., Providence, RI}

\vskip1mm

[Fu4] H. G. J. Fu,   Kinematic formulas in integral
geometry. {\sl Indiana Univ. Math. J.} {\bf 39} (1990), no. 4,
1115-1154

\vskip1mm
[Fu5] H. G. J. Fu,  Curvature measures of subanalytic sets.
{\sl Amer. J. Math.} {\bf 116} (1994), no. 4, 819-880

 \vskip1mm
[Go-MacPh] M. Goresky, R. MacPherson, Stratified Morse
Theory. {\sl Ergebnisse der Mathematik und ihrer Grenzgebiete
(3)}  (1988), {\bf vol. 14},
 {\it Springer-Verlag}
\vskip1mm

 \vskip1mm
[Gr-Sc] P. M. Gruber, R. Schneider,
Problems in geometric convexity. In: Contributions to Geometry, ed.
par J. T\" olke et J. M. Wills,
Birkh\"auser Verlag, Basel, (1979), 225-278
 \vskip1mm

[Had] H. Hadwiger, Vorlesungen \"uber Inhalt, Oberfl\"ache und
Isoperimetrie. Springer-Verlag, Berlin-G\"ottingen-Heidelberg (1957)

 \vskip1mm
[Ha1] R. M. Hardt,    Stratifications of real
analytic mappings and images.  {\sl Invent. Math.} { \bf 28}  (1975),
 193-208
\vskip1mm

[Ha2]  R. Hardt, Semialgebraic local triviality in
semialgebraic mappings,
 {\sl Amer. Journal of Math.}
{\bf 102}  (1980), 291-302,

\vskip1mm

[He-Me1] J. P. Henry, M. Merle, Limites de normales,
conditions de Whitney et \'eclatement d'Hironaka.
{\it Proc. Symp. in  Pure Math. {\bf
40} (1983) (vol. 1), Arcata 1981, Amer. Math. Soc.}, 575-584
\vskip1mm

[He-Me2] J. P. Henry, M. Merle, Conditions de
r\'egularit\'e et
\'eclatements. {\it Ann. Inst. Fourier, Grenoble} {\bf 37}
(1987), 159-190
\vskip1mm

[He-Me-Sa] J. P. Henry, M. Merle, C. Sabbah, Sur la
condition de Thom stricte pour un morphisme analytique
complexe. {\it Ann. Scient. \'Ec. Norm. Sup.} {\bf 17} (1984), 227-268
 \vskip1mm

[Hi1] H. Hironaka, Normal cones in analytic Whitney
stratifications. {\it Publ. Math. I.H.E.S.}  {\bf 36} (1970),
127-138
\vskip1mm

[Hi2] H. Hironaka,  Subanalytic sets. {\sl Number theory,
algebraic geometry and commutative algebra,
in honor of Yasuo Akizuki}, Kinokuniya, Tokyo, (1973), 453-493
\vskip1mm

{[Hi3]} H. Hironaka, Stratifications and flatness. {\it Real and
             complex singularities, Oslo 1976}, Sijthoff and
             Noordhoff,
             (1977)
\vskip1mm
 [Ka1] M. Kashiwara,  Index theorem for a maximally overdetermined
system of linear differential equations.
{\sl Proc. Japan Acad.} {\bf 49} (1973), 803-804,
\vskip1mm

{[Ka2]} M. Kashiwara, $b$-functions and holonomic systems. {\it
             Invent. Math.} {\bf 38},
             (1976), 33-53
\vskip1mm

[Kl]  D. A. Klain, A short proof of Hadwiger's characterization theorem.
{\sl Mathematika } {\bf 42} (1995), 329-339
 \vskip1mm

[Kn-Pi-St]  J. Knight, A. Pillay, C. Steinhorn,
Definable sets in o-minimal strucrures II,
{\sl Trans. Amer. Math. Soc. } {\bf 295} (1986), 593-605
 \vskip1mm

[Kuo] T. C. Kuo, The ratio test for analytic Whitney
stratifications. {\sl Liverpool Singularities Symposium I,
 Lecture Notes in Math. (C. T. C. Wall, \'ed.)} { \bf 192} (1971),
 141-149
\vskip1mm

[Ku] M. Kuppe,
Integralgeometrie Whitney-Stratifizierbarer Mengen,
PhD, M\"uns\-ter-University, (1999)
\vskip1mm

[Ku-Ra] K. Kurdyka, G. Raby, Densit\'e des ensembles sous-analytiques.
{\sl Ann. Inst. Fourier} {\bf 39} (1989), no. 3, 753-771
\vskip1mm

[Ku-Po-Ra] K. Kurdyka, J. -P. Poly, G. Raby, Moyennes des fonctions
sous-analytiques, densit\'e, c\^one tangent et tranches.
(Trento, 1988), 170-177, {\sl Lecture Notes in Math.}, {\bf 1420} (1990),
Springer, Berlin
\vskip1mm

[Laf]  J. Lafontaine,
Mesures de courbure des vari\'et\'es lisses et des poly\`edres
[d'apr\`es Cheeger, M\"uller et Schrader]
{\sl S\'eminaire Boubaki, 1985-1986, Ast\'erisque}
 \vskip1mm

[La1] R. Langevin, Th. Shifrin, Polar varietes and
integral geometry, {\sl Amer. J. of Math.}, {\bf 104}
(1982), no 3, 553-605

\vskip1mm

[La2] R. Langevin,
Introduction to integral geometry.
 Col\'oquio Brasileiro de Matem\'atica.
[21st Brazilian Mathematics Colloquium],
Instituto de Matem\'a\-ti\-ca Pura e Aplicada (IMPA), Rio de Janeiro,
(1997), 160 p.

\vskip1mm
[La3]  R. Langevin,
La petite musique de la g\éom\étrie int\égrale. {\sl  La recherche de la
v\érit\é, \'Ecrit. Math., ACL-\'Ed. Kangourou} (1999), Paris,  117-143

\vskip1mm
[Le]  P. Lelong,  Int\'egration sur un ensemble
analytique complexe.  {\it  Bull. Soc. Math. France}, {\bf
85} (1957), 239-262
\vskip1mm
[L\^e-Te1] L\^e D\~ung Tr\'ang, B. Teissier,
Vari\'et\'es polaires locales et classes de Chern des
vari\'et\'es singuli\`eres. {\it Annals of Math. } {\bf 114} (1981),
  457-491
\vskip1mm

[L\^e-Te2] L\^e D\~ung Tr\'ang, B. Teissier,
Errata \`a ``Vari\'et\'es polaires locales et classes de Chern des
vari\'et\'es singuli\`eres". {\it Annals of Math. } {\bf 115} (1982),
 668-668
\vskip1mm

[L\^e-Te3] L\^e D\~ung Tr\'ang, B. Teissier, Cycles
\'evanescents et conditions de Whitney II.  {\it Proc. Symp.
in  Pure Math. {\bf 40} (1983) (vol. 2), Arcata 1981, Amer. Math.
Soc.}, 65-103
\vskip1mm

[Li]
 J. -M. Lion,  Densit\'e des ensembles semi-pfaffiens.
{\it Ann. Fac. Sci. Toulouse Math.} {\bf  6}, 7, (1998), no. 1, 87-92
\vskip1mm

{[Lio-Ro]} J. M. Lion, J. P. Rolin,  Int\'egration
des fonctions  sous-analytiques et volume des  sous-ensembles
sous-analytiques.  {\it  Ann. Inst. Fourier, Grenoble} {\bf
48},
 (1998), 755-767.
\vskip1mm

[Lo] Ta L\^e Loi, Verdier and strict Thom stratifications
in o-minimal structures. {\sl Illinois J. Math.}
{\bf  42} (1998), no. 2, 347-356
\vskip1mm

[Ma]
J. Mather, Notes on topological stability, Harvard University,
(1970)
\vskip1mm

[MacPh] R. D. MacPherson, Chern classes for singular algebraic
varieties. {\sl Ann. of Math.} (2) {\bf 100} (1974), 423-432

\vskip1mm [Me] M. Merle, Vari\'et\'es polaires, stratifications de
Whitney et classes de Chern des espaces analytiques complexes
(d'apr\ès L\^e-Teissier). {\sl S\'eminaire Bourbaki, Vol.
1982/83}, {\sl Ast\'erisque} {\bf 105-106} (1983), {\sl Soc. Math.
France}  \vskip1mm

\vskip1mm
[McMu-Sc]  P. McMullen, R. Schneider,
Valuations on convex bodies,
{Convexity and its applications, edited by Peter Gruber and
J\"org M. Wills, Boston: Birkh\"auser Verlag }
 (1983)
 \vskip1mm

[Na1]  V. Navarro Aznar, Conditions de  Whitney et
sections planes.  {\it Inv. Math.}
 {\bf 61} (1980), 199-226
\vskip1mm

[Na2] V. Navarro Aznar, Stratifications r\'eguli\`eres
et vari\'et\'es polaires locales.
{\sl Manuscrit}, (1981)
\vskip1mm

[Na-Tr] V. Navarro Aznar, D. Trotman, Whitney
regularity and generic wings.  {\it Ann. Inst. Fourier,
Grenoble}
 {\bf 31}  (1981), 87-111
 
 [Or] P. Orro, Conditions de  r\'egularit\'e, espaces
tangents et fonctions de Morse.  {\it Th\`ese, Orsay,}
 (1984).

 [Or-Tr] P. Orro, D. Trotman, On the regular
stratifications  and conormal structure of subanalytic sets.
{\it Bull. London Math. Soc.} {\bf 18},  (1986),
185-191.

[Pi-St] A. Pillay, C. Steinhorn, Definable sets in o-minimal structures
I  {\it Trans. Amer. Math. Soc}
 {\bf 295} (1986), 565-592

  \vskip1mm
 [Sah]  C. H. Sah, Hilbert's third problem: Scissors congruence,
{\sl Pitman Advanced Publishing Program, San Francisco} (1979)

\vskip1mm
[Sa] L. A. Santalo.
Integral geometry and geometric probability.
 {\sl Encyclopedia of Mathematics and its
Applications} Vol. {\bf 1}  (1976),
Addison-Wesley Publishing Co.,  London-Amsterdam

\vskip1mm
[Sc] W. Schickhoff, Whitneysche Tangentenkegel,
Multiplizit\"atsverhalten, Normal-Pseu\-do\-flach\-heit und
\"Aquisingularit\"atstheorie f\"ur Ramissche R\"aume. {\it
Schriftenreihe des Math. Inst. der Universt\"at M\"unster 2,
Serie}.  Heft {\bf 12}, (1977)

\vskip1mm

 [Sc1]  R. Schneider,
Curvatures measures of convex bodies.
{Ann. Mat. Pura appl. }  {\bf 116},  (1978),
101-134
 \vskip1mm

 [Sc2]  R. Schneider,
A uniqueness theorem for finitely additive invariant measure on a compact homogeneous space.
{Rendiconti del Circolo Matematico di Palermo, }  {\bf XXX},  (1981),
341-344
 \vskip1mm

 [Sc3]  R. Schneider,
Convex bodies: The Brunn-Minkowski Theory
{Encyclopedia of Mathematics and its Applications } {\bf 44} (1993),
{\sl Cambridge University Press}
 \vskip1mm

[Sc4]  R. Schneider,
Integral geometry - Measure theoretic approach and stochastic
applications {Advanced course on integral geometry, CRM }
 (1999)
 \vskip1mm

[Sh] M.  Shiota,
 Geometry of subanalytic and semialgebraic sets.
{\sl Progress in Mathematics} {\bf 150} (1997), Birkh\"auser Boston
\vskip1mm

[St1]  J. Steiner, Von dem Kr\"ummungsschwerpunkte ebener Curven.
{J. Reine angew. math } {\bf 21}, 33-63, (1840). Ges. Werke,
 vol 2. (1882), Reimer, Berlin, 99-159
 \vskip1mm

[St2]  J. Steiner, \"Uber parallele Fl\"achen.
{Monatsber. Preu$\beta$ Akad. Wissen., Berlin }, (1840).
Ges. Werke, vol 2. (1882), Reimer, Berlin
 \vskip1mm

[Te1] B. Teissier, Cycles \'evanescents, sections
planes et conditions de Whitney.
 {\it Ast\'erisque }{\bf 7-8} (1973), {\sl Singularit\'es \`a
Carg\`ese.}, 285-362
\vskip1mm

[Te2] B. Teissier, Vari\'et\'es polaires II:
Multiplicit\'es polaires, sections planes et conditions de
Whitney.
 {\it Actes de la conf\'erence de g\'eom\'etrie alg\'ebrique
\`a la R\`abida. Springer Lecture Notes} {\bf 961}
 (1981), 314-491

\vskip1mm
  [Th] R. Thom, Ensembles et morphismes stratifi\'es.  {\sl Bull.
Amer. Math. Soc.} {\bf 75} (1969), 240-284

\vskip1mm
[Tr1]
 D. Trotman, Counterexamples in
stratification theory : two discordant
horns,  {\sl   Real and complex singularities, Oslo} (1976),
\'ed. P. Holm, Stijthoff-Noordhoff (1977), 679-686

\vskip1mm
[Tr2]
 D. Trotman, Comparing regularity conditions on
stratifications.  {\sl   Proc.  Symp. in Pure
 Math. {\bf 40},  (vol. 2), Arcata 1981, Amer. Math. Soc.}, 
 (1983), 575-586

\vskip1mm
[Va1] G. Valette,
D\étermination et stabilit\é du type m\étrique des singularit\és.
Th\`ese Universit\'e de Provence (2003)

\vskip1mm
[Va2] G. Valette,
Volume, density and Whitney conditions.
{\sl Pr\'epublication}
\vskip1mm

[Ve] J. L. Verdier, Stratifications de Whitney et
th\'eor\`eme de Bertini-Sard.  {\sl  Inv. Math.} {\bf 36} (1976),
295-312
\vskip1mm

[We]  H. Weyl, On the volumes of tubes.
{Amer. J. Math } {\bf 61} (1939), 461-472
 \vskip1mm

{[Za1]} O. Zariski, Studies in equisingularity (I), (II), (III).  {\it
Amer. J. of Math.} {\bf 87}, {\bf 87}, {\bf 90},
 (1965), (1965), (1968), 507-536, 972-1006, 961-1023
 \vskip1mm
{[Za2]} O. Zariski, Some open questions in the theory
of singularities.
 {\it Bull. Amer. Math. Soc.} {\bf 77,} (1971),
481-491
 \vskip1mm

{[Za3]} O. Zariski, On equimultiple subvarieties of
algebroid hypersurfaces.
 {\it Proc. Nat. Acad. Sci. USA} {\bf 72,} no. 4,
(1975), 1425-1426. Correction : {\it Proc. Nat. Acad. Sci. USA}
{\bf 72,} no. 8, (1975)
 \vskip1mm

{[Za4]} O. Zariski, Foundations of a general theory of
equisingularity on r-dimensional algebroid and algebraic
varieties, of embedding dimension r+1.  {\it Amer. J. of
Math.} {\bf 101} (1979), 453-514

 \footnote{}{{\tt  \n GEORGES COMTE, PHILLIPE GRAFTIEAUX, MICHEL MERLE} :
\td  Laboratoire {\tt J. - A.} Dieudonn\'e {\tt UMR CNRS} 6621,
Universit\'e de Nice - Sophia Antipolis, 28 avenue de Valrose,
06108 {\tt NICE C}edex 2, {\tt FRANCE}.
{\tt E}-Mails : comte@math.unice.fr, merle@math.unice.fr, graftiea@math.unice.fr}

\vfill \eject
\end

{\bf Lemme.---} Pour tout c\^one convexe $X$ de $\d R^n$, on a~:
$$  \Lambda_i^{loc}(X)=\sum_{j \geq i}
{1 \over \alpha_{j-i}} {j \choose i}
\sum_{\dim F=j} \Theta_j(F)\Theta_{n-j}(C(X,F)), ,$$
ou en d'autres termes,
$$ {\cal V}_{X\cap B_n}(r)=\sum_{j=0}^n \Big(\sum_{k=0}^j
{\alpha_{n-k}\alpha_j \over \alpha_{j-k}} {j \choose k} r^{n-k}\Big)\Big(
\sum_{\dim F=j}\Theta_j(F)\Theta_{n-j}(C(X,F))\Big)$$
{\bf Preuve~:}
Calculons le volume modifi\'e de $X\cap B_n$.
Pour tout $r>0$, on a
$$ {\cal T} _r(X \cap B_n)=\coprod_{x \in X \cap B_n} x+C(X\cap B_n,x)
\cap B_n(r) =\coprod_{F}\coprod_{x \in F \cap B_n} x+C(X\cap B_n,x)
\cap B_n(r)$$
o\`u $F$ d\'ecrit l'ensemble des faces de $X$.
Or si $x$ appartient \`a l'int\'erieur $F\cap B_n^\circ$ de $F\cap B_n$, alors
$C(X\cap B_n,x)=C(X,F)$, et si $x$ appartient \`a $F\cap S_{n-1}$, alors
on a $C(X\cap B_n,x)=\d R_+ x+C(X,F)$. On en d\'eduit~:
$${\cal T}_r(X)=\coprod_{F}\Big(
F \cap B_n^\circ+C(X,F)\cap B_n(r)\coprod_{x\in F\cap S_{n-1}}
x+(\d R_+ x+C(X,F))\cap B_n(r)\Big)
$$
Or, par homog\'en\'eit\'e, on a pour toute face $F$ de $X$ de dimension $j$
$$ \vol_n\Big( F \cap B_n^\circ+C(X,F)\cap B_n(r)\coprod_{x\in F\cap S_{n-1}}
x+(\d R_+ x+C(X,F))\cap B_n(r) \Big)
={\cal V}_{r,n}(\d R^j)\Theta_j(F)\Theta_{n-j}(C(X,F))\, . $$
Or, comme ${\cal V}_{\d R^j,j}(r)=\alpha_j(1+r)^j$,
on a $\Lambda_k(\d R^j)={\alpha_j\over\alpha_{j-k}}{j\choose k}$ d'o\`u
$${\cal V}_{\d R^j,n}(r)=\sum_{k=0}^j {\alpha_{n-k}\alpha_j\over \alpha_{j-k}}{j\choose k}
r^{n-k} $$
et
$$ \Lambda_k^{loc}(\d R^j)= {1 \over \alpha_{j-k}}{j \choose k}$$
et le r\'esultat.
\ev
{\bf Remarque.}

\vfill\eject
\end

Nous aurons besoin des d\éfinitions suivantes (nous renvoyons \à [Schneider],
[Schneider-McMullen] pour plus d'informations):
\ev
{\bf D\éfinitions (valuations simples et continues). } Soit $\r K^n $
l'ensemble des convexes compacts de $\d R^n$.
Une application $v: \r K^n \to \d R$ est appel\ée une {\sl valuation} sur
$\r K^n$ lorsque:
$$v(\emptyset )=0  \hbox{ et } v(A\cup B)=v(A) + v(B) -v(A\cap B),  $$
pour tout $A, B \subset \r K^n$, tels que  $A\cup B \subset \r K^n$. Nous
dirons que la propri\ét\é pr\éc\édente est la propri\ét\é { \sl additive}
 de $v$.

La  valuation $v$ est dite { \sl simple}  lorsque $v$ s'annule sur les convexes de
dimension $<n $.

La valuation $v$ est dite  {\sl  continue} lorsque $v$ est continue
relativement \à la  m\'etrique de Hausdorff.
\ev
{\bf D\'efinition (condition ($*$)). } Soit $v$ une valuation  sur  $\r K^n$.
On dit que $v$ satisfait
la condition $(*)$ si $v$   est born\ée  uniform\'ement sur les polytopes compacts
de
$B^n_{(0,1)}  $ et
s'il existe $M>0$ tel que pour tout $\u>0$ suffisamment petit, pour
tout polytope compact $P$ et toute transformation affine orthogonale  $f$ tels
que $P$  et $f(P)$ soit inclus dans $B^n_{(0,\u)}$, on ait:
$$v(f(P))=v(P)(1+\alpha\u) , $$
avec $\alpha(=\alpha_{f,P})$  un nombre r\éel v\érifiant $|\alpha|\le M$.
\ev
 Une valuation v\érifiant la condition $(*)$ est presque invariante par les
transformations orthogonales affines. Nous allons prouver dans le lemme qui suit
que l'on peut localiser une telle valuation en une valuation invariante sous
l'action du groupe affine orthogonal.
\ev
{\bf Lemme 2.11. --- } {\sl Soit $v$ une valuation simple sur $\r K^n$,
continue
et v\érifiant la condition $(*)$. Nous d\éfinissons pour tout
entier $p$ la valuation $v_p$ par:
$$ \forall K \in  \r K^n, \ \ v_p(K)=2^{np}v(2^{-p}.K), $$
o\ù le convexe $2^{-p}.K$ est l'image de $K$ par l'homoth\'etie de centre $0$ et de
rapport $2^{-p}$. Alors:
$$\ti v= \lim_{p \to \infty} v_p $$
est une valuation simple sur  $\r K^n$, continue et invariante par l'action du
groupe affine orthogonal de $\d R^n $.}
\ev
{\bf Preuve. } Soit $\u$ suffisamment petit pour que $(*)$ ait lieu et
pour que $M\u<1$. Soit $T$ un simplexe de dimension maximale de $ B^n_{(0,\u)}$.
En coupant en deux toutes les ar\^etes de $T$, on d\éfinit $2^n$ $n$-simplexes qui
sont  tous des images de $2^{-1}.T$ par des isom\étries affines.
En appliquant $(*)$ on obtient:
$$ v(T)=2^nv(2^{-1}.T)(1+\alpha_1\u), \hbox{ avec } | \alpha_1|\le M. $$
Comme $2^{-1}.T $ est contenu dans la boule $B^n_{(0,\u/2)}$, on peut
appliquer \à nouveau le raisonnement \à $2^{-1}.T $ dans $B^n_{(0,\u/2)}$,
au lieu de $T$ dans  $B^n_{(0,\u)}$, et en it\érant ainsi
successivement dans $B^n_{(0,\u/2)}, \ldots, B^n_{(0,\u/2^{p-1})}$, on obtient:
$$ v_p(T)=v(T) \D \prod_{i=1}^{p}{1\over (1+2^{1-i}\alpha_i\u)},
\hbox{ avec } |\alpha_i |\le M, \ \forall i \in \{1, \ldots, p \} .$$
Comme le produit $\Pi =\D  \prod_{i=1}^{\infty}{1\over (1+2^{1-i}\alpha_i\u)} $
converge, $\ti v$ est bien d\éfinie sur les simplexes, et donc sur les polytopes
puisque  $v_p$ est additive et simple. De plus quel que soit le polytope $P$,
$$\ti v(P)=v(P)\cdot \Pi \ \
  \hbox{ \ et \ } $$
$$  | \ti v(P) - v_p(P)|  \le |v(P) |\cdot \Pi \cdot (1-e^{1-p}).   $$

Nous allons maintenant d\éfinir $v$ sur les convexes compacts de $B^n_{(0,\u)}$.
Soit $K\subset B^n_{(0,\u)}$ un convexe compact. Il existe une suite
$(P_k)_{k \in  \N } $ de polytopes de  $B^n_{(0,\u)}$ qui converge vers $K$
pour la m\étrique de Hausdorff.
\Ev
\Ev

\vfill\eject
\end

Soit maintenant $P$ un polytope et $f$ une transformation affine
orthogonale. Montrons que $\vt(P)=\vt(f(P))$.
Pour tout polytope $P$, choisissons $k$ tel que $2^{-k}P$
et  $2^{-k}f(P)$
soit inclu dans la boule de centre $0$ et de rayon $\epsilon$.
Pour tout $p > k$, $(*)$ s'applique \`a $2^{-p}P$ et $2^{-p}f(P)$,
d'o\`u l'existence de $\alpha_p$, $|\alpha_p| < M$, tel que
$$ v_p(P)=2^{np} v(2^{-p}P) = 2^{np} v(2^{-p}f(P))
  (1 + 2^{(p-k)}\epsilon \alpha_p)
= v_p(f(P))  (1 + 2^{-n(p-k)}\epsilon \alpha_p) \, .$$
Le  passage \`a la limite donne l'\'egalit\'ee attendue $\vt(P)=\vt(f(P))$.

La continuit\'e de $\vt$ par rapport \`a la distance
de Hausdorf d\'ecoule de celle de $v_p$ ainsi que de
la convergence uniforme de $v_p$ vers $\vt$
sur les simplexes contenus dans un compact fix\'e.

\vfill\eject
\end